
\input colordvi




\def\updated{3 November 2006}

\count100= 92
\count101= 49

\message{: version 20sep06}

%





\font\titlefont=cmr17
\font\titlei=cmmi10 at 17pt
\font\titlesy=cmsy10 at 17pt
\font\titleit=cmti10 at 17pt
\font\titlesl=cmsl10 at 17pt
\font\titlebf=cmbx10 at 17pt
\font\Bbbt=msbm10 at 17pt
\font\twelverm=cmr12
\font\twelvebf=cmbx10 at 12pt
\font\twelvei=cmmi10 at 12pt
\font\twelvesy=cmsy10 at 12pt
\font\ninerm=cmr9
\font\sevenrm=cmr7
\font\sixrm=cmr6
\font\fiverm=cmr5
\font\ninei=cmmi9
\font\seveni=cmmi7
\font\ninesy=cmsy9
\font\sevensy=cmsy9

\font\sixbf=cmbx6
\font\fivebf=cmbx5
\font\ninebf=cmbx9
\font\nineit=cmti9
\font\ninesl=cmsl9
\font\nineex=cmex9
\font\nineBbb=msbm9

\font\dfont=cmss10
\font\dfont=cmbx10

\font\Bbb=msbm10
\font\Bbbs=msbm7

\def\ninepoint{\def\rm{\fam0\ninerm}
    \textfont0 = \ninerm
    \textfont1 = \ninei
    \textfont2 = \ninesy
    \textfont3 = \nineex
    \scriptfont0 = \sevenrm
    \scriptfont1 = \seveni
    \scriptfont2 = \sevensy
    \scriptscriptfont0 = \fiverm
    \scriptscriptfont1 = \fivei
    \scriptscriptfont2 = \fivesy
    \textfont\itfam=\nineit \def\it{\fam\itfam\nineit}
    \textfont\bffam=\ninebf \scriptfont\bffam=\sixbf
    \scriptscriptfont\bffam=\fivebf \def\bf{\fam\bffam\ninebf}
    \textfont\slfam=\ninesl \def\sl{\fam\slfam\ninesl}
    \let\Bbb\nineBbb
    \baselineskip 10pt}

\def\titlepoint{\def\rm{\fam0\titlefont}
    \textfont0 = \titlefont
    \textfont1 = \titlei
    \textfont2 = \titlesy
    \textfont3 = \nineex
    \scriptfont0 = \twelverm
    \scriptfont1 = \twelvei
    \scriptfont2 = \twelvesy
    \textfont\itfam=\titleit \def\it{\fam\itfam\titleit}
    \textfont\bffam=\titlebf \scriptfont\bffam=\twelvebf
    \def\bf{\fam\bffam\titlebf}
    \textfont\slfam=\titlesl \def\sl{\fam\slfam\titlesl}
    \let\Bbb\Bbbt\titlefont}



\let\CC\Cbb
\def\CCt{\hbox{\Bbbt C}}  
\def\CCs{\hbox{\Bbbs C}}  

\let\RR\Rbb





\def\makebold#1{\mathord{\setbox0=\hbox{#1}%
       \copy0\kern-\wd0%
       \raise\dimen1\copy0\kern-\wd0%
       {\advance\dimen1 by \dimen1\raise\dimen1\copy0}\kern-\wd0%
       \kern\dimen0\raise\dimen1\copy0\kern-\wd0%
       {\advance\dimen1 by \dimen1\raise\dimen1\copy0}\kern-\wd0%
       \kern\dimen0\raise\dimen1\copy0\kern-\wd0%
       {\advance\dimen1 by \dimen1\raise\dimen1\copy0}\kern-\wd0%
       \kern\dimen0\raise\dimen1\copy0\kern-\wd0%
       \kern\dimen0\box0}}




\def \dword#1{{\dfont #1}}

%


\def\inpro#1{\langle#1\rangle}
\def\norm#1{\Vert#1\Vert}


\def\frac#1#2{{#1 \over #2}}




\def\endproofsymbol{\makeblanksquare6{.4}}
\def\eop{\endproofsymbol\nopf}


\def\nopf{\medskip\goodbreak}

\def\makeblanksquare#1#2{
\dimen0=#1pt\advance\dimen0 by -#2pt
      \vrule height#1pt width#2pt depth0pt\kern-#2pt
      \vrule height#1pt width#1pt depth-\dimen0 \kern-#1pt
      \vrule height#2pt width#1pt depth0pt \kern-#2pt
      \vrule height#1pt width#2pt depth0pt
}

\magnification\magstephalf

\hsize6.5truein\vsize8.6truein
\voffset.5truein


\def\title#1{\toneormore#1||||:}
\def\titexp#1#2{\hbox{{\titlefont #1} \kern-.25em%
  \raise .90ex \hbox{\twelverm #2}}\/}
\def\titsub#1#2{\hbox{{\titlefont #1} \kern-.25em%
  \lower .60ex \hbox{\twelverm #2}}\/}

\def\author#1{\bigskip\bigskip\aoneormore#1||||:\smallskip\centerline{\updated}}

\def\abstract#1{\bigskip\bigskip\medskip%
 {\ninepoint
 \narrower{\bf Abstract.~}\rm#1\smallskip
  MSC: \mscnumbers\ifx\keywords\empty\else\smallskip
  Keywords: \keywords\fi\bigskip
 \printtochere}\starttoc\bigskip}

\def\toneormore#1||#2||#3:{\centerline{\titlepoint #1}%
    \def\next{#2}\ifx\next\empty\else\medskip\toneormore#2||#3:\fi}
\def\aoneormore#1||#2||#3:{\centerline{\twelverm #1}%
    \def\next{#2}\ifx\next\empty\else\smallskip\aoneormore#2||#3:\fi}

\newwrite\toc\def\tocone{0}\def\tochalf{.5}\def\toctwo{1}
\countdef\counter=255
\def\diamondleaders{\global\advance\counter by 1
  \ifodd\counter \kern-10pt \fi
  \leaders\hbox to 15pt{\ifodd\counter \kern13pt \else\kern3pt \fi
    .\hss}}
\newdimen\lextent
\def\printtochere{\immediate\closeout\toc\relax%
\begingroup
\def\\##1.  ##2.  {\setbox1=\hbox{##1}\ifnum\wd1>\lextent\lextent\wd1\fi}
\lextent0pt\inputifthere{\jobname.toc}\advance\lextent by 2em\relax
\def\\##1.  ##2.  {\centerline{\hbox to \lextent{\rm##1\def\next{##2}%
\ifx\next\empty\else\diamondleaders\fi\hfil\hbox to 2em{\hss##2}}}}
\inputifthere{\jobname.toc}\endgroup}

\def\starttoc{\immediate\openout\toc=\jobname.toc}
\def\nexttoc#1{{\let\folio=0\edef\next{\write\toc{#1}}\next}}

\def\tocline#1#2#3{\nexttoc{\noexpand\noexpand\noexpand\\\hskip#2truecm #1.  #3.  }}

\def\footnoterule{\kern -3pt \hrule width 0truein \kern 2.6pt}
\def\leftheadline{\ifnum\pageno=\count100 \hfill%
  \else\hfil\it\shortauthor\hfil\llap{\rm\folio}\fi}
\def\rightheadline{\ifnum\pageno=\count100 \hfill%
  \else\hfil\it\shorttitle\hfil\llap{\rm\folio}\fi}

\nopagenumbers
\headline{\ifodd\pageno\rightheadline \else\leftheadline\fi}
\footline{\hfil}
\null\vskip 18pt
\centerline{}
\pageno=\count100
\count102=\count100
\advance\count102 by -1
\advance\count102 by \count101


\def\copyright{\hbox{{\twelverm o}\kern-.61em\raise .46ex\hbox{\fiverm c}}}

\insert\footins{\sixrm
\medskip
\baselineskip 8pt
\leftline{Surveys in Approximation Theory
  \hfill {\ninerm \the\pageno}}
\leftline{Volume 2, 2006.
pp.~\the\pageno--\the\count102.}
\leftline{Copyright \copyright\ 2006 Surveys in Approximation Theory.}
\leftline{ISSN 1555-578X}
\leftline{All rights of reproduction in any form reserved.}
\smallskip
\par\allowbreak}


\def\sect#1{\startsect\edef\showsectno{\the\sectionno}\let\tocindent\tocone%
       \soneormore#1||||:\relax\medskip\noindent\ignorespaces}

\def\soneormore#1||#2||#3:{%
   \leftline{\bf\showsectno\hskip2em #1}
   \def\next{#2}%
   \ifx\next\empty\puttocline{\showsectno\ \ #1}{\folio}%
   \else\puttocline{\showsectno\ \ #1}{}\let\showsectno\skipsectno\let\tocindent\tochalf\soneormore#2||#3:\fi}

\def\puttocline#1#2{\tocline{#1}{\tocindent}{#2}}
\def\skipsectno{\setbox0=\hbox{\the\sectionno}\hskip\wd0}

\def\subsect#1{\formal{#1.}\let\tocindent\toctwo\puttocline{#1}{\folio}}
\def\formal#1{\bigskip{\bf #1}\hskip1em}


\newcount\sectionno\sectionno0
\def\presect{\the\sectionno.}
\newcount\subsectionno
\def\startsect{\ifx\empty\presect\else\restartnums\fi%
               \subsectionno0\global\advance\sectionno by 1\relax

               \goodbreak\bigskip\smallskip}
\def\startsubsect{\global\advance\subsectionno by 1\goodbreak\bigskip}


\def\figinbox#1(#2,#3)#4#5{\centerline{\vbox{\gridbox#2/#3/{
\ifshowfigname\point(0,0){#1}\fi 
  \point(0,0){\epsfxsize=#4truecm \epsfbox{\figsource #1}}#5}}}}

\input epsf
\def\figsource{}
\newif\ifshowfigname

\def\foneormore#1||#2||#3:{\centerline{\ninepoint\rm #1}%
    \def\next{#2}\ifx\next\empty\else\vskip0pt\foneormore#2||#3:\fi}


\def\gridbox#1/#2/#3{
\vbox to #1\gridunits{#3
\ifshowgrid\tickcount=0
  \loop\cgridw%
   \vbox to 0pt{\kern\tickcount \gridunits\hrule width#2\gridunits
       height\gridwidth\vss}
   \nointerlineskip \advance\tickcount by \tickskip
   \ifdim\tickcount pt<#1pt\repeat 
  \hbox to 0pt{\tickcount=0\hbox to 0pt{\tick#1/\hss}\advance\tickcount by \tickskip%
 \loop\ifdim\tickcount pt<#2pt\nexttick\tickskip#1/\advance\tickcount by \tickskip \repeat\hss}
\else \vbox to 0pt{\hrule width#2\gridunits height0pt\vss}
\fi\vfil}\vfil}

\def\ppoint#1#2(#3,#4)#5{\setbox0=\hbox{#5}
   \dimen0=\ht0\advance\dimen0 by\dp0\divide\dimen0 by-2
   \multiply\dimen0 by#1\advance\dimen0 by#3\gridunits
   \dimen1=\wd0\divide\dimen1 by-2\multiply\dimen1 by#2
   \advance\dimen1 by#4\gridunits\dpoint(\dimen0,\dimen1){#5}}

\def\gridunits{truecm}\newcount\tickskip\tickskip1\newcount\majortick\majortick5

\def\point(#1,#2)#3{\dpoint(#1\gridunits,#2\gridunits){#3}}
\def\dpoint(#1,#2)#3{\vbox to 0pt{\kern#1
   \hbox{\kern#2{#3}}\vss}\nointerlineskip}
\newcount\rmndr   
\def\rem#1#2{\rmndr=#1{}\divide\rmndr by#2{}%
\multiply\rmndr by-#2{}\advance\rmndr by #1}
\def\cgridw{\gridwidth\finegridw{}\rem\tickcount\majortick%
   \ifnum\rmndr=0{}\gridwidth\roughgridw\fi} 

\def\tick#1/{\cgridw\vrule width\gridwidth height0pt depth#1\gridunits}
\def\nexttick#1#2/{\hbox to#1\gridunits{\hfil\tick#2/}}
\newcount\tickcount
\newdimen\finegridw\finegridw0.4pt\newdimen\roughgridw\roughgridw1.6pt
\newdimen\gridwidth
\newif\ifshowgrid \showgridtrue





\def\label#1{%
  \ifsamelabel\global\samelabelfalse\else
  \ifmmode\global\advance\eqnum by 1
  \else\global\advance\labelnum by 1
  \fi\fi
  \edef\griff{label:#1}\edef\inhalt{\lastlabel}\definieres%
  \ifmmode\eqno(\inhalt)\else\inhalt\fi
  \ifdraft\ifmmode\rlap{\fiverm #1}\else\marginal{#1}\fi\fi}

\def\eqalignlabel#1{{\def\eqno{}\let\labelnum\eqnum\label{#1}}}

\def\labelplus#1#2{\def\labelsub{#2}\relax\label{#1}\def\labelsub{}}

\def\eqalignlabelplus#1#2{{\def\eqno{}\let\labelnum\eqnum\labelplus{#1}{#2}}}



\newif\ifsamelabel
\def\labelsub{}
\def\lastlabel{\presect\ifmmode\the\eqnum\else\the\labelnum\fi\labelsub}
\def\nextlabel{{\ifmmode\advance\eqnum by 1\else\advance\labelnum by 1\fi\lastlabel}}


\newcount\blackmarks\blackmarks0
\newcount\eqnum
\newcount\labelnum
\def\restartnums{\eqnum0\labelnum0}
\def\singlecount{\let\labelnum\eqnum}

 \newread\testfl
 \def\inputifthere#1{\immediate\openin\testfl=#1
    \ifeof\testfl\message{(#1 does not yet exist)}
    \else\input#1\fi\closein\testfl}

 \inputifthere{\jobname.aux}
 \newwrite\aux
 \immediate\openout\aux=\jobname.aux

\def\plazieres{\expandafter\ifx\csname\griff\endcsname\relax%
  \xdef\esfehlt{\griff}\blackmark\else{\csname\griff\endcsname}\fi}

\def\definieres{\expandafter\xdef\csname\griff\endcsname{\inhalt}%
 \def\blankkk{ }\expandafter\immediate\write\aux{%
 \string\expandafter\def\string\csname%
 \blankkk\griff\string\endcsname{\inhalt}}}

\def\blackmark{\ifnum\blackmarks=0\global\blackmarks=1%
 \write16{============================================================}%
 \write16{Some forward reference is not yet defined. Re-TeX this file!}%
 \write16{============================================================}%
 \fi\immediate\write16{undefined forward reference: \esfehlt}%
 {\vrule height10pt width2pt depth2pt}\esfehlt%
 {\vrule height10pt width2pt depth2pt}}

\def\marginal#1{\strut\vadjust{\kern-\strutdepth%
\vtop to \strutdepth{\baselineskip\strutdepth\vss\llap{\fiverm#1\ }\null}}}
\def\strutdepth{\dp\strutbox}


\newif\ifdraft

\newcount\hour\newcount\minutes
\def\draft{\drafttrue
\headline={\sevenrm \hfill\ifx\filenamed\undefined\jobname\else\filenamed\fi%
(.tex) (as of \ifx\updated\undefined???\else\updated\fi)
 \TeX'ed at {\hour\time\divide\hour by 60{}%
\minutes\hour\multiply\minutes by 60{}%
\advance\time by -\minutes
\the\hour:\ifnum\time<10{}0\fi\the\time\  on \today\hfill}}
}

\def\today{\number\day\space%
\ifcase\month\or January\or February\or March\or April\or May\or June\or
 July\or August\or September\or October\or November\or December\fi%
\space\number\year}



\def\References{\goodbreak\bigskip\centerline{\bf References}%
   \tocline{\skipsectno\ \  References}{\tocone}{\folio}%
   \bigskip\frenchspacing}

\def\bibitem{\smallskip\noindent}


\gdef\formfirstauthor{\the\firstname\  \the\lastname}
\gdef\formnextauthor{, \the\firstname\the\lastname}
\gdef\formotherauthor{ and \the\firstname\the\lastname}
\gdef\formlastauthor{,\formotherauthor}

\gdef\formB{\the\au\ [\yr] ``\the\ti'', \the\pb, \the\pl. \setcitelabel}
\gdef\formD{\the\au\ [\yr] ``\the\ti'', dissertation, \the\pl. \setcitelabel}
\gdef\formJ{\the\au\ [\yr] \the\ti, {\sl\the\jr}\ifx\vl\empty%
\else\ {\bf\vl}\fi, \pp. \setcitelabel}
\gdef\formP{\the\au\ [\yr] \the\ti, in {\sl\the\tit},
\getfirstchar\aut\ifx\firstchar\unknownx\else\the\aut, ed\edsop, \fi
\getfirstchar\pub\ifx\firstchar\unknownx\else\the\pub, \fi \the\pl, \pp. \setcitelabel}
\gdef\formR{\the\au\ [\yr] \the\ti\ifx\is\empty\else, \is\fi. \setcitelabel}

\newtoks\lastname
\newtoks\firstname
\newtoks\au
\newtoks\aut
\newtoks\ti
\newtoks\tit
\newtoks\pb
\newtoks\pub
\newtoks\pl
\newtoks\jr

\newtoks\rhlau

\def\setcitelabel{\edef\griff{cit\rh}\edef\inhalt{\the\rhlau\ \yr}\definieres}
\def\setcitelabel{}

\def\getfirstchar#1{\edef\theword{\the#1}\expandafter\getit\theword:}
\def\getit#1#2:{\def\firstchar{#1}}
\def\unknownx{x}

\newif\ifonesofar
\def\concat#1{\edef\audef{{#1}}\au=\audef}
\def\decodeauthor#1, #2,#3;{\lastname={#1}\firstname={#2}%
\concat{\formfirstauthor}\onesofartrue%
\def\morerhlau{}%
\def\next{#3}\ifx\next\empty\else\def\morerhlau{ et al.}\decodemoreauthor#3;\fi
\edef\morerhlauu{{\the\lastname\morerhlau}}\rhlau=\morerhlauu}
\def\decodemoreauthor#1, #2,#3;{\lastname={#1}\firstname={#2}%
\def\next{#3}\ifx\next\empty\let\formaut=\formlastauthor%
\ifonesofar\ifx\formotherauthor\undefined\else\let\formaut=\formotherauthor%
\fi\fi\concat{\the\au\formaut}%
\else\onesofarfalse\concat{\the\au\formnextauthor}\decodemoreauthor#3;\fi}

\def\refproc #1(#2; #3; {\decodeproc#2; \xdef\yr{#3}}
\def\decodeproc#1), #2 (ed#3.),#4 (#5); {%
 \global\tit={#1}\global\aut={#2}\xdef\edsop{#3}\global
 \pub={#4}\global\pl={#5}}


\def\refB #1; #2; #3 (#4); #5; {\decodeauthor#1,;%
   \ti={#2}\pb={#3}\pl={#4}\def\yr{#5}\bibitem\formB}

\def\refD #1; #2; #3; #4; {\decodeauthor#1,;%
   \ti={#2}\pl={#3}\def\yr{#4}\bibitem\formD}

\def\refJ #1; #2; #3; #4; #5; #6; {\decodeauthor#1,;%
    \ti={#2}\jr={#3}\def\vl{#4}\def\yr{#5}\def\pp{#6}\bibitem\formJ}

\def\refP #1; #2; #3; #4; {{\global\aut={\vrule height15pt width15pt depth0pt}%
 \global\tit={{\bf the specified proceedings does not exist in our files}}%
 \xdef\edsop{}\global\pub={}\def#3{}\input proceed }\decodeauthor#1,;%
        \ti={#2}\def\pp{#4}\bibitem\formP}

\def\refQ #1; #2; (#3; #4; #5; {\decodeproc#3; \decodeauthor#1,;%
   \ti={#2}\def\yr{#4}\def\pp{#5}\bibitem\formP}

\def\refR #1; #2; #3; #4; {\decodeauthor#1,;%
         \ti={#2}\def\is{#3}\def\yr{#4}\bibitem\formR}

\def\presect{} 


\title{Approximation in $\CCt^N$}
\author{Norm Levenberg}

\def\shorttitle{Approximation in $\CC^N$}
\def\shortauthor{Norm Levenberg}

\def\mscnumbers{32-02, 41-02}

\def\keywords{}


\def\sump{\mathop{{\sum}'}} 


\abstract{This is a survey article on selected topics in approximation theory. The topics either use techniques from the theory of several complex variables or arise in the study of the subject. The survey is aimed at readers having an acquaintance with standard results in classical approximation theory and complex analysis but no apriori knowledge of several complex variables is assumed. }

\sect{Introduction and motivation} Let $\CC^N=\{(z_1,\ldots,z_N):z_j\in \CC\}$ where $z_j=x_j+iy_j$ and identify $\RR^N=\{(x_1,\ldots,x_N):x_j \in \RR\}$.  A complex-valued function $f$ defined on an open subset of $\CC^N$ is \dword{holomorphic} if it is separately holomorphic in the appropriate planar region as a function of one complex variable when each of the remaining $N-1$ variables are fixed. This deceptively simple-minded criterion is equivalent to any other standard definition; e.g., $f$ is locally representable by a convergent power series in the complex coordinates; or $f$ is of class $C^1$ and satisfies the Cauchy-Riemann system
$${\partial f\over \partial \bar z_j}:={1\over 2}\bigl ({\partial f\over \partial x_j}+i{\partial f\over \partial y_j}\bigr )=0, \qquad \ j=1,\ldots,N.$$ In particular, holomorphic functions are smooth, indeed, real-analytic; whereas the separately holomorphic criterion makes no apriori assumption on continuity (Hartogs separate analyticity theorem, circa 1906; cf., [Sh] section 6). We make no assumptions nor demands on the reader's knowledge of several complex variables (SCV) but we do require basic knowledge of classical one complex variable (CCV)  theory. An acquaintance with potential theory in CCV, i.e., the study of subharmonic functions, would be helpful in motivating analogies with pluripotential theory, the study of plurisubharmonic functions in SCV, but it is not essential. Sections 2 and 3 provide some background on the important notions of polynomial hulls and plurisubharmonic functions in SCV. Section 4 recalls some classical approximation theory results from CCV. In addition, two short appendices are included (sections 13 and 14) for those interested in a brief discussion of a few specialized topics in SCV: pluripolar sets, extremal plurisubharmonic functions, and the complex Monge-Amp\`ere operator. We highly recommend the texts by
\item {(1)} Ransford [Ra] on potential theory in the complex plane;
\item {(2)} Klimek [K] on pluripotential theory; and
\item {(3)} Shabat [Sha] on several complex variables.

\noindent H\"ormander's SCV text [H\"o] is a classic. Range's book [Ran] is an excellent source for integral formulas in SCV; these will occur at several places in our discussion (cf., sections 3, 7 and 10). Many of the approximation topics we mention are described in the monograph of Alexander and Wermer [AW].

Zeros of holomorphic functions locally look like zero sets of holomorphic polynomials (Weierstrass Preparation Theorem; e.g., [Sha] section 23). In particular, in $\CC^N$ for $N>1$ these sets are never isolated. Consider, for example, $f(z_1,\ldots,z_N)=z_1$: the zero set is a copy of $\CC^{N-1}\subset \CC^N$. This means that, apriori, Runge-type pole-pushing arguments do not exist in SCV. Henceforth the term ``polynomial'' will refer to a \dword{holomorphic polynomial}, i.e., a polynomial in $z_1,\ldots,z_N$, unless otherwise noted. We use the notation ${\cal P}_d={\cal P}_d(\CC^N)$ for the polynomials of degree at most $d$.

Continuing on this theme, rational functions, i.e., ratios of polynomials,  behave quite differently in SCV than in CCV. Consider, in $\CC^2$, the function $r(z_1,z_2):=z_1/z_2$. The ``zero-set'' of $f$ contains the punctured plane $\{z_1=0\}\setminus {(0,0)}$ and the ``pole-set'' contains the punctured plane $\{z_2=0\}\setminus {(0,0)}$, but the point $(0,0)$ itself forms the ``indeterminacy locus'': $f$ is not only undefined at this point, but, as is easily seen by simply considering complex lines $z_2 = tz_1$ through $(0,0)$, $f$ attains all complex values in any arbitrarily small neighborhood of this point.

It is still the case that polynomials are the nicest examples of holomorphic functions and rational functions are the nicest examples of meromorphic functions (which we won't define) in SCV. Thus one  wants to utilize these classes in approximation problems. Many standard tools from CCV either don't exist in SCV or are often more complicated.

In this introductory section, we first recall some classical approximation-theoretic results in the plane with an eye towards generalization, if possible, to $\CC^N, \ N>1$.
Let $K$ be a compact subset of $\CC^N$, and let $C(K)$ denote the uniform algebra of continuous, complex-valued functions endowed with the supremum (uniform) norm on $K$. Let $P(K)$ be the uniform algebra (subalgebra of $C(K)$) consisting of uniform limits of polynomials restricted to $K$. Finally, let $R(K)$ be the uniform closure in $C(K)$ of rational functions $r=p/q$ where $q(z)\not = 0$ for $z\in K$.

As a sample, a question which has a complete and common answer in CCV and SCV, to be given in sections 3 and 4, is:
{\sl For which compact sets $K\subset \CC^N$ is it true that for {\bf any} function $f$ that is holomorphic in a neighborhood of $K$ there exists a sequence $\{p_n\}$ of polynomials which converges uniformly to $f$ on $K$; i.e., $f|_K\in P(K)$? Moreover, for such compacta, estimate $d_n(f,K):= \inf \{\norm{f-p}_K: \deg p \leq n\}$ in terms of the ``size'' of the neighborhood in which $f$ is holomorphic.}

For example, if $N=1$ and $K=\bar \Delta:=\{z:|z|\leq 1\}$ is the closed unit disk, writing $f(z)=\sum_{k=0}^{\infty} a_kz^k$ as a Taylor series about the origin, the Taylor polynomials $p_n(z)=\sum_{k=0}^n a_kz^k$ converge uniformly to $f$ on $K$. More precisely, if $f$ is holomorphic in the disk $\Delta(0,R):=\{z:|z|< R\}$ of radius $R>1$, the Cauchy estimates give
$$|a_k|=\left|{1\over 2\pi i}\int_{|z|=\rho}{f(z)\over z^{k+1}}dz \right| \leq {\sup_{|z|\leq \rho} |f(z)|\over \rho^{k}}\eqno(1)$$
for any $1<\rho <R$ yielding
$$d_n(f,\bar \Delta)\leq \norm{f-p_n}_{\bar \Delta}\leq {1\over (1-1/\rho)}{\sup_{|z|\leq \rho} |f(z)|\over \rho^{n+1}}\eqno(2)$$
so that $\limsup_{n\to \infty} d_n(f,\bar \Delta)^{1/n}\leq 1/R$. On the other hand, taking $K=T:=\partial \Delta :=\{z:|z|= 1\}$ the unit circle, the function $f(z)=1/z$ is holomorphic in $\CC^*=\CC\setminus \{0\}$ but if $p(z)$ is a polynomial with
$|f(z)-p(z)|<\epsilon <1$ on $T$, then, multiplying by $z$, we have $|1-zp(z)|<\epsilon <1$ on $T$ and hence, by the maximum modulus principle, on $\bar \Delta$. This gives a contradiction at $z=0$.

The difference in these sets is explained, and a continuation of our review of classical complex approximation theory proceeds, if we recall a version of the Runge theorem for $N=1$:

\proclaim Theorem (Ru). \ Let $K\subset \CC$ be compact with $\CC\setminus K$ connected. Then for any function $f$ holomorphic on a neighborhood of $K$, there exists a sequence
 $\{p_n\}$ of holomorphic polynomials which converges uniformly to $f$ on $K$.

The condition ``$\CC\setminus K$ connected'' is equivalent, when $N=1$, to $K=\hat K$ where
$$\hat K:=\{z\in \CC^N: |p(z)|\leq \norm{p}_K \ \hbox{for all holomorphic polynomials} \ p\}$$
is the \dword{polynomial hull} of $K$. Clearly a uniform limit on $K$ of a sequence of polynomials yields a holomorphic function on the interior $K^o$ of $K$; this observation motivates one of the conditions in Lavrentiev's result:

\proclaim Theorem (La). \ Let $K\subset \CC$ be compact with $\CC\setminus K$ connected. Then $P(K)=C(K)$ if and only if $K^o=\emptyset$.

In any number of (complex) dimensions, the maximal ideal space of the uniform algebra $C(K)$ is $K$ and that of $P(K)$ is $\hat K$. Thus  a necessary condition that $P(K)=C(K)$ is that $K=\hat K$. Lavrentiev's theorem shows that in the complex plane, removing the only other obvious obstruction yields a necessary and sufficient condition for the density of the polynomials in the space of continuous functions. A nice exposition of these results (and more) in a succinct, clear manner is given in Alexander-Wermer [AW], section 2.  The techniques utilized are elementary functional analysis (Hahn-Banach), classical potential theory (logarithmic potentials) and classical complex analysis (Cauchy transforms).

If we allow $K$ to have interior, then we may ask if functions in $C(K)$ which are holomorphic on $K^o$ are uniformly approximable on $K$ by polynomials. This is the content of Mergelyan's theorem:

\proclaim Theorem (Me). \ Let $K\subset \CC$ be compact with $\CC\setminus K$ connected. Then for any function $f\in C(K)$ which is holomorphic on $K^o$, there exists a sequence $\{p_n\}$
of polynomials which converges uniformly to $f$ on $K$.

What happens in $\CC^N$ for $N>1$? The complex structure plays a major role. As an elementary, but illustrative, example, consider two disks $K_1$ and $K_2$ in $\CC^2 =\{(z_1,z_2):z_1,z_2\in \CC\}$ defined as follows:
$$K_1:=\{(x_1,x_2)\in \RR^2: x_1^2+x_2^2 \leq 1\} \ \hbox{and}$$
$$K_2:=\{(z_1,0): |z_1|\leq 1\}.$$
Both of these sets are ``polynomially convex'' in $\CC^2$; i.e., $\hat K_1=K_1$ and $\hat K_2=K_2$; thus each set satisfies the obvious necessary condition for holomorphic polynomials to be dense in the space of continuous functions on the set. However, $K_2$ lies in the complex $z_1$-plane and $P(K_2)$ can be identified with $P(K)$ where $K$ is the closed unit disk in {\it one} complex variable; the observation made regarding Lavrentiev's theorem shows that $P(K_2)\not = C(K_2)$.

To understand $K_1$ and to motivate an attempt to generalize Lavrentiev's theorem in SCV, we first recall the classical theorem of Stone-Weierstrass:

\proclaim Theorem (SW). \ Let ${\cal U}$ be a subalgebra of $C(K)$ containing the constant functions and separating points of $K$. If $f\in {\cal U}$ implies that $\bar f \in {\cal U}$, then ${\cal U}=C(K)$.

As an immediate corollary, we have the {\it real}  Stone-Weierstrass theorem (which includes the classical Weierstrass theorem for a real interval):

\proclaim Theorem (RSW). \ Let $K$ be a compact subset of\/ $\RR^N\subset \CC^N$. Then $P(K)=C(K)$.

Thus by (RSW), $P(K_1)=C(K_1)$. The difference here is that {\sl the real submanifold $\RR^2=\RR^2+i0$ of $\CC^2$ is totally real; i.e., $\RR^2$ contains no complex tangents}.
We will generalize this example in Theorem (HW) of section 8. The extremely difficult question of determining when $P(K)=C(K)$ will be partially analyzed in the next section.

Recall that $R(K)$ is the uniform subalgebra of $C(K)$ generated by rational functions which are holomorphic on $K$. The Hartogs-Rosenthal theorem gives a sufficient condition for $R(K)=C(K)$ if $K\subset \CC$.

\proclaim Theorem (HR1). \ Let $K$ be a compact subset of $\CC$ with two-dimensional Lebesgue measure zero. Then $R(K)=C(K)$.

\noindent A similar result holds in $\CC^N, \ N>1$. For $\alpha >0$, we let $h_{\alpha}$ denote $\alpha$-Hausdorff measure.

\proclaim Theorem (HRN). \ Let $K$ be a compact subset of $\CC^N$ with $h_2(K)=0$. Then $R(K)=C(K)$.

\noindent This follows since the conjugates $\bar z_j$ of the coordinate functions belong to $R(K)$, by Theorem (HR1); from this it follows trivially that $R(K)$ is closed under complex conjugation. Then Theorem (SW) implies the conclusion.

We turn to a $\CC^N$-version of Theorem (Ru). Note that if we take the ``boundary circles'' of our sets $K_1$ and $K_2$, i.e., take
$$X_1:=\{(x_1,x_2)\in \RR^2: x_1^2+x_2^2 =1\} \ \hbox{and}$$
$$X_2:=\{(z_1,0): |z_1|= 1\},$$
then a higher-dimensional version of Theorem (Ru) is valid for $X_1$ but not for $X_2$, i.e., if $f$ is holomorphic on a neighborhood of $X_1$ (in $\CC^2$!), then there exists a sequence $\{p_n\}$ of  holomorphic polynomials which converges uniformly to $f$ on $X_1$ (e.g., $f|_{X_1}\in C(X_1)$ and $P(X_1)=C(X_1)$ follows from Theorem (RSW)); the analogous statement is not true for $X_2$ (why?). Here the difference can simply be explained by the fact that $X_1$ is polynomially convex while $X_2$ is not (indeed, $\hat X_2 =K_2$). This is the content of the Oka-Weil theorem:

\proclaim Theorem (OW). \ Let $K\subset \CC^N$ be compact with $\hat K=K$. Then for any function $f$ holomorphic on a neighborhood of $K$, there exists a sequence $\{p_n\}$ of polynomials which converges uniformly to $f$ on $K$.

This result was first proved by Andr\'e Weil in 1935 by using a multivariate generalization of the Cauchy integral formula for certain polynomial polyhedra. We sketch his argument in section 3.  In 1936 Kyoshi Oka gave a different proof that made use of his celebrated ``lifting principle'' (cf.\ [AW] Chapter 7).

As a motivational example for what to expect, let $K=\{(z_1,\ldots,z_N):|z_j|\leq 1, \ j=1,\ldots,N\}$ be the closed unit polydisk. If $f$ is holomorphic in a larger polydisk
$D_R:= \{(z_1,\ldots,z_N):|z_j|<R, \ j=1,\ldots,N\}, \ R>1$, then iterating the one-variable Cauchy integral formula, for $\rho <R$ we obtain the formula
$$f(z)=({1\over 2 \pi i})^N\int_{|\zeta_1|=\rho} \cdots \int_{|\zeta_N|=\rho} {f(\zeta_1,\ldots,\zeta_N)\over (\zeta_1-z_1)\cdots (\zeta_N -z_N)}d\zeta_1\cdots d\zeta_N\eqno(3)$$
valid for $z\in D_{\rho}$. We  can write a Taylor series expansion $f(z)=\sum_{|\alpha|=0}^{\infty}a_{\alpha}z^{\alpha}$ where $\alpha =(\alpha_1,\ldots,\alpha_N)$ is a multiindex with $|\alpha|:=\sum_{j=1}^N\alpha_j$ and $z^{\alpha}:=z_1^{\alpha_1}\cdots z_N^{\alpha_N}$ and
$$a_{\alpha}=({1\over 2 \pi i})^N\int_{|z_1|=\rho} \cdots \int_{|z_N|=\rho} {f(z_1,\ldots,z_N)\over z_1^{\alpha_1+1} \cdots z_N^{\alpha_N+1}}dz_1\cdots dz_N.$$
The same estimates as in (1) and (2) show that not only is $f|_K\in P(K)$ but we obtain, using the Taylor polynomials $p_n(z)=\sum_{|\alpha|=0}^n a_{\alpha}z^{\alpha}$, the quantitative estimate $\limsup_{n\to \infty} d_n(f,K)^{1/n}\leq 1/R$. Note that in the Cauchy integral formula (3), the integration takes place over the $N$-dimensional torus $\{(z_1,\ldots,z_N): |z_j|=\rho, \ j=1,\ldots,N\}$, which is a proper subset of the $(2N-1)$-dimensional topological boundary $\partial D_{\rho}$ if $N>1$.

\sect{Polynomial hulls and polynomial convexity}

The condition that $K=\hat K$ occurs in Theorems (Ru), (La), (Me) and (OW); indeed, this condition is implicit in Theorem (RSW): any compact subset of $\RR^N$ is polynomially convex (exercise!). If $K\subset \CC$, $\hat K$ is the union of $K$ with the bounded components of $\CC \setminus K$. For $K\subset \CC^N$ if $N>1$, $\hat K$ contains the union of $K$ with the bounded components of $\CC^N \setminus K$ but it can be much, much more. An elementary example is the connected and simply connected set $K:=K_1\cup K_2$ which is the union of two bidisks
$$K_1:=\{(z_1,z_2): |z_1|\leq 1, \ |z_2|\leq r <1\}$$
and
$$K_2:=\{(z_1,z_2): |z_2|\leq 1, \ |z_1|\leq r <1\}.$$
We'll see in the next section that
$$\hat K = \{(z_1,z_2): |z_1|\leq 1, \ |z_2|\leq 1, \ |z_1z_2|\leq r\}$$
(draw a picture in $|z_1|,|z_2|$-space). Note it is clear that $\hat K$ is contained in the right-hand-side by considering the polynomial $p(z_1,z_2)=z_1z_2$. In general, the polynomial hull of a compact set is difficult to describe.

It follows readily from the maximum modulus principle that if we have a bounded holomorphic mapping $f=(f_1,\ldots,f_N):\Delta \to \CC^N$; i.e., each $f_j:\Delta \to \CC$ is a bounded holomorphic function, with (componentwise) radial limit values $f^*(e^{i\theta}):=\big(f_1^*(e^{i\theta}),\ldots,f_N^*(e^{i\theta})\big)\in K$ for almost all $\theta$, then $f(\Delta) \subset \hat K$. In general, we will say that a set $S\subset \CC^N$ has \dword{analytic structure} if it contains a nonconstant analytic disk $f(\Delta)$. Thus one way to obtain (lots of) points in $\hat K$ is the existence of analytic structure in $\hat K$. Moreover, existence of analytic structure in a compact set $S$ precludes the possibility of $C(S)=P(S)$ for the set.

In 1963, Stolzenberg [Sto] gave an example of a compact set $K$ in the topological boundary
$$\partial (\Delta \times \Delta)=\{(z_1,z_2): |z_1|=1, \ |z_2|\leq 1, \ \hbox{or} \  |z_2|=1, \ |z_1|\leq 1\}$$
of the bidisk $\Delta \times \Delta \subset \CC^2$ such that the origin $(0,0)\in \hat K$ but the projections $\pi_{z_1}(\hat K),\pi_{z_2}(\hat K)$ of $\hat K$ in each coordinate plane contain no nonempty open set; thus $\hat K$ contains no analytic structure. From the lack of analytic structure in $\hat K$ one may be tempted to conjecture that $P(\hat K)=C(\hat K)$. However, there clearly exist $f\in C(\hat K)$ with $|f(0,0)| > \norm{f}_K$ for this set $K$; e.g., $f(z_1,z_2)=1-\max [|z_1|,|z_2|]$. For any $p\in P(\hat K)$, we obviously have $\norm{p}_{\hat K}=\norm{p}_K$. Thus $f\not\in P(\hat K)$.

How can one tell if $\hat K \setminus K$ contains analytic structure? Note that an analytic disk has locally finite Hausdorff two-measure.

\proclaim Theorem (Alexander-Sibony). \ Let $K$ be a compact subset of $\CC^N$ and let $q\in \hat K \setminus K$. If there exists a neighborhood $U$ of $q$ with $h_2(\hat K\cap U)<+\infty$, then $\hat K\cap U$ is a one-dimensional analytic subvariety of $U$.

\noindent This means that $\hat K\cap U$ is essentially a one-dimensional complex manifold (modulo some singular points) and hence looks locally like a nonconstant analytic disk $f(\Delta)$. In particular, if $\hat K\setminus K \not=\emptyset$ and $\hat K \setminus K$ contains no analytic structure, then $h_2(\hat K \setminus K)=+\infty$. A nice discussion of the Alexander-Sibony result can be found in section 21 of [AW]. In [DL] the authors constructed examples of compact sets $K\subset \CC^N$ whose polynomial hull $\hat K$ contains no analytic structure but such that $\hat K \setminus K$ has positive $2N$-Hausdorff measure.

Recall that a compact subset $K$ of $\RR^N$ is automatically polynomially convex and, moreover, $P(K)=C(K)$ for such sets. From Theorem (HRN) and Theorem (OW), we also get the following result.

\proclaim Corollary. \ Let $K=\hat K\subset \CC^N$ with $h_2(K)=0$. Then $P(K)=C(K)$.

\vskip6pt

\noindent {\it Question:} {\sl For an arbitrary compact set
$K\subset \CC^N$, find a ``nice'' condition (C)  on $K$ so that if
$K=\hat K$, then $K$ satisfies (C) if and only if $P(K)=C(K)$;
i.e., find a $\CC^N$-version of Theorem (La)}. \vskip6pt

We return to this matter in section 8. As a final note to reinforce the delicate nature of polynomial hulls in $\CC^N$ for $N>1$,
we mention the curious results of E. Kallin [Ka]. The union of any {\it two} disjoint convex compact sets in $\CC^N$ is polynomially convex. The union of any {\it three} disjoint closed Euclidean {\it balls} is polynomially convex. On the other hand, there exist {\it three} disjoint {\it convex} sets whose union is {\it not} polynomially convex. Moreover, it is unknown whether the union of {\it four or more} disjoint closed {\it balls} is polynomially convex, unless, e.g., the centers of the balls lie on the real subpace $\RR^N$ of $\CC^N$ [Kh].

\sect{Plurisubharmonic functions and the Oka-Weil theorem} We outline the basic notions and sketch a proof of the Oka-Weil theorem, Theorem (OW). First of all, in the complex plane, any domain $D\subset \CC$ is a \dword{domain of holomorphy}; i.e., there exists $f$ holomorphic in $D$ -- we write $f\in{\cal O}(D)$ -- which does not extend holomorphically across any boundary point of $D$. This follows from the classical Weierstrass theorem which allows the construction of a nontrivial holomorphic function with prescribed discrete zero set in $D$. However, in $\CC^N, \ N>1$, there exist domains $D$ with the property that {\it every} $f\in {\cal O}(D)$ extends holomorphically to a larger domain $\tilde D$ (independent of $f$). Products of planar domains, e.g., polydisks, are obviously domains of holomorphy. A simple example of a domain $D\subset \CC^2$ which is not a domain of holomorphy is
$$D=\{(z_1,z_2): |z_1|<2, \ |z_2|<2\}\setminus \{(z_1,z_2): 1\leq |z_1|<2, \ |z_2|\leq 1\}.$$
It is straightforward to see (use Laurent series!) that any $f\in {\cal O}(D)$ extends holomorphically to the bidisk $\{(z_1,z_2): |z_1|<2, \ |z_2|<2\}$. Another example is the interior of the set $K=K_1\cup K_2$ from the previous section: any $f$ holomorphic on
$$D:=  \{(z_1,z_2): |z_1|< 1, \ |z_2|< r\}\cup  \{(z_1,z_2): |z_1|< r, \ |z_2|< 1\}$$
extends holomorphically to the domain
$$\tilde D := \{(z_1,z_2): |z_1|< 1, \ |z_2|< 1, \ |z_1z_2|< r\}.$$
{\bf Exercise}: {\sl For each $f\in {\cal O}(D)$, if $\tilde f$ denotes the holomorphic extension of $f$ to $\tilde D$, then $f(D)=\tilde f(\tilde D)$.}

As a sample of holomorphic extendability situations in SCV, consider the following three results on bounded domains $D\subset \CC^N$. The first two can be reduced to classical one-variable arguments if one utilizes the Cauchy integral formula for polydisks (3).
\vskip4pt
\item {1.} ({\bf Morera type}): {\sl Let $N\geq 1$ and let $S$ be a smooth, real hypersurface in $D$; i.e., (locally) $S=\{\rho =0\}$ where $\rho$ is a smooth, real-valued function on a neighborhood of $S$ and $d\rho \not = 0$ on $S$. If $f\in {\cal O}(D\setminus S)\cap C(D)$, then $f\in {\cal O}(D)$}.
\item {2.} ({\bf Riemann removable singularity type}): {\sl Let $N\geq 1$ and let $A$ be a complex analytic hypersurface in $D$; i.e., (locally)  $A=\{g =0\}$ where $g$ is holomorphic on a neighborhood of $A$. If $f\in {\cal O}(D\setminus A)$ is locally bounded on $A$ (i.e., for each $z\in A$ there is a neighborhood $U$ of $z$ with $f$ bounded on $(D\setminus A)\cap U$), then $f$ has a holomorphic extension $F\in {\cal O}(D)$}.
\item {3.} ({\bf Hartogs type}) {\sl Let $N> 1$ and let $A$ be a complex analytic subvariety of (complex) codimension two in $D$; i.e., $A$ is (locally) the common zero set of two holomorphic functions. If $f\in {\cal O}(D\setminus A)$, then $f$ has a holomorphic extension $F\in {\cal O}(D)$}.
\vskip4pt

We refer the reader to section 32 of [Sha]. As an example of 3., any function holomorphic in a punctured ball
$$D=\{(z_1,z_2)\in \CC^2: 0< |z_1-a_1|^2+|z_2-a_2|^2 < R^2\}$$
in $\CC^2$ extends holomorphically across the puncture $a=(a_1,a_2)$. Note that no  boundedness assumptions on $f$ are required. Indeed, a theorem of Hartogs states that {\sl if $K$ is a compact subset of a domain $D\subset \CC^N$ ($N>1$) such that $D\setminus K$ is connected, then every $f\in {\cal O}(D\setminus K)$ extends holomorphically to $D$}. Thus $D\setminus K$ is not a domain of holomorphy.

A real-valued function $u:D\to [-\infty,+\infty)$ defined on a domain $D\subset \CC^N$ is called \dword{plurisubharmonic} (\dword{psh}) if $u$ is uppersemicontinuous (usc) on $D$ and $u|_{D\cap L}$ is subharmonic on (components of) $D\cap L$ for any complex affine line $L=L_{z_0,a}:=\{z_0+ta:t\in \CC\}$ ($z_0,a\in \CC^N$ fixed).
The canonical examples of such functions are those of the form $u=\log {|f|}$ where $f\in {\cal O}(D)$. The class of psh functions on a domain $D$, denoted $PSH(D)$, forms a convex cone; i.e., if $u,v\in PSH(D)$ and $\alpha,\beta \geq 0$, then $\alpha u+\beta v\in PSH(D)$. The limit function $u(z):=\lim_{n\to \infty}u_n(z)$ of a decreasing sequence $\{u_n\}\subset PSH(D)$ is psh in $D$ (we may have $u\equiv -\infty$); while for any family $\{v_{\alpha}\}\subset PSH(D)$ (resp., sequence $\{v_n\}\subset PSH(D)$) which is uniformly bounded  above on any compact subset of $D$, the functions
$$v(z):=\sup_{\alpha} v_{\alpha}(z) \ \hbox{and} \ w(z):=\limsup_{n\to \infty} v_{n}(z)$$
are ``nearly'' psh: the usc regularizations
$$v^*(z):=\limsup_{\zeta \to z} v(\zeta) \ \hbox{and} \ w^*(z):=\limsup_{\zeta \to z} w(\zeta)$$
are psh in $D$. Finally, if $\phi$ is a real-valued, convex increasing function of a real variable, and $u$ is psh in $D$, then so is $\phi \circ u$.

If $u\in C^2(D)$, then $u$ is psh if and only if for each $z\in D$ and vector $a\in \CC^N$, the Laplacian of $t\mapsto u(z+ta)$ is nonnegative at $t=0$; i.e., {\sl the complex Hessian $[{\partial^2u\over \partial z_j \partial \bar z_k}(z)]$ of $u$ is positive semidefinite on $D$}:
$$\sum_{j,k=1}^N{\partial^2u\over \partial z_j \partial \bar z_k}(z)a_j\bar a_k\geq 0.$$
In particular, the trace of the complex Hessian is nonnegative so that $u$  is $\RR^{2N}-$subharmonic. Indeed, {\sl $u:D\to [-\infty,+\infty)$ is psh if and only if $u\circ A$ is $\RR^{2N}-$subharmonic in $A^{-1}(D)$ for every complex linear isomorphism $A$}; moreover, the notion of a psh function makes sense on any complex manifold. If $u\in C^2(D)$, we call $u$ \dword{strictly} psh if the complex Hessian of $u$ is positive definite; i.e., we have strict inequality in the above displayed equation provided $a\not = (0,\ldots,0)$. The function $u(z)=|z|^2:=|z_1|^2 +\cdots +|z_N|^2$ is strictly psh on $\CC^N$.

A domain $D\subset \CC^N$ is said to be (globally) \dword{pseudoconvex} if $D$ admits a \dword{psh exhaustion function}:  there exists $u$ psh in $D$ with the property that the sublevel sets $D_c:=\{z\in D:u(z) <c\}$ are compactly contained in $D$ for all real $c$. Any planar domain $D\subset \CC$ is pseudoconvex: if $D=\CC$, take $u(z)=|z|$; if $D=\CC\setminus \{z_0\}$, take $u(z)=1/|z-z_0|$; otherwise, take
$$u(z)=|z|+\sup_{z_0\in \partial D} \frac{1}{|z-z_0|} = |z|+{\rm dist}(z,\partial D)^{-1}.$$
It turns out that {\sl a domain $D\subset \CC^N$ is pseudoconvex if and only if the function $$u(z)= -\log {\rm dist}(z,\partial D)$$ is psh in $D$}. Note that, in this case,
$$\exp(u(z))={\rm dist}(z,\partial D)^{-1}$$
is psh since $x\to e^x$ is convex and increasing; if, e.g., $D$ is bounded, both $u(z)$ and $\exp(u(z))$ are continuous psh exhaustion functions for $D$. A slightly more restrictive notion is that of \dword{hyperconvexity}: $D$ is hyperconvex if it admits a negative psh exhaustion function $u$; i.e., $u<0$ in $D$ and the sets $D_c$ are compactly contained in $D$ for all $c<0$. Such a domain is pseudoconvex, for the function $-1/u:=\phi\circ u$, where $\phi(x) = -1/x$ is convex and increasing for $x<0$, is then a psh exhaustion function for $D$. In particular, a fact we will need below is that {\sl if $D$ is pseudoconvex and $u\in PSH(D)$, the components of the sublevel sets $\{z\in D:u(z) <c\}$ are pseudoconvex domains}.

The condition that a domain $D\subset \RR^N$ with smooth (say $C^2$) boundary be convex can be described analytically by the existence of a smooth defining function $r$:
$$D =\{x\in \RR^N: r(x) < 0\}; \quad \partial D =\{x\in \RR^N: r(x) = 0\};  \quad dr\not = 0 \ \hbox{on} \ \partial D,$$
such that the real Hessian of $r$ is positive semidefinite on the real tangent space $T_p(\partial D)$ to $\partial D$; i.e., for $p\in \partial D$,
$$\sum_{j,k=1}^N{\partial^2r\over \partial  x_j  \partial x_k}(p)a_j a_k\geq 0 \quad \hbox{if} \quad \sum_{j=1}^N{\partial r\over  \partial x_j }(p)a_j=0.$$
Strict convexity means that the real Hessian is positive definite on the real tangent space to $\partial D$. The complex analogue of convexity  is \dword{Levi-pseudoconvexity}: a smoothly bounded domain $D\subset \CC^N$ is Levi-pseudoconvex if it admits a defining function $r$ whose {\it complex} Hessian is positive semidefinite on the {\it complex} tangent space $T_p^{\CCs}(\partial D)$ to $\partial D$; i.e., for $p\in \partial D$,
$$\sum_{j,k=1}^N{\partial^2r\over \partial  z_j  \partial \bar z_k}(p)a_j\bar  a_k\geq 0 \quad \hbox{if} \quad \sum_{j=1}^N{\partial r\over  \partial z_j }(p)a_j=0.$$
We say that $\partial D$ is \dword{Levi-pseudoconvex at $p\in \partial D$} if this holds. Note that $T_p^{\CCs}(\partial D)$ is an $(N-1)$-{\it complex}-dimensional linear subspace of the $(2N-1)$-{\it real}-dimensional space $T_p(\partial D)$. \dword{Strict Levi-pseudoconvexity} means that the complex Hessian is positive definite on the complex tangent space to $\partial D$. It turns out that {\sl if $D$ is strictly Levi-pseudoconvex, one can make a holomorphic change of coordinates so that $D$ is (locally) strictly convex}. Precisely, if $p\in \partial D$ is a strictly Levi-pseudoconvex boundary point, then there is a neighborhood $U\subset \CC^N$ of $p$ and a biholomorphic map $\phi: U\to \phi(U)\subset \CC^N$ such that $\phi(U\cap \partial D)$ is strictly convex (in the $\RR^{2N}$-sense) at $\phi(p)$. It is {\it not} the case that if $D$ is merely Levi-pseudoconvex at $p$, then one can make a holomorphic change of coordinates so that in the new coordinates $D$ is convex.

There is a relationship between pseudoconvexity and Levi-pseudoconvex: if $D$ is pseudoconvex, then there exists an increasing  sequence of bounded, strictly Levi-pseudoconvex domains $D_n \subset D$ with smooth boundary which are relatively compact in $D$ such that
$$D = \bigcup_{n=1}^\infty D_n.$$
This follows since once we have a psh exhaustion function $u$ for $D$, we can modify it to get a smooth, strictly psh exhaustion function $\tilde u$; then we can take
$$D_n:=\{z\in D:\tilde u(z)< m_n\}$$
for an appropriate sequence $\{m_n\}$ with $m_n\uparrow \infty$.

We already observed that in $\CC$, every domain is both a domain of holomorphy and a pseudoconvex domain. In $\CC^N$ for $N>1$ this no longer holds, but these two notions are equivalent: a domain $D\subset \CC^N$ is a  domain of holomorphy if and only if it is pseudoconvex. The ``if'' direction is deep and is one version of the so-called \dword{Levi problem}. Thus from now on we will simply use the terminology  ``pseudoconvex domain.'' Products of pseudoconvex domains are pseudoconvex domains. Euclidean balls and, more generally, convex domains (in the $\RR^{2N}$-sense) are pseudoconvex; but there are many non-convex pseudoconvex domains.

The statement that canonical examples of psh functions are those of the form $u=\log {|f|}$ where $f\in {\cal O}(D)$ can now be made precise:

\proclaim Theorem (B). \ A psh function $u$ on a pseudoconvex
domain can be written in the form
$$u(z)=\bigl[ \limsup_{j\to \infty} a_j\log {|f_j(z)|}\bigr]^*$$
where $a_j\geq 0$ and $f_j\in {\cal O}(D)$.

\noindent Here $u^*(z)$ is the usc regularization of the function $u$. The result is false if $D$ is not pseudoconvex [Le].

\noindent {\bf Sketch of proof}. The domain
$$\tilde D := \{(z,w)\in D\times \CC\subset \CC^{N+1}:z\in D, \ |w|<e^{-u(z)}\}$$ is pseudoconvex, for the function $\tilde u(z,w):=u(z)+\log {|w|}$ is psh in $D\times \CC$ and
$$\tilde D =\{(z,w)\in D\times \CC:\tilde u(z,w) <0\}.$$
Since $\tilde D$ is a pseudoconvex domain and hence a domain of holomorphy, there exists $F\in {\cal O}(\tilde D)$ which is not holomorphically extendible across any boundary point of $\tilde D$; expanding $F$ in a \dword{Hartogs series}, i.e., a series of the form
$$F(z,w)=\sum_{j=0}^{\infty}f_j(z)w^j$$
where $f_j\in {\cal O}(D)$, the radius of convergence in $w$ (as a function of $z$) is given by the usual formula
$$R(z)=\bigl [ \limsup_{j\to \infty} |f_j(z)|^{1/j}\bigr]^{-1}.$$
Since $F$ is not holomorphically extendible across any boundary point of $\tilde D$, $R(z)=e^{-u(z)}$ almost everywhere; i.e., $u(z) =-\log {R(z)}$ a.e.; taking usc regularizations (making the right-hand-side psh) the result follows. \eop

\noindent This result is due to Bremermann but the proof given is due to Sibony. An important remark, which we will use later, is that a slight refinement of this argument yields a local uniform approximability result (cf.\ [JP], Proposition 4.4.13): {\sl if $u$ is psh and continuous in a bounded domain $D$, then for any compact set $K\subset D$ and $\epsilon >0$ there exist finitely many $f_j\in {\cal O}(D), \ j=1,\ldots,m$ and positive constants $a_1,\ldots,a_m$ such that}
$$\bigl |u(z)-\max_{j=1,\ldots,m}a_j\log {|f_j(z)|} \bigr | < \epsilon  \quad \hbox{for}  \quad z\in K.\eqno(4)$$

Given a polynomially convex compact set $K\subset \CC^N$; i.e., $K=\hat K$, and given a bounded open neighborhood $U$ of $K$, by compactness of $K$ and the definition of $\hat K$ we can find finitely many polynomials $q_1,\ldots,q_m$ such that
$$K\subset \Pi :=\{z\in U:|q_j(z)|< 1, \ j=1,\ldots,m\} \subset U.\eqno(5)$$
Note that necessarily $m\geq N$. We call $\Pi$ a \dword{polynomial polyhedron}. By slightly modifying the polynomials $q_j$, we may assume that $\Pi$ is a \dword{Weil polyhedron} which simply means that the ``faces''
$$\sigma_j:=\{z\in U: |q_j(z)|=1, \ |q_k(z)|\leq 1, \ k\not = j\}$$
are real $(2N-1)$-dimensional manifolds and the intersection of any $s$ distinct faces for $2\leq s\leq N$ has dimension at most $2N-s$. The set of $N$ dimensional ``edges''  $\sigma_{i_1,\ldots,i_N}:=\sigma_{i_1}\cap \cdots \cap \sigma_{i_N}$ form the ``skeleton'' of $\Pi$. For the unit polydisk $P:=\{(z_1,\ldots,z_N):|z_j|< 1, \ j=1,\ldots,N\}$, the skeleton, or distinguished boundary, is the $N$-torus $T^N:=\{(z_1,\ldots,z_N):|z_j|= 1, \ j=1,\ldots,N\}$. We have the following generalization of the Cauchy integral formula for polydisks (cf.\ [Sha] section 30).

\proclaim Theorem (Weil integral formula). \  Let $\Pi$ be a Weil polyhedron. Then for any $f\in {\cal O}(\Pi)\cap C( \bar \Pi)$,
$$f(z)=({1\over 2 \pi i})^N \sump_{i_1,\ldots,i_N}\int_{\sigma_{i_1,\ldots,i_N}} {f(\zeta)\det \{[P_m^{i_n}(\zeta,z)]_{m,n=1,\ldots,N}\} \over \prod_{\nu =1}^N(q_{i_{\nu}}(\zeta)-q_{i_{\nu}}(z))}d\zeta_1\cdots d\zeta_N$$
for $z\in \Pi$.

\noindent Here, $\sum_{i_1,\ldots,i_N}'$ refers to a summation over increasing multiindices: $i_1< \cdots < i_N$. The functions $P_s^t$ are polynomials satisfying
$$q_i(\zeta)-q_i(z)=\sum_{j=1}^N(\zeta_j-z_j)P_j^i(\zeta,z).$$
This formula is a special case of a more general result known as Hefer's formula; here, we are merely rearranging the Taylor expansion of $q_i$ at the point $z$. As an example, for the polydisk $P$, one can take $m=N$ and $q_j(z)=z_j$ in which case $P_s^t=\delta_{st}$.

Returning to the setting of the Oka-Weil theorem, Theorem (OW), given a function $f$ holomorphic on a neighborhood $U$ of $K=\hat K$, we construct a Weil polyhedron satisfying (5). We  obtain a Taylor-like expansion
$$f(z)=\sum_{|\kappa|=0}^{\infty}{\sum_I}'A^I_\kappa(z)q_I(z)^\kappa$$
where $\kappa=(k_1,\ldots,k_N), \ I=(i_1,\ldots,i_N), \ q_I(z)^\kappa=q_{i_1}(z)^{k_1}\cdots q_{i_N}(z)^{k_N}$, and
$$A^I_\kappa(z)=({1\over 2 \pi i})^N \int_{\sigma_{i_1,\ldots,i_N}} {f(\zeta)\over q_{i_1}(\zeta)^{k_1+1}\cdots q_{i_N}(\zeta)^{k_N+1}}\times $$
$$\det \{P_m^{i_n}(\zeta,z)_{m,n=1,\ldots,N}    \}d\zeta_1\cdots d\zeta_N$$
are polynomials. Using truncations of this expansion, as with the polydisk at the end of section 1, one concludes that $f|_K\in P(K)$, completing the outline of the proof of Theorem (OW). An alternate proof of Theorem (OW), using Lagrange interpolation at generalized Fekete points (see section 9) has been given by Siciak [Si].

From (4), (5) and Theorem (OW), it follows that if $K$ is a compact set in $\CC^N$ and $D$ is an open neighborhood of the polynomial hull $\hat K$, then $\hat K$ can just as well be constructed as a ``hull'' with respect to holomorphic or continuous psh functions; i.e., $\hat K$ coincides with both
$$\hat K_{{\cal O}(D)}:=\{z: |f(z)|\leq \norm{f}_K \ \hbox{for all} \ f\in {\cal O}(D)\}$$
and
$$\hat K_{PSH(D)}:=\{z: u(z)\leq \sup_{\zeta \in K}u(\zeta) \ \hbox{for all} \ u\in PSH(D)\cap C(D)\}.$$
The reader may now verify the claim in the previous section about the polynomial hull of $K_1 \cup K_2$ utilizing the above observation with $D$ being a dilation of the domain
$$\tilde D := \{(z_1,z_2): |z_1|< 1, \ |z_2|< 1, \ |z_1z_2|< r\}$$
and the exercise at the beginning of this section.

\sect{Quantitative approximation theorems in $\CC$} Before jumping to a
quantitative Runge-type theorem in $\CC^N$, we recall the example of the closed unit disk $\bar \Delta$ to review the one-variable story. In the introduction we proved one direction of the following.
 \proclaim Theorem. \
Let $f$ be continuous on $\bar \Delta = \{ z \in \CC : |z| \le 1
\} $,
and $R > 1$.
Then
$$
  \limsup_{n \to \infty} d_n(f,\bar \Delta)^{1/n}\leq 1/R
\eqno(6)
$$
if and only if $f$ is the restriction to $\bar \Delta$ of a function
holomorphic
in $\Delta(0,R)= \{ z \in \CC : |z| < R \}$.

\noindent {\bf Proof}. For the only if direction we note that for any nonconstant polynomial $p$, the function $$u(z):={1\over \deg p}\log {|p(z)|\over \norm{p}_{\bar \Delta}}-\log {|z|}$$
is subharmonic on $\CC\setminus \bar \Delta$, bounded at $\infty$, and nonpositive on $T=\partial \Delta$. By the maximum principle $u\leq 0$ on $\CC\cup \{\infty\} \setminus \bar \Delta$ which gives
$$|p(z)| \le  \norm{p} _{\bar \Delta}|z|^{\deg p}, \ |z|\geq 1,$$
hence
$$|p(z)| \le  \norm{p} _{\bar \Delta}\max (|z|,1)^{\deg p}=\norm{p} _{\bar \Delta}(e^{\log^+{|z|}})^{\deg p},\quad z\in \CC.$$
This is the \dword{Bernstein-Walsh inequality}. In particular,
$$
  |p(z)| \le  \norm{p} _{\bar \Delta} \ \rho ^{\deg p},
  \qquad
  |z| \le \rho.
\eqno(7)
$$
Let $f$ be a continuous function on $\bar \Delta$ such that
(6) holds and choose a polynomial $p_n$ of degree at most $n$ satisfying
$d_n=\norm{f-p_n}_{\bar \Delta}$. We claim that the series
$p_0 + \sum_1^{\infty} (p_n-p_{n-1})$ converges
uniformly on compact subsets of $\{ z : |z| < R \}$ to a holomorphic
function $F$ which agrees with $f$ on
$\bar \Delta$. For if
$1< R' < R$, by hypothesis the polynomials
$p_n$ satisfy
$$
  \norm{f-p_n}_{\bar \Delta} \leq {M \over  {R'} ^n}, \qquad n=0,1,2,\ldots,
\eqno(8)
$$
for some $M > 0$.  Then for $1 < \rho <  R'$, simply apply
(7) to $p_n - p_{n-1}$:
$$\eqalign{
  \sup _{|z| \le \rho} |p_n(z) - p_{n-1} (z)|
  \leq &\;
  \rho^n \norm{p_n-p_{n-1}}_{\bar \Delta}\cr
  \leq &\;
  \rho^n ( \norm{p_n -  f}_{\bar \Delta} + \norm{f - p_{n-1}}_{\bar
\Delta} )
  \leq
  \rho^n  {M ( 1 + R' ) \over {R'} ^n}.\cr}
$$
From (8), $F=f$ on $\bar \Delta$.
 \eop

We generalize this. Let  $K$ be a compact subset of $\CC$ with $\CC\setminus K$ connected. Recall this is equivalent to the condition that $K=\hat K$; i.e., $K$ is polynomially convex. We say that $K$ is \dword{regular} if there is a continuous
function $g_K : \CC \to [0,+\infty)$ which is identically equal
to zero on $K$, harmonic on $\CC \setminus K$, and has a logarithmic
singularity at infinity in the sense that
$g_K(z) - \log |z|$ is harmonic at infinity (this is equivalent to $\CC \setminus K$ being a regular domain for the Dirichlet problem). We call $g_K$ the classical \dword{Green function} for $K$. For $K=\bar \Delta$, we have $g_{\bar \Delta}(z)=\log^+{|z|}$. In the general case, if $p$ is any nonconstant polynomial, then the function
$$V := {1 \over \deg p} \log {|p|\over \norm{p}_K} - g_K$$ is subharmonic on
$\CC \setminus K$, bounded at $\infty$,
and continuously assumes
nonpositive values on $\partial K$.  By the maximum principle we
have $V \le 0$ on $\CC \cup \{ \infty \} \setminus K$, yielding the Bernstein-Walsh property
$$|p(z)| \leq \norm{p}_K (e^{g_K(z)})^{ \deg p }.$$
In particular, if $R > 1$ and
$$
  D_R := \{ z : g_K (z) < \log R \},
\eqno(9)
$$
then
$$
  |p(z)| \leq \norm{p}_K R^{ \deg p },
  \qquad   z \in D_R.
\eqno (10)
$$
Then a similar argument proves one-half of the following univariate Bernstein-Walsh theorem.

\proclaim Theorem (BW1). Let $K$ be a
regular compact subset of the plane with Green function $g_K$.  Let $R > 1$, and define $D_R$
by {\rm (9)}.
Let $f$ be continuous on $K$. Then
$
  \limsup_{n\to \infty} d_n(f,K)^{1/n}\leq 1/R
$
if and only if $f$ is the restriction to $K$ of a function
holomorphic in $D_R$.

\noindent {\bf Proof}. We prove the other half using
duality. We suppose $f$ is holomorphic on  $D_R$; and
we will rewrite the numbers $d_n$ in such a way that we can
estimate them. Let $1 < r < \rho < R$. To get a global $C^{\infty}$ extension $F$ of $f$ that agrees with $f$ on
a neighborhood of $K$, we let $\phi$ be a smooth
cut-off function which is identically equal to $1$ on $\bar D_{\rho}$
and has compact support in $D_R$.  We
then set $F=\phi f$ in $D_R$ and let $F$ be identically $0$ outside of
$D_R$.

For $n$ fixed, by the Hahn-Banach theorem
there exists a complex measure $\mu = \mu_n$ supported in $K$ with total
variation $|\mu|(K)=1$, such that $\mu$ annihilates the vector space ${\cal
P}_n$ of holomorphic polynomials of degree at most $n$
\ (that is, $\int_K p_n \, d\mu =0$ for all $p_n \in {\cal P}_n$) and
$$
  d_n = \int_K f \, d\mu.
$$
Since $F=f$ on $K$,  we can write
$$
  d_n = \int_K F \, d\mu = (\mu*\check F)(0).
\eqno(11)
$$
where $\check F(z):=F(-z)$. Now form the convolution
$$
  \mu * \check F = (\mu* \check F) * \delta = (\mu* \check F) * {\partial \over \partial \bar z} E = {\partial \over \partial \bar z} \check F * (\mu*E),
\eqno (12)
$$
where  $\delta$ is the point mass at $0$
and $E(z) := 1/(\pi z)$ is the Cauchy
kernel.  Associativity of the triple convolution holds since each term has compact support. We write
$$
  \hat \mu (z)
  :=
  (\mu*E)(z) = {1 \over \pi} \int_K {d\mu(\zeta) \over z - \zeta},
$$
the Cauchy transform of $\mu$.  Note that $\mu*E$ is
holomorphic outside $K$.  From (11) and (12) we then obtain
(note that $\partial F/\partial \bar z=0$ on $D_\rho$)
$$
  d_n
  =
  \int_{D_R \setminus D_{\rho}} (\mu*E)(z) {\partial \over \partial \bar z} \check F(z) \, d A (z),
\eqno(13)
$$
where $A$ is Lebesgue measure in $\CC$.

In order to utilize formula (13) for $d_n$, we
need estimates for $\hat \mu$.  We first note that
since  $|\mu| (K) = 1$,  we have
$$
  |(\mu*E)(z)| \leq M, \quad   z \in \partial D_r,
\eqno(14)
$$
for some constant $M > 0$ depending only on the distance from
$K$ to $\partial D_r$.  In addition we have the growth estimate
$$
  |(\mu*E)(z)| = O (1 / |z|^{n+1}) \quad {\rm as\ } |z| \to \infty;
\eqno(15)
$$
this follows from noting that for $z$  sufficiently large
we have
$1/(z - \zeta) =  \sum_k \zeta^k / z^{k+1}$
uniformly for $\zeta \in K$, and then using the fact that
$\mu$ satisfies
$\int_K \zeta ^k \, d\mu(\zeta)=0$ for $0 \le k \leq n$.
We now consider the function
$$
  u(z)
  :=
  g_{K} (z) + {1\over n} \log \left( {|\mu*E(z)| \over M} \right).
$$
Using (14) and (15), we see that $u(z)$ is subharmonic in
$\CC \setminus \bar{D}_r$,
bounded at $\infty$, and continuously assumes
values which are at most $\log r$ on
$\partial D_r$.  By the maximum
principle we have $u(z) \leq \log r$ on $\CC \setminus D_r$; that is,
$$
  |(\mu*E)(z)|
  \leq
  M[ e^{\log r - g_{K} (z)}]^n, \quad z \in \CC \setminus D_r.
\eqno(16)
$$

From (13) and (16) we conclude that
$\displaystyle{
\limsup _{n \to \infty} d_n (f,K) ^{1/n} \leq r/\rho
}$.
Now let $r \downarrow 1$ and $\rho \uparrow R$. \eop

Note that for each non-constant polynomial $p$
with  $\norm{p}_K \leq 1$, we have $ {1 \over \deg p} \log|p(z)|\leq g_K(z)$ so that
$$\max \left\{ 0 ,
    \sup _p \left\{ {1 \over \deg p} \log|p(z)| \right\}
    \right\}\leq g_K(z).\eqno(17)$$
It turns out that equality holds in (17). We use this as a starting point in jumping to several complex variables in the next section.

\sect{The Bernstein-Walsh theorem in $\CC^N, \ N>1$} For a compact set $K\subset \CC^N$, we may define
$$
  V_K(z)
  := \max \left\{ 0 ,
    \sup _p \left\{ {1 \over \deg p} \log|p(z)| \right\}
    \right\}$$
where the supremum is taken over all non-constant polynomials $p$
with  $\norm{p}_K \leq 1$.  This is a generalization of the one-variable Green function $g_K$. Note that from the definition of the polynomial hull $\hat K$, we have
$$V_K=V_{\hat K}.$$
The function $V_K$ is lower semicontinuous, but it need not be
upper semicontinuous.  The upper semicontinuous regularization
$$
  V_K^*(z) = \limsup_{\zeta \to z} V_K (\zeta)
$$
of $V_K$ is either identically $+\infty$ or else $V_K^*$ is
plurisubharmonic. The first case occurs if the set $K$ is too ``small''; precisely if $K$ is \dword{pluripolar}: this means that there exists a psh function $u$ defined in a neighborhood of $K$ with $K\subset \{z:u(z)=-\infty\}$ (see section 13 for more on pluripolar sets). We say
that $K$ is \dword{$L$-regular} if $V_K=V_K^*$, that is, if $V_K$ is
continuous.  For example, if $\CC^N\setminus K$ is regular with respect to $\RR^{2N}$-potential theory, then $K$ is $L$-regular. A simple example is a closed Euclidean ball $K=\{z\in \CC^N: |z-a|\leq R\}$; in this case, $V_K(z)=V_K^*(z)=\max [0,\log |z-a|/R]$. For a product $K=K_1\times \cdots \times K_N$ of planar compact sets $K_j\subset \CC$, $V_K(z_1,\ldots,z_N)=\max_{j=1,\ldots,N}g_{K_j}(z_j)$. In particular, for a polydisk
$$P:=\{(z_1,\ldots,z_N): |z_j-a_j|\leq r_j, \ j=1,\ldots,N\},$$
$V_K(z_1,\ldots,z_N)=\max_{j=1,\ldots,N}[0,\log |z_j-a_j|/r_j]$. Any compact set $K$ can be approximated from above by the decreasing sequence of $L$-regular sets $K_n:=\{z:{\rm dist}(z,K)\leq 1/n\}$. The reason for the ``$L$'' is that
the class of plurisubharmonic
functions $u$ in $\CC^N$ \dword{of logarithmic growth}, i.e., such that
$u(z)\leq \log|z| +C, \ |z|\to \infty$, is called the class $L=L(\CC^N)$. The functions ${1 \over \deg p} \log|p(z)|$ for a polynomial $p$ clearly belong to $L$; historically, for any Borel set $E$, the
function
$$
  V_E(z):= \sup \{ u(z): u\in L, \ u\leq 0 \ {\rm on}  \ E\}
$$
was called the \dword{$L$-extremal function of $E$} and it was proved that for compact sets $K$, this upper envelope
coincides with that in the beginning of this section. We sketch a proof of this. An important feature of the proof is the correspondence between psh functions in $L(\CC^N)$ and   ``homogeneous'' psh functions in $\CC^{N+1}$. We remind the reader of the standard correspondence between polynomials $p_d$ of degree $d$ in $N$ variables and homogeneous polynomials $H_d$ of degree $d$ in $N+1$ variables via
$$p_d(z_1,\ldots,z_N)\mapsto H_d(w_0,\ldots,w_N):=w_0^dp_d(w_1/w_0,\ldots,w_N/w_0).$$

Clearly $ V_K(z)\leq V(z):=\sup \{ u(z): u\in L, \ u\leq 0 \ {\rm on}  \ K\}$ and to prove the reverse inequality, by approximating $K$ from above, if necessary, we may assume $K$ is $L$-regular. We consider $h(z,w)$ defined for $(z,w)\in \CC^{N+1}=\CC^N\times \CC$ as follows:
$$h(z,w):= \cases{|w|\exp V(z/w),&  $w\not =0$;\cr
          \limsup_{(z',w')\to (z,0)}h(z',w'), & $w=0$.}$$

This is a nonnegative {\it homogeneous} psh function in $\CC^{N+1}$; i.e.,  $h(tz,tw)=|t|h(z,w)$ for $t\in \CC$. We say that the function $\log h$ is \dword{logarithmically homogeneous}: $\log h(tz,tw) = \log |t| + \log h(z,w)$. Fix a point $(z_0,w_0)\not =(0,0)$ with $z_0/w_0\not\in K$ and fix $0<\epsilon <1$. Using the fact that the polynomial hull coincides with the hull with respect to continuous psh functions (see the end of section 3), it follows that the compact set
$$E:=\{(z,w)\in \CC^{N+1}:h(z,w)\leq (1-\epsilon)h(z_0,w_0)\}$$
is polynomially convex. Moreover, $E$ is \dword{circled}: $(z,w)\in E$ implies $(e^{it}z,e^{it}w)\in E$ for all real $t$.
\vskip6pt

\noindent {\bf Exercise}. {\sl Given a compact, circled set $E\subset \CC^N$ and a polynomial $p_d = h_d +h_{d-1} + \cdots + h_0$ of degree $d$ written as a sum of homogeneous polynomials, we have $\norm{h_j}_E \leq \norm{p_d}_E, \ j=0,\ldots,d$.} Hint: Fix a point $b\in E$ at which $|h_j(b)|=\norm{h_j}_E$ and use Cauchy's estimates on $\lambda \mapsto p_d(\lambda b)=\sum_{j=0}^d \lambda^j h_j(b)$.
\vskip6pt

From the exercise, the polynomial hull of our circled set $E$ is the same as the hull obtained using only homogeneous polynomials. Since $E=\hat E$ and $(z_0,w_0)\not \in E$, we can find a homogeneous polynomial $h_s$ of degree $s$ with $|h_s(z_0,w_0)|> \norm{h_s}_E$. Define
$$p_s(z,w):={h_s(z,w)\over \norm{h_s}_E}\cdot [(1-\epsilon)h(z_0,w_0)]^s.$$
Then $|p_s(z,w)|^{1/s}\leq |h(z,w)|$ for $(z,w)\in \partial E$ and by homogeneity of $|p_s|^{1/s}$ and $h$ we have $|p_s|^{1/s}\leq h$ in all of $\CC^{N+1}$. At $(z_0,w_0)$, we have
$$|p_s(z_0,w_0)|^{1/s}> (1-\epsilon)h(z_0,w_0);$$
since $\epsilon >0$ was arbitrary, as was the point $(z_0,w_0)$ (provided $z_0/w_0\not\in K$), we get that
$$h(z,w)=\sup_s \{|p_s(z,w)|^{1/s}: p_s \ \hbox{homogeneous of degree} \ s, \ |p_s|^{1/s}\leq |h|\}.$$
At $w=1$, we obtain
$$\exp V(z)=h(z,1)=\sup_s \{|Q_s(z)|^{1/s}: Q_s \ \hbox{of degree} \ s, \ |Q_s|^{1/s}\leq \exp V\}$$
which proves the result (note $V\leq 0$ on $K$).

If the compact set $K \subset \CC^N$ is $L$-regular,
then for each $R > 1$ we define the set
$$
  D_R := \{ z : V_K (z) < \log R \};
\eqno(18)
$$
this is an open neighborhood of $\hat K$ and we clearly have the Bernstein-Walsh inequality
$$
  |p(z)| \leq \norm{p}_K R^{ \deg p }=\norm{p}_{\hat K} R^{ \deg p },
  \qquad   z \in D_R
\eqno (19)
$$
for every polynomial $p$ in $\CC^N$. Theorem (BW1) goes over
exactly to several complex variables:

\proclaim Theorem (BWN). \  Let $K$ be an $L$-regular
compact set in $\CC^N$.
Let $R > 1$, and let $D_R$ be defined by {\rm (18)}.
Let $f$ be continuous on $K$. Then
$$
  \limsup_{n\to \infty} d_n(f,K)^{1/n}\leq 1/R
$$
if and only if $f$ is
the restriction to $K$ of a function holomorphic in $D_R$.

\noindent {\bf Sketch of proof}. The ``only if'' direction follows since $K$ satisfies the
Bernstein-Walsh inequality (19). For the converse, we may assume that $K=\hat K$. Fix $f\in {\cal O}(D_R)$ and $\epsilon >0$. Since $\partial D_R$ is compact and $V_K(z)=\max \left\{ 0 ,
    \sup _p \left\{ {1 \over \deg p} \log|p(z)| \right\}
    \right\}$, we can find finitely many polynomials $p_1,\ldots,p_m$ of degree $d$, say, with $\norm{p_j}_K\leq 1$ such that
    $$\max_j\{{1\over d} \log {|p_j(z)|}\}> \log R -\epsilon  \quad\hbox{on}  \quad\partial D_R.$$
By slightly modifying the $p_j$'s, if necessary, we may assume that
$$\Pi:=\{z: |q_j(z)|<1, \ j=1,\ldots,m\}$$
is a Weil polyhedron where $q_j(z):={p_j(z)/(e^{-\epsilon d}R^d)}$; and clearly $K\subset \Pi$ and $\Pi$ approximates $D_R$; i.e., $D_{R-\delta}\subset \Pi \subset D_R$ for small $\delta$. Now use the Weil integral formula from section 3 to expand $f(z)$ and truncate the series to get good polynomial approximators.  \eop
\vskip6pt

\noindent This sketch follows the outline of the proof given by Siciak [Si].  Zaharjuta [Z2] gave the first
proof of Theorem (BWN).

There is an interesting related result, due to Tom Bloom, which applies pluripotential theory to multivariate approximation theory. To motivate this, we return to the one-variable situation and let $K\subset \CC$ be compact with $\CC\setminus K$ connected and $g_K$ continuous. Let $W(K)$ denote the closure (in the uniform norm on $K$) of the functions holomorphic on a neighborhood of $K$. Given $f\in W(K)$, let $B_d(z)=b_dz^d+\cdots$ be the best approximant to $f$ (in sup-norm on $K$) from ${\cal P}_d(\CC)$. Wojcik [W] showed that $f$ has a holomorphic extension to $D_R$ for some $R>1$ if and only if
$$\limsup_{d\to \infty} |b_d|^{1/d} \leq {1\over R{\rm cap}(K)}.$$
Here, ${\rm cap}(K):=\lim_{|z|\to \infty}|z|\exp(-g_K(z))$ is the \dword{logarithmic capacity of $K$}. Equivalently,
$${\rm cap}(K)=\lim_{d\to \infty} \inf \{\norm{p_d}_K^{1/d}: p_d(z) =z^d + \cdots\},$$
the \dword{Chebyshev constant} of $K$.

Now in several complex variables, what should replace the ``leading coefficient'' $b_d$ of a univariate polynomial? Moreover, what replaces the asymptotics of the Green function $g_K$? To address the first question, recall that any polynomial $P_d$ of degree $d$ may be written as the sum $P_d=H_d +H_{d-1}+\cdots +H_0$ where $H_j$ is a homogeneous polynomial of degree $j$; i.e., $H_j(tz)=t^jH_j(z)$ for $t\in \CC, \ z\in \CC^N$. Going backwards, {\it given} a homogeneous polynomial $H_d$ of degree $d$ and a compact set $K\subset \CC^N$, we define the \dword{Chebyshev polynomial of $H_d$ relative to $K$}, denoted ${\rm Tch}_KH_d$, to be a polynomial of the form $H_d +R_{d-1}$ with $\norm{H_d+R_{d-1}}_K$ minimal among all $R_{d-1}\in {\cal P}_{d-1}(\CC^{N})$. Such a polynomial need not be unique if $N\geq 2$ but the number $\norm{{\rm Tch}_KH_d}_K$ is well-defined. Note if $N=1$ and $H_d(z)=z^d$, ${\rm Tch}_KH_d$ is just the classical Chebyshev polynomial for $K$ of degree $d$.

For the second question, we make some preliminary definitions. Given a psh function $u\in L(\CC^N)$ we define the \dword{Robin
function} of $u$ to be
$$
\rho_u(z):=\limsup_{|\lambda|\to \infty} \left[u(\lambda z)-
\log |\lambda| \right].
$$
Note that for $\lambda \in \CC$,
$\rho_u(\lambda z)=\log {|\lambda|} + \rho_u(z)$; i.e., $\rho_u$
is logarithmically
homogeneous. It is known ([Bl6], Proposition 2.1) that for
$u\in L(\CC^N)$, the Robin function $\rho_u(z)$ is plurisubharmonic
 in ${\CC}^N$; indeed, either
$\rho_u\in L(\CC^N)$ or $\rho_u\equiv -\infty$. As an example, if $p$
is a polynomial of degree $d$ so that $u(z):={1\over d}\log
{|p(z)|}\in L(\CC^N)$,
 then $\rho_u(z)={1\over
d}\log {|\hat p(z)|}$ where $\hat p$ is the top degree $(d)$
homogeneous
 part of $p$. For a compact set $K$, we denote by $\rho_K$ the Robin
function of $V_K^*$; i.e., $\rho_K := \rho_{V_K^*}$.

We can now state the beautiful result of Bloom:

\proclaim Theorem ([Bl6]). \ Let $K$ be an $L$-regular, polynomially convex compact set in $\CC^N$. Let $f\in W(K)$ and let $B_d:= H_d + {\rm lower \ degree \ terms}$, $d=1,2,\ldots$, be a sequence of best approximating polynomials to $f$ on $K$. For $R>1$, the following are equivalent:
\item {1.} $f$ extends holomorphically to $D_R$;
\item {2.} $\limsup_{d\to \infty} \norm{f-B_d}_K^{1/d}\leq 1/R$;
\item {3.} $\limsup_{d\to \infty} \norm{{\rm Tch}_KH_d}_K^{1/d}\leq 1/R$;
\item {4.} $\limsup_{d\to \infty} {1\over d}\log {|H_d(z)|} \leq \rho_K(z) -\log R$ for $z\in \CC^N\setminus \{0\}$.

Of course, the equivalence of 1. and 2. is the Bernstein-Walsh theorem. The deep part of this result is the implication 4. implies 3.; this follows from

\proclaim Theorem ([Bl6]). \ Let $K$ be an $L$-regular, polynomially convex compact set in $\CC^N$. Let $\{H_d\}$ be a sequence of homogeneous polynomials with $\deg H_d=d$ and assume that $\limsup_{d\to \infty} {1\over d}\log {|H_d(z)|} \leq \rho_K(z)$ for $z\in \CC^N\setminus \{0\}$. Then
$$\limsup_{d\to \infty} \norm{{\rm Tch}_KH_d}_K^{1/d}\leq 1.$$

To prove this theorem, Bloom constructs polynomials $W_j, \ j=1,\ldots,s$ with $\norm{W_j}_K\leq 1$ such that a Weil polyhedron $\{z\in \CC^N:|\hat W_j(z)|< R_j, \ j=1,\ldots,s\}$ utilizing the top degree homogeneous polynomials $\hat W_j$ of $W_j$ contains $K$. Each $H_d$ may be expanded in a series involving the $\hat W_j$'s:
$$H_d (z)= \sum_{|M|=0}^{\infty}{\sum_I }'A_M^I(z)[\hat W_I(z)]^M.$$
Replacing each $\hat W_I$ by $W_I$ creates polynomials
$$P_d (z):= \sum_{|M|=0}^{\infty}{\sum_I}' A_M^I(z)[W_I(z)]^M$$
of degree $d$ which are competitors for ${\rm Tch}_KH_d$ and whose sup norms on $K$ can be estimated.

\sect{Quantitative Runge-type results in multivariate
approximation} The duality proof presented of the one-variable
Walsh theorem, Theorem (BW1) of section 4, may be extended to
yield a quantitative Runge theorem for harmonic functions in
$\RR^N$, where $N \ge 2$. To state this we let ${\cal H}^N_n$ be
the vector space of all harmonic, real-valued polynomials of $N$
variables of degree at most $n$. If $f$ is a continuous
real-valued function on a compact set $K \subset \RR^N$, we now
define
$$
  d_n(f,K):= \inf \{ \norm{f-h_n}_K: h_n\in {\cal H}^N_n\}.
$$

\proclaim Theorem ([A],[BL2]).
\ Let $K$ be a compact subset of $\RR^N$ such that
$\RR^N \setminus K$ is connected.  \vskip0pt
{\rm (a) \ } Let $\Omega$ be an open neighborhood of $K$.  Then
there is a constant $\rho \in (0,1)$, depending only on $K$ and
$\Omega$, with the following property:  if $f$ is
harmonic on $\Omega$, then $\limsup_{n\to \infty}
d_n (f,K) ^{1/n} \leq \rho$.
\vskip0pt
{\rm (b) \ } Suppose that when we regard
$K \subset \RR^N = \RR^N + i0 \subset \CC^N$,
the set $K$ is $L$-regular. If $f$ is a real-valued continuous function
on  $K$ such that
$\limsup _{n \to \infty} d_n (f,K) ^{1/n} < 1$, then
$f$ extends to a harmonic function on an open neighborhood of $K$.

\noindent {\bf Sketch of proof of (a)}. We follow the outline of the duality proof of the Walsh theorem, replacing the Cauchy kernel by the fundamental solution
$E(x-y)=c_N|x-y|^{2-N}$ for the Laplace operator $\Delta$
(here $c_N$ is a constant depending only on the dimension).
Analogous to (13), we can write
$d_n=\int_L(\mu*E)(z)\Delta F(z)dA(z)$
where $\mu=\mu_n$ is a complex measure supported on $K$ with
$|\mu|(K)=1$ which annihilates ${\cal H}_n^N$, $L$ is a compact subset of $\Omega\setminus K$, and
$F$ is a $C^\infty$ function.
The details of this proof may be found in [BL2], but
the main difference here is the problem of estimating the \dword{Newtonian potential}
$$
  (\mu *E)(x)= \int_K c_N |x-y|^{2-N} \, d\mu(y)
$$
on compact subsets of $\Omega \setminus K$ under the assumption that
we have the decay estimate
$$
  |(\mu*E)(x)| = O(|x|^{-n}) \quad  {\rm as}  \quad |x| \to \infty
\eqno(20)
$$
analogous to (15).  We utilize the Kelvin transform
$T(x) := x/|x|^2 $, under which a harmonic function $h(x)$ on a
domain $G \subset \RR^N \setminus \{ 0 \}$ is
transformed into a harmonic function $\tilde h(x):= |x|^{2-N}
h(x/|x|^2)$ on $T(G)$.
The condition (20) of rapid decay at
infinity is transformed into a condition
of flatness near the origin, and the problem
of estimating $\mu * E$  under the hypothesis (20)
is reduced to proving the following
\dword{Schwarz lemma} for harmonic functions [BL2].

\proclaim Lemma (HSL). \
Let $\Omega$ be a bounded domain in $\RR^N$ and
let $a\in \Omega$. If $K$ is a compact subset of $\Omega$,
then there exist constants $C>1$ and $\rho \in (0,1)$,
depending only on $K$ and $\Omega$, with the following property:
if $f$ is a harmonic function on $\Omega$ satisfying
$|f| \leq 1 \ {\rm in\ } \Omega$,
and if $D^{\alpha}f(a)=0$ whenever $|\alpha|<n$, then
$\norm{f}_K\leq C \rho ^n$.

Here $D^{\alpha}f={\partial^{|\alpha|}f\over \partial x_1^{\alpha_1}\cdots  \partial x_N^{\alpha_N}}$. The proof of Lemma (HSL) is based on techniques
from the theory of functions of several complex variables.
In the preceding section we introduced the
function $V_K$ in $\CC^N$ as a substitute for the
ordinary Green function with pole at infinity in $\CC^1$.
To prove Lemma (HSL) we introduce in $\CC^N$
a substitute for the Green function with a finite pole
in $\CC^1$.  Following Klimek [K], section 6.1, we define
for each domain $\widetilde \Omega \subset \CC^N$
and each point $a \in \widetilde \Omega$
the \dword{pluricomplex Green function}
$$
  G_{\widetilde \Omega}(z;a)
  :=
  \sup_{u} u(z),
$$
where the supremum is taken over all
nonpositive plurisubharmonic functions $u$ on
$\widetilde \Omega$ such that
$u(z) - \log|z-a|$ has an upper bound in some neighborhood of $a$.
It is known that if, e.g., $\widetilde \Omega$ is bounded, or, more generally, hyperconvex (see section 3), then $G_{\widetilde \Omega} (\cdot;a) $ is a nonconstant, negative
plurisubharmonic function in $\widetilde \Omega$
(see [K], [BL2]).    This fact leads to the following
Schwarz lemma for holomorphic functions of several variables
(see [Bi], [BL2]).

\proclaim Lemma (SL). \
Let $\widetilde \Omega$ be a bounded domain in $\CC^N$ and
let $a\in \widetilde \Omega$. Let $\tilde K$ be a compact subset of
$\widetilde \Omega$.
If $f$ is a holomorphic function on $\widetilde \Omega$ satisfying
$|f|\leq 1 \ {\rm in} \ \widetilde \Omega$,
and if $\partial^{\alpha}f(a)=0$ whenever  $|\alpha|<n$,  then
$\norm{f}_{\tilde K}\leq \rho ^n$, where
$$
  \rho := \sup_{\tilde K} \exp(G_{\widetilde \Omega}(\cdot;a))< 1.
$$

Here $\partial^{\alpha}f={\partial^{|\alpha|}f\over \partial z_1^{\alpha_1}\cdots  \partial z_N^{\alpha_N}}$.
The inequality $\rho < 1$ is clear from the fact that
$G_{\widetilde \Omega} (\cdot;a) $
is a negative function on $\widetilde\Omega$ which is
{\it subharmonic} as a function of
$2N$ real variables.  The rest of the lemma follows from
the fact that the function $u(z) :=  {1 \over n} \log |f(z)|$
is one of the competitors in the definition
of $G_{\widetilde \Omega} (\cdot;a) $.

\bigskip
A harmonic
function is real analytic and thus about each point $x_0$ in our domain
$\Omega\subset \RR^N$
we can get a power series expansion $\sum a_{\alpha}(x-x_0)^{\alpha}$ of $f$ which converges {\it as a
holomorphic function} $\sum a_{\alpha}(z-x_0)^{\alpha}$ in a neighborhood of this point in
$\CC^N$. We prove Lemma (HSL) by covering the compact set $K$ by finitely many
real balls $B$ and the union of the complex balls $\widetilde B$ gives a
neighborhood $\widetilde \Omega$ of $K$
in $\CC^N$ to which we can apply Lemma (SL).  \eop
\vskip6pt

In [BL3], the authors prove Bernstein theorems for solutions of more general elliptic
partial differential equations. Let $p(x) := \sum _{|\alpha| = m} a_{\alpha} x^{\alpha}$
be a non-constant homogeneous polynomial in $\RR^N$,
with complex coefficients, which is never equal to zero
on $\RR^N \setminus \{ 0 \}$; here $N \ge 2$. Then the partial
differential operator
$p(D) := p( \partial /\partial x_1, \dots, \partial /\partial x_n)$
is elliptic.
We let ${\cal L}_n$ be the vector space
of polynomials $q$ of degree at most $n$ in
$N$ variables which are solutions of the equation $p(D)q=0$.
If $f$ is a continuous function on a compact set
$K \subset \RR^N$, define
$$
  d_n(f,K)
  =
  \inf \{\norm{f-p_n}_K: p_n\in {\cal L}_n\}.
$$

\proclaim Theorem ([BL3]).
\ Let $K$ be a compact subset of $\RR^N$ such that
$\RR^N \setminus K$ is connected. \vskip0pt
{\rm (a) \ } Let $\Omega$ be an open neighborhood of $K$.
Then there is
a constant $\rho \in (0,1)$, depending only on $p(D)$, $K$, and
$\Omega$, with the following property:  if $f$ is
a solution of $p(D)f = 0$ on $\Omega$, then
$\limsup_{n\to \infty} d_n (f,K) ^{1/n} \leq \rho$.
\vskip0pt
{\rm (b) \ } Suppose that when we regard
$K \subset \RR^N = \RR^N + i0 \subset \CC^N$,
the set $K$ is $L$-regular. If $f$ is a real-valued continuous function
on  $K$ such that
$\limsup _{n \to \infty} d_n (f,K) ^{1/n} < 1$, then
$f$ extends to a solution $F$ of $p(D)F = 0$
on an open neighborhood of $K$.

A duality proof, utilizing a fundamental solution $E(x-y)$ for the operator $p(D)$, and yet another tool from pluripotential theory, the \dword{relative extremal function} $$\omega^*(z,F,\Omega):=[\sup \{u(z):
 u \  \hbox{psh in} \ \Omega, \ u\leq 0, \ u|_F \leq -1\}]^*$$ of a set $F\subset \Omega$ relative to $\Omega$, may be found in [BL3]. See section 13 for more on $\omega^*(z,F,\Omega)$. The need for this function arises as we don't have a Kelvin transform in this general setting; but we do get a ``transfer of smallness'' result analogous to Lemma (SL):

\proclaim Lemma. \ Let $\Omega$ be a bounded domain in $\CC^N$. Let $F\subset \Omega$ be nonpluripolar and let $K\subset \Omega$ be compact. Then there is a constant $a\in (0,1]$ such that for any holomorphic function $g$ on $\Omega$ with $|g|\leq M$ on  $\Omega$ and $|g|\leq m<M$ on $F$, we have
$$|g|\leq m^aM^{1-a} \ \hbox{on} \ K.$$

The proof is trivial: simply observe that $\log |g|$ is psh and by the definition of $\omega(z,F,\Omega)$, it follows that
$$u(z):={\log {(|g(z)|/M)} \over \log {(M/m)}}\leq \omega (z,F,\Omega).$$
Then the constant $a$ can be chosen to be $a(\Omega,F,K):= -\sup_K  \omega(z,F,\Omega)$. The nonpluripolarity of $F$ insures that $a>0$ (see Proposition ($\omega$) in section 13).

Finally, we mention that Jackson-type approximation theorems for solutions to $p(D)f=0$ on $\Omega$ which are continuous on $\bar \Omega$ can be found in [BBL1] and [BBL2].

\sect{Mergelyan property and solving $\bar \partial$} Recall the Cauchy-Green formula: Let $\Omega$ be a bounded domain in $\CC$ with $C^1$-boundary and let $f\in C^1(\bar \Omega)$. Then
$$f(z)={1\over 2\pi i}\int_{\partial \Omega} {f(\zeta)\over \zeta - z} d\zeta -{1\over \pi} \int_{\Omega} {\partial f \over \partial \bar z}(\zeta) \cdot ({1\over \zeta -z}) dA(\zeta)\eqno(21)$$
where $dA$ denotes Lebesgue measure. In particular,
\item {1.} if $f \in {\cal O}(\Omega)\cap C(\bar \Omega)$,
$$f(z)={1\over 2\pi i}\int_{\partial \Omega} {f(\zeta)\over \zeta - z} d\zeta \ \hbox{(Cauchy integral formula)};\eqno(22)$$
\item {2.} if $f\in C^1_0(\Omega)$,
$$f(z)=-{1\over \pi} \int_{\Omega} {\partial f \over \partial \bar z}(\zeta) \cdot ({1\over \zeta -z}) dA(\zeta).$$

\noindent As an immediate corollary, if $g\in C_0(\Omega)$, then
$$G(z):=-{1\over \pi} \int_{\Omega} g(\zeta) \cdot ({1\over \zeta -z}) dA(\zeta)$$ solves the inhomogeneous Cauchy-Riemann equation
$\partial G / \partial \bar z=g$
 in $\Omega$. Moreover, we see that
$$\sup_{\Omega} |G| \leq [\sup_{\Omega} |g|]  \sup_{z\in \Omega}\bigl[{1\over \pi}\int_{{\rm supp}g}  |{1\over \zeta -z}| dA(\zeta)\bigr]\eqno(23)$$
so that $\sup_{\Omega} |G| \leq C \sup_{\Omega} |g|$.
More generally, if $\mu$ is a measure with compact support in $\Omega$, the Cauchy transform of $\mu$,
$$\hat \mu (z) :=-{1\over \pi} \int_{\Omega}{1\over \zeta -z}d\mu(\zeta),$$
satisfies $\partial \hat \mu / \partial \bar z=\mu$ in the sense of distributions on $\Omega$.

Now suppose $K=\bar \Omega$ where $\Omega$ is a simply connected domain with boundary of class $C^1$. An elementary proof of Mergelyan's theorem (Theorem (Me) from section 1) for such $K$ goes as follows:
\vskip4pt
\item {1.} {\sl Smooth approximation}. Cover $\partial K$ by finitely many open sets $U_1,\ldots,U_n$ such that for each $j=1,\ldots,n$ there is a vector $t_j$ transverse to $\partial \Omega$ at each point of $\partial \Omega \cap U_j$ pointing outward (into $\CC\setminus \bar \Omega$). Take a partition of unity $\phi,\phi_1,\ldots,\phi_n$ for a neighborhood of $\bar \Omega$ subordinate to the cover consisting of $\Omega,U_1,\ldots,U_n$. Given $f\in {\cal O}(\Omega)\cap C(\bar \Omega)$, for sufficiently large $j$,
$$g_j(z):=\phi(z) f(z) + \sum_{k=1}^n \phi_k(z) f(z- t_k/j)$$
is defined and $C^{\infty}$ on a neighborhood $\Omega_j$ of $\bar \Omega$. Since $f\in C(\bar \Omega)$, $g_j\to f$ uniformly on $\bar \Omega$.

\item {2.} {\sl Holomorphic correction}. We have that
$$ {\partial g_j \over \partial \bar z}(z)= f(z)\cdot {\partial \phi \over \partial \bar z} (z)+ \sum_{k=1}^n  f(z- t_k/j)\cdot {\partial \phi_k \over \partial \bar z}(z)$$
is uniformly small on $\Omega_j$ since $\phi + \sum_k \phi_k =1$ there; say $\sup_{\Omega_j}  |{\partial g_j \over \partial \bar z}|\leq \delta_j$ where $\delta_j \to 0$. From the previous discussion, utilizing (23), for each $j$ we can find $G_j \in C^{\infty}(\Omega_j)$ with
$${\partial G_j \over \partial \bar z} = {\partial g_j \over \partial \bar z}$$
in $\Omega_j$ and $\sup_{\Omega_j} |G_j| \leq C_j \delta_j \to 0$. Then $f_j:=g_j - G_j\in {\cal O}(\Omega_j)$ and $f_j \to f$ uniformly on $\bar \Omega$.
\vskip4pt

In higher dimensions, the smooth approximation step works fine. However, things get tricky in step 2 for two reasons:
\vskip4pt
\item {(i)} we need to solve a $\bar \partial$-equation; more precisely,
\item {(ii)} we need to solve a $\bar \partial$-equation {\it with uniform (sup-norm) estimates}.
\vskip4pt

\noindent Let's make this precise. Given a $C^1$ function $u$, the $1$-form $du$ can be written as $du=\partial u +\bar \partial u$ where
$$\partial u:=\sum_{j=1}^N {\partial u\over \partial z_j}dz_j$$
is a form of bidegree $(1,0)$ and
$$\bar \partial u:=\sum_{j=1}^N{\partial u\over \partial \bar z_j}d\bar z_j$$
is a form of bidegree $(0,1)$. In general, a differential form $\phi$ of bidegree $(p,q)$ is a sum
$$\phi =\sump_{|I|=p, \ |J|=q}c_{I,J} dz^I \wedge d\bar z^J$$
where $c_{I,J}$ are functions ($0$-forms) and
$$dz^I = dz_{i_1}\wedge \cdots \wedge dz_{i_p}; \ d\bar z^J = d\bar z_{j_1}\wedge \cdots \wedge d\bar z_{j_q};$$
the prime means the indices are increasing. We define
$$\bar \partial \phi = \sump_{|I|=p, \ |J|=q}\bar \partial c_{I,J} \wedge dz^I \wedge d\bar z^J;$$
this is a form of bidegree $(p,q+1)$.

We can extend the operator $\partial$ to $(p,q)$-forms as well (the result is a $(p+1,q)$-form). Since any form $\omega$ of degree $r\in \{0,1,\ldots,2N\}$ can be written as a sum of forms $\omega_{p,q}$ of bidegree $(p,q)$ where $0\leq p, q \leq r$ and $p+q=r$, we extend $\bar \partial$ and $\partial$ to general forms by linearity. Note then as differential operators on the space of smooth forms, we have
$$d^2 =d\circ d =0 =(\partial +\bar \partial)\circ  (\partial +\bar \partial) =\partial^2 + \partial \circ \bar \partial + \bar \partial \circ \partial + \bar \partial^2;$$
by bidegree considerations
$$\partial^2 =\bar \partial^2=0; \quad \partial \circ \bar \partial =- \bar \partial \circ \partial$$
(e.g., if $\phi$ is a $(p,q)$ form, $d^2\phi=0$ is a form of total degree $p+q+2$; $ \bar \partial^2 \phi$ is of bidegree $(p,q+2)$ and there are no other terms with this bidegree). In particular, let $\phi =\sum_{j=1}^N \phi_j d \bar z_j$ be a smooth $(0,1)$ form on a domain $\Omega$ in $\CC^N$. {\it If we want to be able to find a function $u\in C^{\infty}(\Omega)$ with $\bar \partial u = \phi$ in $\Omega$, a necessary condition is that $\bar \partial \phi=0$}. This condition is vacuous if $N=1$ (there are no $(0,2)$ forms in $\CC$). Note that the inhomogeneous Cauchy-Riemann equation $\bar \partial u = \phi$ is a system of $2N$ (real) partial differential equations in $\RR^{2N}$ for the (two) unknown functions $\Re u$ and $\Im u$.  This is an overdetermined system if $N>1$.

Are there integral formulas providing solutions to $\bar \partial$? For $N=1$, define, for $z\in \CC$, the $(1,0)$-form in $\zeta$
$$\omega_{BM}(\zeta - z) := {1\over 2\pi i}{d\zeta \over \zeta -z} = {1\over 2\pi i}{ \bar \zeta - \bar z\over |\zeta -z|^2}d\zeta.$$
Then if $\Omega\subset \CC$ is a bounded domain with $C^1$-boundary and $f\in {\cal O}(\Omega)\cap C(\bar \Omega)$, we have
$$f(z)=\int_{\partial \Omega} f(\zeta) \omega_{BM}(\zeta - z)$$
for $z\in \Omega$. This is the Cauchy integral formula (22). In SCV, we have the following.

\proclaim Proposition (Bochner-Martinelli formula). \ Define the $(N,N-1)$-form
$$\omega_{BM}(\zeta - z) := {(N-1)!\over (2\pi i)^N}\sum_{j=1}^N{ (-1)^{j-1}(\bar \zeta_j - \bar z_j)\over |\zeta -z|^{2N}}d\bar \zeta[j]\wedge d\zeta.$$
If $\Omega\subset \CC^N$ is a bounded domain with $C^1$-boundary and $f\in {\cal O}(\Omega)\cap C(\bar \Omega)$, then
$$f(z)=\int_{\partial \Omega} f(\zeta) \omega_{BM}(\zeta - z)\eqno(24)$$
for $z\in \Omega$.

Here $d\bar \zeta[j]=d\bar \zeta_1\wedge \cdots \wedge \hat {d\bar \zeta_j} \wedge \cdots \wedge d\bar \zeta_N$ (omit $d\bar \zeta_j$) and $d\zeta = d\zeta_1 \wedge \cdots \wedge d\zeta_N$. The reader will note that the $\zeta_j$ partial derivative of ${1\over |\zeta -z|^{2N-2}}$, a fundamental solution for the Laplacian in $\RR^{2N}$, is
$${\partial \over \partial \zeta_j}{1\over |\zeta -z|^{2N-2}}={ (1-N)(\bar \zeta_j - \bar z_j)\over |\zeta -z|^{2N}},$$
which is, up to a constant, the coefficient of $d\bar \zeta[j]\wedge d\zeta$ in $\omega_{BM}(\zeta - z)$.
As a generalization of the Cauchy-Green formula (21), for any $f\in C^1(\bar \Omega)$ we have
$$f(z)=\int_{\partial \Omega} f(\zeta) \omega_{BM}(\zeta - z)-{(N-1)!\over \pi ^N}\int_\Omega \sum_{j=1}^N {\partial f \over \partial \bar z_j}{\bar \zeta_j -\bar z_j \over |\zeta -z|^{2N}}dA(\zeta) \eqno(25)$$
for $z\in \Omega$. Here $dA(\zeta)$ is Lebesgue measure in $\CC^N$, i.e., $$dA(\zeta)=(i/2)^Nd\zeta_1 \wedge d \bar \zeta_1\wedge \cdots \wedge d\zeta_N \wedge d \bar \zeta_N.$$

However, {\it only if $N=1$ are the coefficients of this Bochner-Martinelli kernel $\omega_{BM}(\zeta - z)$ holomorphic in $z$}; thus only for $N=1$ can this formula be used to construct solutions to $\bar \partial u =\phi$. That is, given a smooth $(0,1)$ form $\phi=\sum_{j=1}^N \phi_j d \bar z_j$ on $\Omega \subset \CC^N$ with $\bar \partial \phi=0$, from (25) we'd {\it like} to define
$$u(z)=-{(N-1)!\over (2\pi i)^N}\int_\Omega \sum_{j=1}^N \phi_j{\bar \zeta_j -\bar z_j \over |\zeta -z|^{2N}}d\bar \zeta \wedge d\zeta$$
to solve $\bar \partial u = \phi$ in $\Omega$, but since
$(\bar \zeta_j -\bar z_j)/ |\zeta -z|^{2N}$ is not holomorphic in $z$ if $N>1$, this doesn't work.  A major area of research in SCV was the attempted construction of integral formulas with {\it holomorphic} kernels. For {\it strictly pseudoconvex domains} (recall section 3), this was done by Henkin and Ramirez. The article [He] of Henkin is a nice historical survey on the subject.

We say that a domain $\Omega$ has the \dword{Mergelyan property} if every $f\in {\cal O}(\Omega)\cap C(\bar\Omega)$ can be approximated uniformly on $\bar\Omega$ by functions in ${\cal O}(\bar\Omega)$. This definition is natural, given the one-variable Mergelyan theorem, Theorem (Me), from section 1.

\proclaim Theorem. \ A smoothly bounded, strictly pseudoconvex domain satisfies the Mergelyan property.

\noindent The reason is that Henkin, Kerzman and Lieb showed that one can solve $\bar \partial$ with {\it uniform} estimates by constructing holomorphic kernels in this case.

What if $D$ is pseudoconvex, but not strictly pseudoconvex? Let
$$D=\{(z,w)\in \CC^2: 0<|z|<|w|<1\},$$ the so-called \dword{Hartogs triangle} (draw a picture in $|z|,|w|-$space to explain the terminology). Then $D$ is pseudoconvex. However, any function holomorphic on a neighborhood of $\bar D$ must necessarily extend holomorphically to the unit bidisk $P=\Delta \times \Delta$ (exercise; or see [Sha] sections 7 and 40). Now
consider  the function $f(z,w):=z^2/w$. Note that
$$\limsup_{(z,w)\to (0,0), \ (z,w) \in D}|f(z,w)|\leq \limsup_{(z,w)\to (0,0), \ (z,w) \in D}|w|^2/|w|= 0$$
so that $f$ extends continuously to $(0,0)$. Hence $f\in {\cal O}(D)\cap C(\bar D)$. We show that $f$ is not uniformly approximable on $\bar D$ by holomorphic functions on $\bar D$. For suppose $\{f_j\} \subset {\cal O}(\bar D)$ converge uniformly to $f$ on $\bar D$. In particular, uniform convergence of $\{f_j\}$ on
$$T^2=\partial \Delta \times \partial \Delta=\{(z,w):|z|=|w|=1\}$$ implies uniform convergence of $\{f_j\}$ on  $\bar P$ (recall the multivariate Cauchy integral formula (3) from section 1). This limit function, call it $g$, is necessarily holomorphic on $P$; in particular, it is holomorphic at $(0,0)$. Moreover, $g$ must coincide with $f$ on $D$. But $f$ does not extend holomorphically to $(0,0)$.

There are examples of smoothly bounded pseudoconvex domains which do not satisfy the Mergelyan property. It is conjectured that {\sl if $D$ is a smoothly bounded pseudoconvex domain, then $D$ has the Mergelyan property if and only if there are pseudoconvex domains $D_j$ with $\bar D \subset D_j$ such that $\bar D = \cap_j D_j$}. This latter property fails for the Hartogs triangle. It also fails in all of the known examples of smoothly bounded pseudoconvex domains which fail to satisfy the Mergelyan property. We refer the reader to the article of Bedford and Fornaess [BF] for a more detailed discussion.

\sect{Approximation on totally real sets} Recall from section 1 the example of the two polynomially convex disks $K_1$ and $K_2$ in $\CC^2$ defined as
$$K_1:=\{(x_1,x_2)\in \RR^2: x_1^2+x_2^2 \leq 1\} \ \hbox{and}$$
$$K_2:=\{(z_1,0): |z_1|\leq 1\}.$$
Here $P(K_1)=C(K_1)$ but $P(K_2)\not = C(K_2)$. The difference is that $K_1$ lies in the totally real submanifold $\RR^2$ of $\CC^2$.

\proclaim Definition. \ Let $\Sigma$ be a submanifold of class $C^1$ of an open set $D\subset \CC^N$. We say $\Sigma$ is \dword{totally real} if for each $p\in \Sigma$, the tangent space $T_p\Sigma$ contains no complex lines; i.e., no complex linear subspaces of positive dimension.

In particular, the dimension of such a (real) submanifold is at most $N$. A top-dimensional example is the torus $T^N=\{(z_1,\ldots,z_N):|z_j|=1, \ j=1,\ldots,N\}$. Returning to the simpler example of $\RR^N=\RR^N+i0$, note that
$$u(z):={\rm dist}(z,\RR^N)^2= y_1^2+\cdots + y_N^2$$
is of class $C^2$ and strictly psh (compute the complex Hessian!). Of course, directly from the definition, $\RR^N=\{z\in \CC^N: u(z)=0\}$. Indeed, {\sl if $\Sigma$ is a totally real submanifold of class $C^2$ of an open set $D\subset \CC^N$, then there exists a neighborhood $\omega$ of $\Sigma$ such that $u(z)={\rm dist}(z,\Sigma)^2$ is of class $C^2$ and strictly psh on $\omega$} (cf. [AW], Lemma 17.2). The main result on approximation on totally real submanifolds is the following generalization of the real Stone-Weierstrass theorem due to Harvey and Wells [HW]. It was first proved under stronger regularity hypotheses by H\"ormander and Wermer [H\"oW]. A very enlightening proof has recently been given by Berndtsson [Ber].

\proclaim Theorem (HW). \ Let $\Sigma$ be a totally real submanifold of class $C^1$ in an open set in $\CC^N$ and let $K\subset \Sigma$ be compact and polynomially convex. Then $P(K)=C(K)$.

This is related to the question stated towards the end of section
2: {\sl Give a ``nice'' condition (C) on a compact set $K\subset
\CC^N$ so that if $K=\hat K$ then $K$ satisfies (C) if and only if
$P(K)=C(K)$}. In Lavrentiev's theorem, Theorem (La), in the
complex plane, $K^o=\emptyset$ was a necessary and sufficient
condition for a polynomially convex compact set $K$ to have this
property. Since uniform limits of holomorphic objects like
(holomorphic) polynomials should be, in some sense, holomorphic,
we seek a condition (C) which prohibits $K$ from having any type
of ``analytic structure''. Note that $P(K)$ is a uniform algebra,
i.e., a closed subalgebra of $C(K)$. There was a famous conjecture
known as the \dword{peak point conjecture}: Suppose $A$ is a
uniform algebra on its maximal ideal space $X$ such that every
point $x\in X$ is a peak point for $A$, i.e., there exists $f\in
A$ such that $f(x)=1$ and $|f(y)|<1$ for all $y\not=x$. Does it
follow that $A$ coincides with the algebra $C(X)$? In  case $X$ is
a polynomially convex compact set in $\CC^N$ and $A=P(X)$, we are
asking if this ``peak point property'' suffices as a condition
(C). A counterexample given by Cole in 1968 shows that the answer
to the general peak point conjecture is no. Anderson, Izzo and
Wermer [AIW1], [AIW2] have shown that {\sl if $\Sigma$ is a
compact polynomially convex {\bf real analytic} variety in $\CC^N$
such that every point in $\Sigma$ is a peak point for $P(\Sigma)$,
then $P(\Sigma)=C(\Sigma)$.} Recall that a relatively closed
subset $V$ of an open set $U$ in $\CC^N$ is a real analytic
subvariety of $U$ if for each $z_0\in V$ there exists a
neighborhood $U'\subset U$ of $z$ and real valued, real analytic
functions $f_1,...,f_m$ in $U'$ with
$$V\cap U'=\{z\in U':f_1(z)=\cdots = f_m(z)=0\}.$$

The unit sphere in $\CC^N$ for $N>1$ is a smooth submanifold in $\CC^N$ which definitely {\it has} complex tangents (i.e., is {\it not} totally real) as the dimension of the sphere is $2N-1 >N$; however, it is straightforward to see that for any compact subset $K$ of the sphere, $P(K)$ has the peak point property (e.g., at $(1,0,\ldots,0)$, take $f(z)=z_1$).   Despite this, Izzo [I] has constructed examples of the following:
\vskip4pt
\item {1.} {\sl There exists a compact polynomially convex subset $K$ of the unit sphere in $\CC^3$ such that $P(K)\not= C(K)$}.
\item {2.} {\sl There exists a $C^{\infty}$-embedding $F: \CC^2\rightarrow \CC^5$ such that the set $K = F(\{(z,w): |z|\leq 1, |w| = 1\})$ is a compact polynomially convex subset of the unit sphere in $\CC^5$ which satisfies $P(K)\not= C(K)$}.
\vskip4pt
Note this last example is a compact, polynomially convex $C^{\infty}$ submanifold $K$ for which every point is a peak point for $P(K)$ but $P(K)\not = C(K)$; the Anderson, Izzo and Wermer theorem shows that such an example cannot occur in the real analytic category. Stout has recently strengthened the Anderson, Izzo and Wermer theorem to eliminate the hypothesis on peak points:

\proclaim Theorem ([St2]). \ Let $K$ be a compact, polynomially convex real analytic subvariety of $\CC^N$. Then $P(K)=C(K)$. \par

The aforementioned results of Harvey-Wells and/or H\"ormander-Wermer essentially reduce questions of approximation on subsets of real submanifolds on $\CC^N$ to approximation on the points where the tangent space to the manifold contains a complex line. The H\"ormander-Wermer approach to Theorem (HW), which requires some additional regularity hypotheses on $\Sigma$, can be summarized as follows: given $f\in C(K)$, we can clearly approximate $f$ uniformly on $K$ by a global smooth function; i.e., we may assume $f\in C^{\infty}(\CC^N)$. From Theorem (OW), since $K=\hat K$, it suffices now to approximate $f$ uniformly on $K$ by functions holomorphic on a neighborhood of $K$. It is straightforward to construct a function $F$ of class $C^1$ on $\CC^N$ which agrees with $f$ on $K$ and such that
$$|{\partial F\over \partial \bar z_j}(z)|=O({\rm dist}(z,\Sigma)^m), \ j=1,\ldots,N$$
if, say, $\Sigma$ is of class $C^{2m+1}$. Next, using the function $u(z)={\rm dist}(z,\Sigma)^2$, we can construct a bounded, pseudoconvex neighborhood $\omega$ of $K$ in $\CC^N$ to which we can apply standard several complex variables machinery -- solvability of $\bar \partial$ in $\omega$ -- to construct a function $G$ in $\omega$ with
$${\partial G\over \partial \bar z_j}={\partial F\over \partial \bar z_j}, \quad j=1,\ldots,N$$
 and with $|G|$ very small in $\omega$. Then $G-F$ is holomorphic in $\omega$ and approximates $f$ very well on $K$.

The $\bar \partial$ machinery utilized in the previous paragraph is the  H\"ormander $L^2$-theory. If $D$ is a smoothly bounded, pseudoconvex domain in $\CC^N$, then there is a constant $C$ depending only on $D$ such that for any $(0,1)$-form $\phi=\sum_{j=1}^N\phi_j d\bar z_j$ with $L^2(D)$-coefficients satisfying $\bar \partial \phi =0$, there exists $u\in L^2(D)$ with $\bar \partial u =\phi$ in $D$ and
$$\int_D|u|^2 dA \leq C\int_D \sum_{j=1}^N |\phi_j|^2 dA$$
(cf., [H\"o], Chapter 4 or [AW] section 16). Note that this is a global $L^2$-norm estimate. From this, one gets local interior regularity of solutions sufficient to derive the required estimate on $G$ in the previous paragraph. Berndtsson uses a ``weighted'' version of the global $L^2$-norm estimate in his work.

The Harvey-Wells approach uses integral kernels to solve $\bar \partial$ and is at least similar in spirit  to our outline of the proof of Theorem (OW) using the Weil integral formula. Extensions of the Harvey-Wells result have been made by Range and Siu [RS] as well as by Bruna and Burg\'es [BB]. These papers deal with approximation in H\"older norms on a totally real {\it compact} subset $X\subset \CC^N$: we assume there exists a strictly plurisubharmonic $C^2$ function in a neighborhood of $X$ whose zero set is $X$.

Using a generalization of the Bochner-Martinelli kernel, a suitably constructed {\it Cauchy-Fan\-tap\-pi\`e-Leray kernel} (see section 10), Weinstock [Wei] has proved an interesting perturbation of the Stone-Weierstrass theorem. We discuss this briefly. Given any compact set $K\subset \CC^N$ and functions $f_1,..,f_m\in C(K)$, let $[f_1,\ldots,f_m]$ denote the algebra generated by these functions. Note then we always have $[z_1,\ldots,z_N,\bar z_1,\ldots,\bar z_N]$ is dense in $C(K)$. Suppose $N$ functions $R_1,\ldots,R_N$ are given. Let $A=
[z_1,\ldots,z_N,\bar z_1+R_1,\ldots,\bar z_N+R_N]$. Under what conditions on $R:=(R_1,\ldots,R_N)$ is $A$ dense in $C(K)$? Assume each $R_j$ is defined and continuous in a neighborhood $U$ of $K$.

\proclaim Theorem (PSW). \ If there exists $0\leq k <1$ with
$$|R(z)-R(z')|\leq k|z-z'| \quad \hbox{for} {\ }{\quad} z,z'\in U,$$
then $A$ is dense in $C(K)$.

For $N=1$ this result is due to Wermer [We]. Even in this case, it is the Lipschitz norm of the perturbation $R$ that matters, {\it not} the supremum norm. For example, if $K=\bar \Delta$, the closed unit disk in $\CC$, and
$$R(z):=-\bar z \ \hbox{if} \ |z|\leq \epsilon; \ R(z):={-\epsilon \bar z\over |z|} \hbox{if} \ \epsilon \leq |z| \leq 1,$$
then $|R(z)|\leq \epsilon$ on $\bar \Delta$ but $\bar z +R(z)\equiv 0$ on the disk $\{z:|z|\leq \epsilon\}$. Thus each function in $A$ must be holomorphic on $\{z:|z|<\epsilon\}$. A nice exposition of the one and several variable results can be found in chapter 14 of [AW].

\sect{Lagrange interpolation and orthogonal polynomials} A natural
way to construct polynomials which approximate a given function is
to use interpolating polynomials. Let $K$ be a polynomially convex
$L$-regular compact subset of $\CC^N$. Let $m_n={N+n\choose n}$
denote the dimension of the complex vector space ${\cal P}_n$ of
polynomials in $N$ complex variables of degree at most $n$.

For each integer $n\geq 1$, let $a_{n1},\ldots,a_{nm_n}\in K$. Thus we have a doubly indexed array $(a_{nj})_{n=1,\ldots; \ j=1,\cdots, m_n}$ of points in $K$. Given a function $f$ holomorphic in a neighborhood of $K$, under what conditions on the array do the Lagrange polynomials $L_nf$ interpolating $f$ at the points $(a_{nj})_{{j=1},\cdots,m_n}$ converge uniformly to $f$ on $K$?

In one variable, Walsh gave a necessary and sufficient condition on the array in order to guarantee uniform convergence of $\{L_nf\}$ to $f$ on $K$ for {\it all} such $f$. In several variables $(N\geq 2)$, much less is known because there is no analogue of the Hermite remainder formula  used in the proof of Walsh. We remind the reader of the
Hermite remainder formula
for interpolation of a holomorphic function of one variable. This is a simple consequence of the Cauchy integral formula.
Let $z_1,\ldots,z_n$ be $n$
distinct points in the plane and let $f$ be a
function which is defined at these points. The functions
$$l_j(z):=\prod_{k\not= j}(z-z_k)/(z_j-z_k), \quad j=1,\ldots,n,$$
are polynomials of degree $n-1$ with $l_j(z_k)=\delta_{jk}$, called the
\dword{fundamental Lagrange interpolating polynomials}
associated to $z_1,\ldots,z_n$.  Then $(L_nf)(z):=\sum_{j=1}^n
f(z_j)l_j(z)$ is the unique polynomial of degree at most $n$
satisfying $(L_nf)(z_j)=f(z_j), \ j=1,\ldots,n;$ we call it the \dword{Lagrange
interpolating polynomial}
associated to $f,z_1,\ldots,z_n.$ If $\Gamma$ is a
rectifiable Jordan curve such that
the points $z_1,\ldots,z_n$ are inside $\Gamma$, and
$f$ is holomorphic inside and on $\Gamma$, we can estimate
the error in our approximation of $f$ by $L_nf$ at points
inside $\Gamma$ using the following formula.

\proclaim Lemma (Hermite Remainder Formula).  For any $z$ inside
$\Gamma$,
$$
  f(z)-(L_nf)(z)
  =
  {1 \over 2\pi i} \int_{\Gamma} {\omega(z) \over \omega(t)}
  {f(t) \over (t-z)} \, dt, \eqno(26)
$$
where $\omega (z):=\prod _{k=1} ^n (z-z_k).$

Note that if $f(z)=1/(t-z)$, then
$$f(z)-(L_nf)(z)= {\omega(z)\over \omega(t)} {1\over t-z}.\eqno(27)$$
The necessary and sufficient condition of Walsh on the univariate array $(a_{nj})\subset K\subset \CC$ so that $L_nf \to f$ uniformly on $K$ for any $f$ holomorphic on $K$ is that the polynomials $\omega_{n+1}(z):=\prod _{j=1} ^{n+1} (z-a_{nj})$ satisfy
$$\lim_{n\to \infty} \norm{\omega_{n+1}}_K^{{1\over n+1}} ={\rm cap}(K);$$
equivalently, the subharmonic functions
$$u_n(z):={1\over n+1} \log {|\omega_{n+1}(z)|\over \norm{\omega_{n+1}}_K}$$ converge locally uniformly to $g_K$ on $\CC\setminus K$.

For a survey of some results in the several variable case, we refer the reader to [BBCL],  [Bl3] and [BlL1]. We outline the elementary positive results. Let $e_1,\ldots,e_{m_n}$ form a basis for ${\cal P}_n$. Given
$A_n=\{a_{n1},\ldots,a_{nm_n}\}\subset K$ we form the \dword{generalized
Vandermonde determinant}
$$
  V_n(A_n) := \det[e_i(a_{nj})]_{i,j=1,\ldots,m_n}.
$$
If $V_n(A_n)\not= 0$, we can form the polynomials
$$
  l_{nj}(z):= {V_n(a_{n1},\ldots,z,\ldots,a_{nm_n}) \over V_n(A_n)},
  \qquad j=1,\ldots,m_n
$$
satisfying $l_{nj}(a_{ni})= \delta_{ji}$.
We call
$$
  \Lambda_n:= \sup_{z\in K} \sum_{j=1}^{m_n} |l_{nj}(z)|
$$
the $n$-th \dword{Lebesgue constant} for $K,A_n.$ For $f$ defined on $K$,
$$
  (L_nf)(z):= \sum_{j=1}^{m_n} f(a_{nj})l_{nj}(z)
$$
is the Lagrange interpolating polynomial for $f$ at the points $A_n$.
We say that $K$ is \dword{determining} for $\bigcup {\cal P}_n$
if whenever $h\in \bigcup {\cal P}_n$ satisfies
$h=0$ on $K$, it follows that $h\equiv 0.$ For these sets we can find
point sets $A_n$ for each $n$ with $V_n(A_n)\not= 0.$ We have
the following elementary result.

\proclaim Proposition. \
Let $K$ be determining for $\bigcup {\cal P}_n$ and let
$A_n\subset K$ be sets of points satisfying $V_n(A_n)\not= 0$ for each $n.$ Given $f$ bounded
on $K$, if $\limsup \Lambda_n^{1/n}=1$, then $\limsup \norm{f-L_nf}_K^{1/n}=
\limsup d_n^{1/n}$ where
$$d_n=d_n(f,K)=\inf\{\norm{f-p_n}_K: p_n\in {\cal P}_n\}.$$

\noindent {\bf Proof.} Fix $\epsilon>0$ and choose, for each $n$, a
polynomial $p_n\in {\cal P}_n$ with $$\norm{f-p_n}_K^{1/n}\leq d_n^{1/n}+
\epsilon.$$
Since $p_n\in {\cal P}_n$, we have $L_np_n=p_n$ and
$$\eqalign{\norm{f-L_nf}_K\;=\;&\norm{f-p_n+L_np_n-L_nf}_K\cr
  \;\leq\;& \norm{f-p_n}_K + \Lambda_n\norm{f-p_n}_K=
(1+\Lambda_n)\norm{f-p_n}_K.\cr}$$
Using the hypothesis $\limsup \Lambda_n^{1/n}=1$, we obtain the conclusion. \eop

Let $p(D) := p( \partial /\partial x_1, \dots, \partial /\partial x_n)$ be an elliptic partial
differential operator as at the end of section 6, and let ${\cal L}_n$ be the vector space
of polynomials $q$ of degree at most $n$ in
$N$ variables which are solutions of the equation $p(D)q=0$. The same proof shows:  {\sl Let $K\subset \RR^N$ be determining for $\cup_n{\cal L}_n$ and let $A_n\subset K$ satisfy $V_n(A_n)\not= 0$ for each $n.$ Given $f$ bounded
on $K$, if $\limsup \Lambda_n^{1/n}=1$, then $\limsup \norm{f-L_nf}_K^{1/n}=
\limsup d_n^{1/n}$}. Here, we replace $m_n$ by $\tilde m_n := {\rm dim}{\cal L}_n$, $A_n= \{a_{n1},\ldots,a_{n\tilde m_n}\}\subset K$ and $d_n$ is defined just as before,
but now with respect to ${\cal L}_n$.

Arrays of points $\{A_n\},n=1,2,\ldots$ satisfying $\limsup \Lambda_n^{1/n}=1$
can be constructed by taking, e.g., $A_n$ to be a set
of \dword{$n$-Fekete points for $K$}: for each $n$, choose $A_n\subset K$ so
that
$$\max_{X_n\subset K} |V_n(X_n)|=|V_n(A_n)|.$$
Since $|V_n(X_n)|$ is a continuous function on $K^{m_n}$, such points
exist. Moreover, from the definition of Fekete
points, $\norm{l_{nj}}_K=1$ so that
$\Lambda_n\leq m_n$. It is easy to see that $\lim m_n^{1/n}=1$. The problem is that in several variables these points are essentially impossible to construct. The first known {\it explicit} example of an array $\{A_n\},n=1,2,\ldots$ satisfying $\limsup \Lambda_n^{1/n}=1$ associated to a compact set $K\subset \CC^N, \ N>1$, has been recently discovered by Bos, et al [BCDVX]. The set $K$ is the unit square $[-1,1]\times [-1,1]$ in $\RR^2$. In this example, the Lebesgue constants have minimal possible asymptotic growth: $\Lambda_n =O([\log n]^2)$.

In CCV, for a non-polar compact set $K\subset \CC$, the normalized
counting measures associated to Fekete arrays satisfy $$ \mu_n:={1\over n+1}
\sum_{j=1}^{n+1} \delta_{A_{nj}}\to {1\over 2\pi} \Delta g_K$$ in the weak$^*$-topology as measures. Here, $\Delta g_K$, the Laplacian of $g_K$, is to be
interpreted as a positive distribution, i.e., a positive measure. Indeed, for any array  $\{A_n\},n=1,2,\ldots$ satisfying $\limsup \Lambda_n^{1/n}=1$ the same conclusion holds (cf. [BBCL]). In SCV, for $K\subset \CC^N$ nonpluripolar, there is a conjecture that for Fekete arrays, the normalized discrete measures
$$\mu_n:={1\over m_n} \sum_{j=1}^{m_n} \delta_{A_{nj}}$$ converge weak-* to the \dword{Monge-Amp\'ere measure} $\mu_K :=(dd^cV_K^*)^N$ of the $L$-extremal function $V_K^*$. For a discussion of the complex Monge Amp\'ere operator $(dd^c\cdot)^N$, see the appendix. To this date, {\it nothing} is known if $N>1$. Some special situations, which really reduce to one-variable problems, can be found in [GMS] and [BlL2].

A widely studied topic in classical approximation theory is the study of orthogonal polynomials. Let $\mu$ be a positive Borel measure with compact support $K={\rm supp} (\mu) \subset \CC^N$. Assume that the set $K$ is determining for $\cup{\cal P}_n$. If $N=1$, this just means that $K$ contains infinitely many points; for $N>1$, $K$ being nonpluripolar is sufficient (but not necessary). Then the standard basis monomials $\{e_{\alpha}(z):= z_1^{\alpha_1}\cdots z_N^{\alpha_N}\}$ are linearly independent in $L^2(\mu)$ and one can form the orthonormal polynomials $\{p_{\alpha}(z,\mu)\}$. For an introduction to this topic in $\CC^N$, we recommend Tom Bloom's paper [Bl5] which concerns the relationship between the so-called $n$-th root asymptotic behavior of the orthonormal polynomials $\{p_{\alpha}(z,\mu)\}$ and the ``pluripotential theory'' of the set $K$. This is presented in a systematic manner analogous to the one-variable study developed in the book of H. Stahl and V. Totik [ST]. See also [Bl8].

The pair $(K,\mu)$ is said to have the \dword{Bernstein-Markov property} if for each $\epsilon >0$ there exists a positive constant $M=M(\epsilon)$ such that
$$\|p\|_K\leq M(1+\epsilon)^{\deg p}\|p\|_{L^2(\mu)}$$ for all polynomials $p=p(z)$. For an $L$-regular compact set $K$, such measures always exist; e.g., the Monge-Amp\`ere measure $\mu_K=(dd^cV_K)^N$ (this is the Laplacian $\Delta g_K$ if $N=1$). One can even find such a measure $\mu$ which is rather ``sparse'' in the sense that there exists a countable subset $E\subset K$ with $\mu(E)=\mu(K)$. Returning to the setting of the Bernstein-Walsh theorem, given such a measure, best $L^2(\mu)$-approximants to $f\in C(K)$ have optimal behavior.

\proclaim Proposition. \ Let $K$ be a polynomially convex $L$-regular compact set in $\CC^N$ and let $\mu$ be a measure supported on $K$ such that $(K,\mu)$ satisfies the Bernstein-Markov property. If $f\in C(K)$ satisfies $\limsup_{n\to \infty}d_n(f,K)^{1/n}=\rho <1$, and if $\{p_n\}$ is a sequence of best $L^2(\mu)$-approximants to $f$, then $\limsup_{n\to \infty} \norm{f-p_n}_K^{1/n}=\rho$.

The proof follows trivially from the fact that if $\rho <r<1$ and $\{q_n\}$ are polynomials with $\norm{f-q_n}_K\leq Mr^n$ for some $M$ (independent of $n$), then
$$\norm{f-p_n}_{L^2(\mu)}\leq \norm{q_n-f}_{L^2(\mu)}\leq \norm{q_n-f}_K\mu(K)^{1/2}\leq Mr^n\mu(K)^{1/2}.$$
For simplicity we take $\mu(K)=1$. Then we have $\norm{p_n-p_{n-1}}_{L^2(\mu)}\leq Mr^n(1+1/r)$ which shows that $p_o+\sum_{n=1}^{\infty} (p_n-p_{n-1})$ converges to $f$ in $L^2(\mu)$ and pointwise $\mu$-a.e.\ to $f$ on $K$. By the Bernstein-Markov property, for each $\epsilon < 1/r - 1$ there exists $\tilde M>0$ with
$$\norm{p_n-p_{n-1}}_K\leq \tilde M (1+\epsilon)^n\norm{p_n-p_{n-1}}_{L^2(\mu)}\leq \tilde M [(1+\epsilon)r]^nM(1+1/r)$$
showing that $p_o+\sum_{n=1}^{\infty} (p_n-p_{n-1})$ converges uniformly to a continuous function $g$ on $K$ (holomorphic on the interior of $K$). Since $f$ and $g$ are continuous and $g=f$ $\mu$-a.e.\ on $K$, $g=f$ on $K$. Then
$$\norm{f-p_n}_K =\norm{\sum_{k=n+1}^{\infty} (p_k-p_{k-1})}_K \leq \tilde M [(1+\epsilon)r]^{n+1}M{(1+1/r) \over [1-(1+\epsilon)r]}$$ showing that $\limsup_{n\to \infty} \norm{f-p_n}_K^{1/n}\leq (1+\epsilon)r$. Again, a similar result holds in the elliptic partial
differential operator case if $K\subset \RR^N$.

\sect{Kergin interpolation} A more promising type of interpolation procedure has been successfully applied to many approximation problems by Tom Bloom and his collaborators. A natural extension of Lagrange interpolation to $\RR^s, \ s>1$ was discovered by P. Kergin (a student of Bloom) in his thesis. Indeed, Kergin interpolation acting on ridge functions (a univariate function composed with a linear form) is Lagrange interpolation. The Kergin interpolation polynomials generalize to the case of $C^m$ functions in $\RR^N$ both the Lagrange interpolation polynomials and those of Hermite.

As brief motivation, given $f\in C^m([0,1])$, say, and given $m+1$ points $t_0< \cdots < t_m\in [0,1]$, if one constructs the Lagrange interpolating polynomial $L_mf$ for $f$ at these points, then there exist (at least) $m-1$ points between pairs of successive $t_j$ at which $f'$ and $(L_mf)'$ agree; then there exist (at least) $m-2$ points between triples of successive $t_j$ at which $f''$ and $(L_mf)''$ agree, etc. Given a set $A=[A_0,A_1,\ldots,A_m]\subset \RR^N$ of $m+1$ points and $f$ a function of class $C^m$ on a neighborhood of the convex hull of these points, there exists a unique polynomial ${\cal K}_{A}(f)={\cal K}_{A}(f)(x_1,\ldots,x_N)$ of total degree $m$ such that ${\cal K}_{A}(f)(A_j)=f(A_j)$, $j=0,1,\ldots,m$, and such that for every integer $r$, $0\leq r\leq m-1$, every subset $J$ of $\{0,1,\ldots,m\}$ with cardinality equal to $r+1$, and every homogeneous differential operator $Q$ of order $r$ with constant coefficients, there exists $\xi$ belonging to the convex hull of the $(A_j)$, $j\in J$, such that $Qf(\xi)=Q{\cal K}_{A}(f)(\xi)$. In [Bl1], Bloom gives a proof of this result by using a formula due to Micchelli and Milman [MM] which gives an explicit expression for ${\cal K}_{A}(f)$.
If $f=u+iv$ is holomorphic in a convex region $D$ in $\CC^N$, and if $A=[A_0,A_1,\ldots,A_m]\subset D \subset \CC^N=\RR^{2N}$, then we can construct ${\cal K}_{A}(u)$ and ${\cal K}_{A}(v)$. It turns out (cf.\ [Bo2]) that ${\cal K}_{A}(u)+i{\cal K}_{A}(v)$ is a holomorphic polynomial.

An alternate description, which we give in the holomorphic setting, is as follows (cf.\ [BC2]). Let $D$ be a $\CC$-convex domain in $\CC^N$, i.e., the intersection of $D$ with any complex line is connected and simply connected. Note that in $\RR^N$ this is the same condition as convexity if we replace ``complex line'' by ``real line.'' For any set ${\cal A}=[A_0,\ldots,A_d]$ of (not necessarily distinct) $d+1$ points in $D$ there exists a unique linear projector ${\cal K}_{{\cal A}}: {\cal O}(D)\to {\cal P}_d$ (recall that ${\cal O}(D)$ is the space of holomorphic functions on $D$ and ${\cal P}_d$ is the space of polynomials of $N$ complex variables of degree less than or equal to $d$) such that
\vskip4pt
\item {(i)} ${\cal K}_{{\cal A}}(f)(A_j)=f(A_j)$ for $j=0,\cdots,d$,
\item {(ii)} ${\cal K}_{{\cal A}}(g\circ {\lambda})={\cal K}_{{\lambda}({\cal A})}(g)\circ {\lambda}$ for every affine map ${\lambda}: \CC^N\to \CC$ and $g\in {\cal O}({\lambda}(D))$, where ${\lambda}({\cal A})=({\lambda}(A_0),\ldots,{\lambda}(A_d))$, \item {(iii)} ${\cal K}_{{\cal A}}$ is independent of the ordering of the points in ${\cal A}$, and
\item {(iv)} ${\cal K}_{\cal B}\circ {\cal K}_{{\cal A}}={\cal K}_{\cal B}$ for every subsequence $\cal B$ of ${\cal A}$.
\vskip4pt
\noindent The operator ${\cal K}_{{\cal A}}$ is called the Kergin interpolating operator with respect to ${\cal A}$.

Set ${\cal K}_d:={\cal K}_{{\cal A}_d}$ with ${\cal A}_d=[A_{d0},\ldots,A_{dd}]$ and $A_{dj}$ in a compact subset $K$ of $D\subset \CC^N$ for every $j=0,\ldots,d$ and $d=1,2,3,\ldots$. {\sl Under what conditions on the array $\{{\cal A}_d\}_{d=1,2,\ldots}$ is it true that ${\cal K}_d(f)$ converges to $f$ uniformly on $K$ as $d\to {\infty}$ for every function $f$ holomorphic in some neighborhood of $\bar D$}? Bloom and Calvi [BC2] attacked this problem with the aid of an integral representation formula for the remainder $f-{\cal K}_d(f)$ proved by M. Andersson and M. Passare [AP]. Their solution reads as follows. Assume that the measures  $\mu_d=(d+1)^{-1}\sum_{j=0}^d{\delta}_{A_{dj}}$ converge weak-* as $d\to {\infty}$ to a measure $\mu$. In one variable, the answer comes from potential theory: one considers the logarithmic potential
$$V_{\mu}(u):=\int_K \log {|u-t|}d\mu(t)$$
and the required condition is that
$$\{u\in \CC: V_{\mu}(u) \leq \sup_K V_{\mu}\} \subset D.$$
For $N>1$, given a linear form $p:\CC^N\to \CC$, define $\mu^p=p_*\mu$ as the push-forward of $\mu$ to $\CC$ via $p$, i.e., for $f\in C_0(\CC)$,
$$\mu^p(f):=\int_{\CCs} fd\mu^p=\mu(f\circ p):=\int_{\CCs^N}(f\circ p)d\mu.$$
Set
$${\Psi}_{\mu}(p,u):=\mu^p(\log|u-\cdot|)=\int_{\CCs} \log |u-\zeta|d\mu^p(\zeta),$$
and let $M_{\mu}(p)$ be the maximum of $u\mapsto  {\Psi}_{\mu}(p,u)$ on $p(K)$. {\sl If $D$ has $C^2$ boundary and $\{u\in \CC: {\Psi}_{\mu}(p,u)\leq M_{\mu}(p)\} \subset p(D)$ for every linear form $p$ on $\CC^N$, then ${\cal K}_d(f)$ converges to $f$ uniformly on $K$ as $d\to {\infty}$ for every function $f$ holomorphic in some neighborhood of $\bar D$}.

We call an array $\{{\cal A}_d\}_{d=1,2,\ldots}$ \dword{extremal} for $K$ if ${\cal K}_d(f)$ converges to $f$ uniformly on $K$ for {\bf each} $f$ holomorphic in a neighborhood of $K$. In the setting of subsets $K$ of $\RR^N$, Bloom and Calvi proved the following striking result.

\proclaim Theorem ([BC3]). \ Let $K\subset \RR^N, \ N\geq 2$, be a compact, convex set with nonempty interior. Then $K$ admits extremal arrays if and only if $N=2$ and $K$ is the region bounded by an ellipse.

For the Andersson-Passare remainder formula one needs an integral formula with
\vskip4pt
\item {1.} a {\it holomorphic} kernel; moreover, one with
\item {2.} a kernel that is the composition of a univariate function with an affine function.
\vskip4pt
\noindent Together with property (ii) of the Kergin interpolating operator, this allows a reduction of the multivariate problem to a univariate setting.

For $a,b\in \CC^N$, we write $\inpro{a,b}:=\sum_{j=1}^N a_jb_j$. Let $D\subset \CC^N$ be a bounded domain with smooth boundary and fix $N$ functions $w_j(\zeta), \ j=1,\ldots,N$ which are defined and smooth on $\partial D$ and satisfy
$$\inpro{w(\zeta),\zeta-z}=\sum_{j=1}^N w_j(\zeta)(\zeta_j -z_j)\not = 0\eqno(28)$$
for all $z\in D$ and $\zeta \in \partial D$. We give examples of such $w_j$ below. Define
$$\Omega(s,t):={(N-1)!\over (2\pi i)^N}\sum_{j=1}^N{ (-1)^{j-1} t_j \over \inpro{s,t}^{N}}d t[j]\wedge ds.$$
Here $dt[j]=dt_1\wedge \cdots \wedge \hat {dt_j} \wedge \cdots \wedge dt_N$ (omit $dt_j$) and $ds = ds_1 \wedge \cdots \wedge ds_N$.
Note that for fixed $z$, $\Omega(\zeta-z,\bar \zeta -\bar z)$ is simply the Bochner-Martinelli kernel $$\omega_{BM}(\zeta - z) := {(N-1)!\over (2\pi i)^N}\sum_{j=1}^N{ (-1)^{j-1}(\bar \zeta_j - \bar z_j)\over |\zeta -z|^{2N}}d\bar \zeta[j]\wedge d\zeta$$
which is used in the Bochner-Martinelli formula (24) from section 7. If $f\in {\cal O}(D)\cap C(\bar D)$, we have the following generalization of (24):
$$f(z)= \int_{\partial D} f(\zeta) \Omega(\zeta-z,w(\zeta))$$
$$={(N-1)!\over (2\pi i)^N}\int_{\partial D} {f(\zeta)\over [\sum_{j=1}^N w_j(\zeta)\cdot (\zeta_j -z_j)]^N}\cdot \sum_{j=1}^N (-1)^{j-1}w_j(\zeta) dw[j]\wedge d\zeta\eqno(29)$$
for $z\in D$. Note here $dw[j]=dw_1\wedge \cdots \wedge \hat {dw_j} \wedge \cdots \wedge dw_N$; thus it is the $(0,1)-$piece of each $1-$form $dw_j=dw_j(\zeta)$ that is important. This is known as a \dword{Cauchy-Fantappi\`e-Leray (CFL) formula}. Weinstock's proof of Theorem (PSW) in section 8 hinged on a judicious choice of the $w_j$'s. The Henkin, Kerzman and Lieb results mentioned in section 7 also utilize CFL-type kernels.

Let $D=\{\zeta \in \CC^N: \rho(\zeta) <0\}$ where $\rho\in C^1(\bar D)$ with $d\rho \not = 0$ on $\partial D$. Suppose that at each point $\zeta \in \partial D$ the complex tangent plane $T_p^{\CCs}(\partial D)$ lies outside of $D$, i.e.,
$$\sum_{j=1}^N {\partial \rho \over \partial \zeta_j}(\zeta)\cdot (\zeta_j -z_j)\not =0$$
for $\zeta \in \partial D$ and $z\in D$. Such a domain is called {\sl lineally convex}; convex domains are special examples. In the smoothly bounded category, lineally convex domains are the same as $\CC-$convex domains (cf. \ [APS], Chapter 2). The functions  $$w_j(\zeta)={\partial \rho \over \partial \zeta_j}(\zeta)$$
satisfy (28) and we obtain the following special case of the CFL  formula (29):
$$f(z)={(N-1)!\over (2\pi i)^N}\int_{\partial D} {f(\zeta)\over [\sum_{j=1}^N {\partial \rho \over \partial \zeta_j}(\zeta)\cdot (\zeta_j -z_j)]^N}\cdot \sum_{j=1}^N {\partial \rho \over \partial \zeta_j}(\zeta) d\bar \zeta[j]\wedge d\zeta.\eqno(30)$$
(see [Sha], chapter III for details). For example, if $D$ is the unit ball and $\rho(\zeta) = \sum_{j=1}^N \zeta_j \bar \zeta_j - 1$, we have $w_j(\zeta)=\bar \zeta$ and we get
$$f(z)={(N-1)!\over (2\pi i)^N}\int_{\partial D} {f(\zeta)\over [1-\sum_{j=1}^N \bar \zeta_j z_j]^N}\cdot \sum_{j=1}^N \bar \zeta_j d\bar \zeta[j]\wedge d\zeta.$$
Let's write $\rho'(\zeta):=({\partial \rho \over \partial \zeta_1},\ldots,{\partial \rho \over \partial \zeta_N})$. The Andersson-Passare remainder formula reads as follows.

\proclaim Theorem ([AP]). \ Let $D=\{z\in \CC^N: \rho(z) <0\}$ be a $\CC-$convex domain with $C^2$-boundary and let $f\in {\cal O}(D)\cap C(\bar D)$. Let $p_0,\ldots,p_d$ be $d+1$ points in $D$. Then
$$(f-{\cal K}_df )(z)={1\over( 2\pi i)^N}\int_{\partial D}\sum_{|\alpha|+\beta=N-1}\bigl (\prod_{j=0}^d{\inpro{\rho'(\zeta),z-p_j}\over \inpro{\rho'(\zeta),\zeta -p_j}}\bigr)$$
$$\times {f(\zeta)\partial \rho(\zeta)\wedge (\bar \partial \partial \rho(\zeta))^{N-1}\over [\prod_{j=0}^d\inpro{\rho'(\zeta),\zeta -p_j}^{\alpha_j}]\inpro{\rho'(\zeta),\zeta -z}^{\beta +1}}\eqno(31)$$
for $z\in D$ where $\alpha = (\alpha_0,\ldots,\alpha_d)$ is a multiindex and $\beta$ is a nonnegative integer.

\noindent Here, $\bar \partial \partial \rho$ is a $(1,1)$-form and
$$(\bar \partial \partial \rho)^{N-1}=\bar \partial \partial \rho \wedge \cdots \wedge \bar \partial \partial \rho \ (N-1 \ \hbox{times}).$$
Thus $\partial \rho(\zeta)\wedge (\bar \partial \partial \rho(\zeta))^{N-1}$ is an $(N,N-1)$-form. Equation (31) follows from (30) in a manner analogous to that of obtaining the Hermite Remainder Formula (26) from the Cauchy Integral Formula and an explicit formula for the remainder between the Cauchy kernel and its Lagrange interpolant, formula (27). One explicitly computes the Kergin interpolant of
$$\left[\sum_{j=1}^N {\partial \rho \over \partial \zeta_j}(\zeta)\cdot (\zeta_j -z_j)^N\right]^{-1},$$ the portion of the CFL kernel depending on the $z$-variables, using the fact that this is the composition of a univariate function with an affine function on $\CC^N$.

Bloom and Calvi [BC1] also considered what happens to multivariate Lagrange (or Hermite) interpolants $L_df$ to a given function $f$ of some minimal smoothness fixing the degree $d$ and letting the interpolation points coalesce. They give both a geometric condition and an algebraic condition sufficient for the interpolants to converge to the Taylor polynomial of the function at the point of coalescence. The proof makes an interesting use of Kergin interpolation.

There has been lots of work done in the holomorphic category; for results on Kergin interpolation of entire functions, see [Bl2] and [Bl4]. We finish this section with an interesting ``real'' result. The nature and definition of Kergin interpolants requires $C^n$ smoothness of a real-valued function $f$ in order to construct an interpolant to $f$ associated with $n+1$ points. Bos and Waldron [BW] have observed that for $n+1$ points in $\RR^N$ in general position, the Kergin polynomial interpolant of $C^n$ functions may be extended to an interpolant on all functions of class $C^{N-1}$. In particular, in $\RR^2$ one can construct Kergin interpolants of all degrees, provided points are in general position, for any $C^1$ function. Using $(n+1)$-st ``roots of unity''  $A_n$ on the unit circle in $\RR^2$, Bos  and Calvi [BC] proved that {\sl for any $f\in C^2(U)$, where $U$ is a neighborhood of the closed unit disk $K$ in $\RR^2$, $\lim_{d\to\infty}\norm{f-{\cal K}_{A_n}(f)}_K=0$}. This is a natural generalization of the analogous fact for $C^1$-functions on the interval using Lagrange interpolants at the Chebyshev nodes.

\sect{Rational approximation in $\CC^N$} Suppose $f$ is holomorphic in a neighborhood of the origin in $\CC^N$. We say that a sequence $r_ 1,r_ 2,\ldots$ of rational functions (with the degree of $r_ k$ not greater than $k$) \dword{rapidly approximates} $f$ if the $k$th root of $|f-r_ k|$ converges to zero in measure. Let $R^0$ be the class of all $f$ that admit a rapid approximation near the origin. If $N=1$ Sadullaev [Sa2] characterized the class $R^0$ in terms of Taylor coefficients.

\proclaim Theorem ($R^0$). \ Let $f(z)=\sum_{k=0}^{\infty} a_kz^k$ be holomorphic in a neighborhood of the closed unit disk in $\CC$. Define
$$A_{j_1,\ldots,j_k}:= |\det [a_{j_n+m}]_{n=1,\ldots,k; \ m=0,\ldots,k-1}|$$
and $V_k:=\sup_{j_1,\ldots,j_k}A_{j_1,\ldots,j_k}$. Then $f\in R^0$ if and only if $\lim_{k\to \infty} V_k^{1/k^2} =0$.

Sadullaev used this condition to show that {\sl a holomorphic function $f$ in a neighborhood of the origin in $\CC^N$ for $N>1$ is rapidly approximable if and only if its restriction to every complex line $L$ through the origin is rapidly approximable}. The idea of the only if direction is very simple and utilizes Hartogs series, which we used back in section 3. Via a preliminary complex-linear transformation, we may assume $f\in R^0(\CC^N)$ is holomorphic in a neighborhood of the unit polydisk and that $L=\{(z',z_N):=(z_1,\ldots,z_{N-1},z_N): z'=0\}$. We want to show that $g(z_N):=f(0,\ldots,0,z_N)$ is in $R^0(L)$. We expand $f$ in a Hartogs series
$$f(z)= \sum_{j=0}^{\infty} a_j(z') z_N^j.$$
We get a sequence of functions $V_k(z')$ defined in the closed polydisk
$$\bar U':=\{z':=(z_1,\ldots,z_{N-1}): |z_j|\leq 1\}$$
in $\CC^{N-1}$. Since each $V_k$ is the supremum of the moduli of holomorphic functions in $U'$, each function
$$u_k(z'):={1\over k^2} \log V_k(z')$$
is psh in $U'$. Since $f$ is holomorphic in a neighborhood of the unit polydisk in $\CC^N$ the coefficients $a_j$ are uniformly bounded on $\bar U'$; i.e., $|a_j(z')|\leq C$ for $z\in \bar U'$ for each $j=0,1,\ldots$.  Hence the sequence $\{u_k\}$ of psh functions is uniformly bounded above on $\bar U'$.

The key step is to show that
$$\lim_{k\to \infty} \int_{U'} u_k(z')dA (z')= -\infty.$$
To be brief, this is achieved using the fact that $f\in R^0(\CC^N)$ together with estimates on the size of certain sets involving a notion of a Chebyshev constant $T(K)$ associated to a compact set $K$. This Chebyshev constant will be defined in section 13. Since $u_k$ is psh in $U'$, it is $\RR^{2N-2}$-subharmonic; by the subaveraging property and the fact that $\{u_k\}$ is uniformly bounded above on $\bar U'$,
$$u_k(0)\leq {1\over A(U')}\int_{U'} u_k(z')dA (z')\to -\infty \ \hbox{as} \  k\to \infty.$$
This shows that $\lim_{k\to \infty} V_k(0)^{1/k^2} =0$; thus by Theorem ($R^0$), $g(z_N)=f(0,\ldots,0,z_N)\in R^0(L)$.

Gonchar [G1] showed that if $f$ is rapidly approximable then the maximal region to which $f$ continues analytically is single-sheeted  and the rapid approximation persists in this region. Taking this a step further, Sadullaev [Sa2] showed that {\bf every} {\sl holomorphic function on a domain $D$ is rapidly approximable if and only if the complement of the {\it envelope of holomorphy} of $D$ is a pluripolar set}. In this setting, the envelope of holomorphy $\tilde D$ of $D$ is the smallest domain of holomorphy containing $D$. In particular, all $g\in {\cal O}(D)$ extend holomorphically to $\tilde D$.

We remark that Bloom proved that rapid convergence in measure of a sequence $\{r_n \}$ of rational functions to a holomorphic function $f$ on an open set $\Omega \subset \CC^N$ implies rapid convergence in relative capacity (this will also be defined in section 13) on the natural domain of definition of $f$. This has the consequence that for a meromorphic function $f$ on $\CC^N$ which is holomorphic on a neighborhood of the origin the Gonchar-Pad\'{e} approximants $\{\pi_n (z,f,\lambda)\}$ converge rapidly in capacity to $f$. We refer the reader to [Bl7] for definitions and details.

In [C1], Chirka proved a ``meromorphic'' version of the Bernstein-Walsh theorem. We first describe the one-variable result. For an open set $D\subset \CC$ and a nonnegative integer $m$, let ${\cal M}_m(D)$ denote the class of meromorphic functions in $D$ which have at most $m$ poles (counted with multiplicities). Recall for a compact set $K$ in $\CC^N, \ N\geq 1$, and $R>1$, we write $D_R:=\{z\in \CC^N:V_K(z)< \log R\}$. For $f\in C(K)$ and nonnegative integers $m$ and $n$, let
$$r_{m,n}=r_{m,n}(f,K):=\inf\{\norm{f-p/q}_K: p\in {\cal P}_n, \ q\in {\cal P}_m\}.$$
The following result is due to Gonchar [G2]; a special case was proved earlier by Saff [S].

\proclaim Theorem (RBW1). \ Let $K\subset \CC$ be a regular compact set and let $R>1$. Given a continuous function $f: K\to \CC$ and a fixed integer $m\geq 0$, the following conditions are equivalent:
\item {(i)} $\limsup_{n\to \infty}(r_{m,n})^{1/n}\leq 1/R$;
\item {(ii)} there exists a function $F\in {\cal M}_m(D_R)$ with $F|_K =f$. \par

If $N>1$, the definitions of ${\cal M}_m(D)$ and the approximation numbers $r_{m,n}(f,K)$ need to be modified. We define ${\cal M}_m(D)$ to be the class of all functions in $D$ of the form $h/q_m$ where $h\in {\cal O}(D)$ and $q_m\in {\cal P}_m$.

\proclaim Theorem (RBWN). \ Let $K\subset \CC^N$ be compact and $L-$regular and let $R>1$. Given a continuous function $f: K\to \CC$ and a fixed integer $m\geq 0$, for $n\geq 1$ let
$$r^*_{m,n}=r^*_{m,n}(f,K):=\inf\{\norm{qf-p}_K: p\in {\cal P}_n, \ q\in {\cal P}_m, \ \norm{q}_K=1\}.$$
The following conditions are equivalent:
\item {(i)} $\limsup_{n\to \infty}(r^*_{m,n})^{1/n}\leq 1/R$;
\item {(ii)} there exists a function $F\in {\cal M}_m(D_R)$ with $F|_K =f$. \par

The proof of (ii) implies (i) is immediate from the standard Bernstein-Walsh theorem, Theorem (BWN) in section 5. Let $f=h/g_m\in {\cal M}_m(D_R)$. Since $h=fg_m \in {\cal O}_m(D_R)$, by Theorem (BWN) there exists a sequence $\{p_n\}$ of polynomials, $p_n\in {\cal P}_n$, with
$$\limsup_{n\to \infty} \norm{fg_m-p_n}_K^{1/n} \leq 1/R$$
which gives (i). The other implication is much deeper. Much in the spirit of Sadullaev's proof of Theorem ($R^0$), Chirka needs SCV-type capacity estimates as well as univariate arguments and techniques to achieve his goal.

Chirka constructs some interesting examples in [C1] to explain the difference between Theorems (RBW1) and (RBWN). The first example utilizes the Hartogs triangle from section 7. Precisely, let
$$K:=\bar D=\{(z,w)\in \CC^2: |z|\leq |w|\leq 1\}$$
be the closure of the Hartogs triangle $D=\{(z,w)\in \CC^2: 0<|z| < |w|< 1\}$. The polynomial hull $\hat K$ is the closed unit bidisk:
$$\hat K=\bar \Delta \times \bar \Delta=\{(z,w)\in \CC^2: |z|\leq 1, \ |w|\leq 1\}$$
since $K\subset \bar \Delta \times \bar \Delta$; $K$ contains the torus $T^2=\partial \Delta \times \partial \Delta$; and the polynomial hull of the torus $T^2$ is clearly the closed bidisk $\bar \Delta \times \bar \Delta$. Thus the $L-$extremal function $V_K$ coincides with that of the bidisk:
$$V_K(z,w)=\max[\log^+|z|,\log^+|w|]$$
so that the sublevel sets $D_R$ are larger bidisks. In particular, the set $K$ is $L-$regular so that we may apply Theorem (RBWN) to the function $f(z,w):=z^2/w$ (recall that
$f\in {\cal O}(D) \cap C(K)$). By its very definition, $f\in {\cal M}_1(D_R)$ for all $R>1$ so that
$$\limsup_{n\to \infty}(r^*_{1,n})^{1/n}\leq 1/R$$
for all $R>1$ and hence $(r^*_{1,n})^{1/n}\to 0$.
However, we cannot even uniformly approximate $f$ on $K$ by a rational function $p/q$ with $q\not =0$ on $K$, for $p/q \in {\cal O}(\bar D)$ and, as we saw in section 7, $f$ is not uniformly approximable by functions in ${\cal O}(\bar D)$. Thus for each $m$ the sequence $\{r_{m,n}\}$ does not tend to zero.

In this example, the set $K$ is not polynomially convex. However, Chirka constructs another example in which the set $K$ is a small ball $\{z\in \CC^N: |z|\leq \delta\}$ (and hence $\hat K=K$). He constructs a function $f=h/g_m\in {\cal M}_m(D_R)$ for certain $m>1$ and $R>1$ with the property that there does not exist a sequence $\{p_n/ q_m\}$ with $p_n\in {\cal P}_n$ and $q_m\in {\cal P}_m$ so that
$$\limsup_{n\to \infty} \norm{f-p_n/q_m}_K^{1/n} \leq 1/R.$$
The problem is that the ``pole-set'' of $f$ in this example cannot be written in the form $\{z\in D_R: q(z)=0\}$ for a polynomial $q$ (in this sense, the ``pole-set''  is ``nonalgebraic''). On the other hand, since $h=fg_m \in {\cal O}_m(D_R)$, we do have $\limsup_{n\to \infty}(r^*_{m,n})^{1/n}\leq 1/R$.

In [C2] Chirka utilized Jacobi series to prove holomorphic extension results. As a sample, let $g=p/q$ be a rational function in $\CC$. Let $G_r$ be a connected component of the set $\{z:|g(z)|\leq r\}$ ($G_r$ is a \dword{rational lemniscate}). If $f$ is holomorphic in a neighborhood of $\bar G_r$ then
$$F(z,w):={1\over 2\pi i} \int_{\partial G_r}{f(\zeta)\over g(\zeta)-w}\cdot {g(\zeta)-g(z)\over \zeta -z}d\zeta$$ is a holomorphic function on $G_r\times \{|w|<r\}$ which satisfies $F(z,g(z))=f(z)$ (from the Cauchy integral formula). Thus we can expand $F$ in a Taylor series in $w$ and set $w=g(z)$ to obtain a \dword{Jacobi series} for $f$:
$$f(z):=\sum_{k=0}^{\infty} C_k(z) [g(z)]^k$$
where
$$C_k(z)= {1\over 2\pi i} \int_{\partial G_r}f(\zeta)\cdot {g(\zeta)-g(z)\over [g(\zeta)]^{k+1}(\zeta -z)}d\zeta$$
are rational functions with poles at the poles of $g(z)$. The following result is proved in [C2].

\proclaim Theorem. \ Let $f$ be holomorphic in the polydisk $U'\times \{|z_N|<r\}$ where $U'$ is a polydisk in $\CC^{N-1}$. Suppose for each fixed point $p\in E \subset U'$, where $E$ is nonpluripolar in $\CC^{N-1}$, $f(p,\cdot)$ extends to a function holomorphic in $\CC=\CC_{z_N}$ except perhaps for a finite number of singularities. Then $f$ extends holomorphically to $(U'\times \CC)\setminus A$ where $A$ is an analytic variety.

Far-reaching generalizations of these ``extension'' results exist throughout the literature. For a start, consult [Iv].

\sect{Markov inequalities} The classical Bernstein-Markov inequalities  say that  for $p:\RR \to \RR$ a real
polynomial such that $\norm{p}_{[-1,1]}=\sup_{x\in [-1,1]}|p(x)|\le 1$,
$$\left|{p'(x)\over\sqrt{1-p^2(x)}}\right|\le(\deg p){1\over\sqrt{1-x^2}}, \ x\in (-1,1);$$
and, for a uniform estimate,
$$\norm{p'}_{[-1,1]}\leq (\deg p)^2\norm{p}_{[-1,1]}.$$
Equivalently, for a trigonometric polynomial $t=t(\theta)$ on the unit circle $T$,
$$\sup_{\theta} |t'(\theta)|\leq (\deg t)\sup_{\theta} |t(\theta)|. \eqno(32)$$
These estimates are useful in inverse theorems in univariate approximation theory. More generally, let $K$ be a compact set in $\CC^N$. We say that \dword{$K$ satisfies a Markov inequality with exponent $r$} if there exist constants $r\geq 1$ and $M>0$ depending only on $K$ such that
$$
\norm{\frac {\partial p}{\partial z_j}}_K \leq M (\deg p)^r \norm{p}_K, \quad j=1,\ldots,N$$
for all polynomials $p$. Convex sets in $\RR^N\subset \CC^N$ satisfy a Markov inequality with $r=2$ (cf.\ [BP]) while a closed Euclidean ball in $\CC^N$ satisfies a Markov inequality with exponent $r=1$. This last statement follows from the fact that (recall section 5) $V_K(z)=\max [0,\log |z-a|/R]$ if $K=\{z:|z-a|\leq R\}$ so that $V_K=V_K^*$ is Lipschitz, together with the following observation:

\proclaim Proposition. \ Let $K\subset \CC^N$ satisfy a H\"older continuity property (HCP):
$$\exp [V_K(z)]\leq 1+M\delta^m \ \hbox{if dist}(z,K)\leq \delta \leq 1$$ where $m,M>0$ are independent of $\delta>0$. Then $K$ satisfies a Markov inequality with exponent $r=1/m$.

\noindent {\bf Proof}. Fix a polynomial $p$ of degree $n$, say, and let $\norm{\frac {\partial p}{\partial z_j}}_K=|\frac {\partial p}{\partial z_j}(a)|$. Applying Cauchy's inequalities on a polydisk $P$ centered at $a\in K$ of (poly-)radius $R>0$ (recall the Cauchy integral formula (3)), we have
$$\norm{\frac {\partial p}{\partial z_j}}_K\leq \norm{p}_P/R.$$
From the Bernstein-Walsh inequality (19) we have
$$\norm{p}_P\leq \norm{p}_K\exp [n \sup_PV_K]\leq \norm{p}_K (1+MR^m)^n.$$
Thus
$$\norm{\frac {\partial p}{\partial z_j}}_K\leq \norm{p}_K \frac{(1+MR^m)^n}{R}.$$
Choosing $R= 1/n^{1/m}$ gives the result. \eop

To this date, {\sl there are no known examples of compact sets in $\CC^N$ which satisfy a Markov inequality but which do not satisfy (HCP)}.

One of the most beautiful applications of multivariate Markov inequalities is due to Plesniak [Pl1]. Recall that a $C^{\infty}$ function on a compact set $E\subset \RR^N$ is a function $f:E\to \RR$ such that there exists $\tilde f\in C^{\infty}(\RR^N)$ with $\tilde f|_E =f$. We write $f\in C^{\infty}(E)$. We say that $E$ is \dword{$C^{\infty}$-determining} if $g\in C^{\infty}(\RR^N)$ with $g|_E=0$ implies $D^{\alpha}g|_E =0$ for all multiindices $\alpha$. Plesniak [Pl1] has shown the following.

\proclaim Theorem ([Pl1]). \
Let $E\subset \RR^N$ be $C^{\infty}$-determining. Then $E$ satisfies a Markov inequality if and only if there is a continuous linear extension operator
$$L:(C^{\infty}(E),\tau_1)\to (C^{\infty}(\RR^N),\tau_0)$$
such that $L(f)|_E =f$ for each $f\in C^{\infty}(E)$.
\par

Here $\tau_0$ is the standard Fr\'echet space topology on $C^{\infty}(\RR^N)$ generated by the seminorms $\{\norm{\tilde f}_K^d:=\max_{|\alpha|\leq d}\norm{D^{\alpha}\tilde f}_K\}$ where $K$ ranges over compact subsets of $\RR^N$ and $d=0,1,\ldots$; and $\tau_1$ is the quotient topology on $C^{\infty}(\RR^N)/I(E)$ where $I(E):=\{f\in C^{\infty}(\RR^N):f|_E=0\}$. This is proved in [Pl1] using Lagrange interpolation operators corresponding to Fekete points (see section 9): the operator $L$ is of the form
$$L(f):=u_1L_1(f)+ \sum_{d=1}^{\infty}u_d(L_{d+1}(f)-L_d(f))$$
where $u_d$ are standard cut-off functions and $L_d(f)$ is the Lagrange interpolating polynomial of $f$ at a set of $d$-Fekete points of $E$.

There is an extensive literature on Markov inequalities. Baran and Plesniak and their students have produced many of the results related to multivariate approximation; cf.\ [Ba1], [Ba2], [Ba3], [Ba4]; and see  [Pl2] for a nice survey article. Markov inequalities have been used to construct natural pseudodistances on compact subsets of $\RR^N$ (see [BLW]). As a final application, observe that (32) can be interpreted in the following manner: setting $x=\cos \theta$ and $y=\sin \theta$, for any bivariate polynomial $p(x,y)$ on $\RR^2$, the unit tangential derivative $D_{\tau}p(x,y)$ on the unit circle $T\subset \RR^2$ satisfies
$$|D_{\tau}p(x,y)|_T\leq (\deg p) \norm{p}_T,  \quad (x,y)\in T.$$
Of course no estimate on normal derivatives is possible as there exist nonzero polynomials (e.g., $p(x,y)=x^2+y^2-1$) which vanish on $T$. Note that $T$ is an algebraic submanifold of $\RR^2$. Using a deep result of Sadullaev [Sa3] on the $L$-extremal function of compact subsets of algebraic sets in $\CC^N$, the following characterization of algebraicity is known:

\proclaim Theorem ([BLMT]). \ Let $K$ be a smooth, $m$-dimensional
submanifold of $\RR^N$ without boundary where $1\leq m \leq N-1$.
Then $K$ is algebraic if and only if $K$ satisfies a tangential
Markov inequality with exponent one: there exists a positive
constant $M$ depending only on $K$ such that for all polynomials
$p$ and all unit tangential derivatives $D_{\tau}$,
$$|D_{\tau}p(x_1,\ldots,x_N)|\leq M(\deg p)\norm{p}_K, \ \ \ (x_1,\ldots,x_N)\in K.$$
\par

Note that the finite-dimensionality of the vector space of polynomials of degree at most $n$ implies that there is a constant $C_n$ depending on $n$ and $K$ with
$$|D_{\tau}p(x_1,\ldots,x_N)|\leq C_n\norm{p}_K, \ \ \ (x_1,\ldots,x_N)\in K$$
for all such polynomials $p$. The content of the above theorem is that one can take $C_n=Mn$ where $M$ depends only on $K$. Indeed, a stronger version of the ``if'' implication is known: {\sl if $K$ satisfies a tangential Markov inequality
$$|D_{\tau}p(x_1,\ldots,x_N)|\leq M(\deg p)^r\norm{p}_K, \ \ \ (x_1,\ldots,x_N)\in K$$
with exponent $r< (m+1)/m$, then $K$ is algebraic}. We refer the reader to [BLMT] for details.

\sect{Appendix on pluripolar sets and extremal psh functions} In CCV, polar sets play an essential role. A subset $E\subset \CC$ is \dword{polar} if there exists
a {\it subharmonic} function $u$ defined in a neighborhood of $E$ with $E\subset \{z:u(z)=-\infty\}$; whereas a subset $E\subset \CC^N$ is \dword{pluripolar} if there exists a {\it plurisubharmonic} function $u$ defined in a neighborhood of $E$ with $E\subset \{z:u(z)=-\infty\}$. The neighborhood may be taken to be all of $\CC^N$. Apriori, there is a local notion in each case: $E$ is \dword{locally (pluri-)polar} if for each point $z\in E$ there exists an open neighborhood $U$ of $z$ and a (pluri-)subharmonic function $u$ in $U$ such that
$$E\cap U \subset \{z\in U: u(z)=-\infty\}.$$
It is easy if $N=1$ and much harder if $N>1$ to verify that the local notions are equivalent to the global ones. For $N>1$ this was first proved by Josefson [J].  We remark that since the notion of psh function makes sense on a complex manifold $M$ (see section 3), the notion of a locally pluripolar set in $M$ can be defined.
\vskip4pt
\item {1.} {\sl Nonpluripolar sets can be small}: Take a non-polar Cantor
set $E \subset \RR\subset \CC$ of Hausdorff dimension $0$ (for the idea behind the construction of such sets, see [Ra] section 5.3). Then
$E\times \cdots \times E$ is nonpluripolar in $\CC^N$ (in general,
$E_1 \times \cdots \times E_j\subset \CC^{m_1}\times  \cdots \times
\CC^{m_j}$ is nonpluripolar in $\CC^{m_1+\cdots +m_j}$ if and
only if $E_k\subset \CC^{m_k}$ is nonpluripolar in $ \CC^{m_k}$ for $k=1,\ldots,j$) and
has Hausdorff dimension $0$.
\item {2.} {\sl Pluripolar sets can be big}: A complex hypersurface
$S=\{z:f(z)=0\}$ associated to a holomorphic function $f$ is a pluripolar set (take
$u=\log {|f|}$) which has Hausdorff dimension $2N-2$. Recall that a psh
function is, in particular, subharmonic in the $\RR^{2N}$ sense;
hence a pluripolar set is Newtonian polar and thus the Hausdorff
dimension of a pluripolar set cannot exceed $2N-2$ (cf.\ [Ca], section IV).
\item {3.} {\sl Size doesn't matter}: In $\CC^2$, the totally real
plane $\RR^2=\{(z_1,z_2): \Im z_1 = \Im z_2=0\}$ is nonpluripolar
(why?) but the {\it complex} plane $\CC= \{(z_1,0): z_1 \in {\bf
C}\}$ is pluripolar (take $u=\log {|z_1|}$). Also, there exist
$C^{\infty}$ arcs in $\CC^N$ which are not pluripolar (cf.\ [DF]); while such a
{\it real-analytic} arc must be pluripolar (why?).
\vskip4pt

One can easily construct examples of nonpluripolar sets $E\subset \CC^N$ which intersect every affine complex line in finitely many points (hence these intersections are polar in these lines). Indeed, take
$$E:=\{(z_1,z_2)\in \CC^2: \Im (z_1+z_2^2)=\Re (z_1+z_2+z_2^2)=0\}.$$
Then for any complex line $L:=\{(z_1,z_2):a_1z_1 +a_2z_2 =b\}, \ a_1,a_2,b\in \CC$, $E\cap L$ is the intersection of two real quadrics and hence consists of at most four points. However, $E$ is a totally real, two-(real)-dimensional submanifold of $\CC^2$ and hence -- as is the case with $\RR^2=\RR^2+i0 \subset \CC^2$ in 3. -- is not pluripolar. Thus pluripolarity cannot be detected by ``slicing'' with complex lines. In this example, $E$ intersects the one-(complex)-dimensional analytic variety $A:=\{(z_1,z_2):z_1+z_2^2=0\}$ in a nonpolar set. Nevertheless, one can construct a nonpluripolar set $E$ in $\CC^N, \ N>1$, which intersects every one-dimensional complex analytic subvariety in a polar set [CLP].

In certain instances, however, slicing can detect pluripolarity. We define a set $E\subset \CC^N$ to be \dword{pseudoconcave} if for each point $p\in E$ there is a neighborhood $U$ of $p$ such that $U\setminus E$ is open and pseudoconvex in $\CC^N$ (see section 3). Canonical examples are zero sets of a holomorphic function, or, more generally, zero sets of multiple-valued holomorphic functions; e.g., $\{(z_1,z_2)\in \CC^2: z_2^2=z_1\}$ is pseudoconcave. This notion is related to the work in section 11 with the class $R^0$. Sadullaev [Sa1] has shown the remarkable result that if $E$ is a closed, pseudoconcave set in $\CC^N\setminus \{0\}$, then $E$ is pluripolar if and only if $E\cap L$ is polar in $L$ for each complex line $L$ passing through $0$. Moreover, for this class of pseudoconcave sets, pluripolarity is equivalent to $\RR^{2N}$-(Newtonian) polarity!

There are several distinct notions of capacities in SCV. Given a compact set $K\subset \CC$, we recall the definition of the extremal psh function $V_K^*(z)$, the usc regularization of
$$
  V_K(z)
  := \max \left\{ 0 ,
    \sup _p \left\{ {1 \over \deg p} \log|p(z)| \right\}
    \right\}.$$
The \dword{Siciak} or \dword{Robin} capacity of $K$ is the number
$$c(K):=\exp\bigl (-\limsup_{|z|\to \infty} [V_K^*(z) - \log {|z|}] \bigr).$$
Unlike in CCV, the limit (usually) does not exist. Indeed, the Robin function $\rho_K :=\rho_{V_K^*}$ (see section 5) associated to $V_K^*$ provides information on the asymptotic behavior of this function on complex lines through the origin.

The quantity $c(K)$ coincides with a Chebyshev constant $\tilde T(K):=\lim_{n\to \infty} \tilde M_n(K)^{1/n}$ where
$$\tilde M_n(K):= \inf \{\norm{p_n}_K: p_n=\hat p_n+ \ \hbox{lower degree terms} \ , \ \norm{\hat p_n}_{\bar B}\geq 1\}$$
([Sa4], section 10). Here $\bar B$ is the closed unit (Euclidean) ball. Other normalizations may be used to define other Chebyshev constants. For example, defining
$$M_n(K):= \inf \{\norm{p_n}_K: \ \norm{p_n}_{\bar B}\geq 1\},$$
the Chebyshev constant $T(K):=\lim_{n\to \infty} M_n(K)^{1/n}$ coincides with $\exp \bigl (-\sup_{z\in B }V_K^*(z)\bigr)$. Although the numbers $T(K)$ and $\tilde T(K)$ are, in general, different, they are comparable; i.e. for all compact sets $K$, $T(K)$ is bounded above and below by a constant multiple of $\tilde T(K)$. Note that the subsets of ${\cal P}_n$ used to define $M_n(K)$ and $\tilde M_n(K)$ are not multiplicative classes (as in the case of univariate monic polynomials); thus an (elementary) argument is needed to verify the existence of the limits $T(K)$ and $\tilde T(K)$.

Next we give the definition of the transfinite diameter $d(K)$. Details may be found in Zaharjuta's paper [Z1]. Let $e_1(z),\ldots,e_j(z),\ldots$ be a listing of the monomials
$\{e_i(z)=z^{\alpha(i)}=z_1^{\alpha_1}\cdots z_N^{\alpha_N}\}$ in
$\CC^N$ indexed using a lexicographic ordering on the multiindices $\alpha(i)\in {\bf N}^N$, but with $\deg e_i=|\alpha(i)|$ nondecreasing. For
$\zeta_1,\ldots,\zeta_n\in \CC^N$, let
$$VDM(\zeta_1,\ldots,\zeta_n)=|\det [e_i(\zeta_j)]_{i,j=1,\ldots,n}| $$
and for a compact subset $K\subset \CC^N$ let
$$V_n =V_n(K):=\max_{\zeta_1,\ldots,\zeta_n\in K}VDM(\zeta_1,\ldots,\zeta_n).$$
Define
$h_d= \#\{i:\deg e_i\leq d\}$ and $l_d=\sum_{i=1}^{h_d}(\deg e_i)$. Then
$$d(K):= \limsup_{d\to \infty}V_{h_d}^{1/l_d} \eqno(33)$$
is the \dword{transfinite diameter} of $K$.

If $N=1$, it is well-known and trivial that the sequence $\{V_{h_d}^{1/l_d} \}$ is monotone decreasing and hence has a limit; moreover, in this case, $d(K)$ coincides with the logarithmic capacity of $K$ and the Chebyshev constant of $K$ defined back in section 5. Zaharjuta [Z1] proved the highly nontrivial result that the limit in (33) exists in the case when $N>1$. In this setting the numbers $c(K),T(K)$ and $d(K)$ are not generally equal; however it is the case that {\sl for $K\subset \CC^N$ compact, $K$ is pluripolar if and only if $c(K)=T(K)=d(K)=0$} (see [LT]).

There are several extremal psh functions in SCV. Recall the relative
extremal function introduced at the end of section 6:
for $E$ a subset of $D$, define
$$\omega(z,E,D):=\sup \{u(z): u \ \hbox{psh in} \ D, \  u \leq
0 \ \hbox{in} \ D, \ u|_E \leq -1\}.$$
The usc regularization $\omega^*(z,E,D)$ is called the \dword{relative
extremal function} of $E$ relative to $D$.

\proclaim Proposition ($\omega$). \ Either $\omega^*\equiv 0$ in
$D$ or else $\omega^*$ is a nonconstant psh function in $D$. We have $\omega^*\equiv
0$ if and only if $E$ is pluripolar.

\noindent {\bf Proof}. If $\omega^*(z^0)=0$ at some point $z^0\in D$, then
$\omega^*\equiv 0$ in $D$ by the maximum principle. Hence we can find a
sequence $z^j \to z^0, \ z^j \in D$, with $\omega(z^j,E,D)\to 0$. By
subaveraging,
$$\omega(z^j,E,D)\leq {1\over {\rm
vol}(B(z^j,r))}\int_{B(z^j,r)}\omega(z,E,D)dA(z)$$
for $r$ sufficiently small so that $B(z^j,r)\subset D$. We conclude that $\omega(z,E,D)=0$ a.e.\ in a neighborhood
of
$z^0$. Fix a point
$z'$ with $\omega(z',E,D)=0$ and take a sequence of psh functions $u_j$ in $D$ with $u_j \leq
0$ in $D$, $u_j|_E \leq -1$, and $u_j(z')\geq -1/2^j$. Then $u(z):=\sum
u_j(z)$ is psh in $D$ (the partial sums form a decreasing sequence of psh
functions) with $u(z')\geq -1$ (so $u\not \equiv -\infty$) and $u|_E=
-\infty$; thus $E$ is pluripolar.

Conversely, if $E$ is pluripolar, there exists $u$ psh in $D$ with
$u|_E=-\infty$; since $D$ is bounded we may assume $u\leq 0$ in $D$. Then
$\epsilon u\leq \omega(z,E,D)$ in $D$ for all $\epsilon >0$ which implies
that $\omega(z,E,D)=0$ at all points $z\in D$ where $u(z)\not = -\infty$. Since
pluripolar sets have measure zero (why?), $\omega(z,E,D)=0$ a.e.\ in $D$
and hence $\omega^*(z,E,D)\equiv 0$ in $D$. \eop

Using Proposition ($\omega$), Bedford and Taylor [BT2] gave a simple proof of Josefson's result that locally pluripolar sets are globally pluripolar. Similar to this proposition, one can show that for a bounded set $E\subset \CC^N$, $V_E^*\equiv +\infty$ if and only if $E$ is pluripolar. As mentioned in section 1, a nice introduction to pluripotential theory is the book of Klimek [K]. There is also a  developing theory of {\it weighted} pluripotential theory. We refer the reader to the book of Saff-Totik [SaT] for an introduction to one-variable weighted potential theory in $\CC$. An introduction to the SCV setting can be found in the appendix of [SaT] written by Tom Bloom (see also [BlL3]).

We remark that the quantity
$$C(E,D):=\sup \{\int_E (dd^cu)^N: u \ \hbox{psh in} \ D, \ 0\leq u \leq
1 \ \hbox{in} \ D\}$$
is called the \dword{relative capacity} of $E$ relative to $D$. The precise definition of this complex Monge-Amp\`ere operator $(dd^cu)^N$ will be given in the next section. Alexander and Taylor [AT] proved the following comparison between the relative capacity $C(K,D)$ and the Chebyshev constant $T(K)$ for a compact set $K$:

\proclaim Theorem ([AT]). \ Let $K$ be a compact subset of the
unit ball $B\subset \CC^N$. We have $$\exp
\bigl[-A/C(K,B)\bigr]\leq T(K) \leq
\exp\bigl[-(c_N/C(K,B))^{1/N}\bigr]$$ where the right-hand
inequality holds for all $K\subset B$ and the left-hand inequality
holds for all $K\subset B(0,r):=\{z:|z|<r\}$ where $r<1$ and
$A=A(r)$ is a constant depending only on $r$. \par

\noindent Here, $c_N$ is a dimensional constant. As a corollary, Alexander and Taylor construct normalized polynomials that are small where a given holomorphic function is small; or, more generally, where a psh function is very negative:

\proclaim Proposition ([AT]). \ Let $u$ be a negative psh function in the unit ball $B$ with $u(0)\geq -1$. For $r<1$ and $A>1$, let $K$ be a compact subset of
$$\{z\in \CC^N: |z|\leq r, \ u(z)< -A\}.$$
Then there exists a sequence $\{p_d\}$ of polynomials with deg$p_d\leq d$ and
$\norm{p_d}_B=1$ such that
$$\norm{p_d}_K\leq \exp \bigl[-C(A)^{1/N}d\bigr]$$
where $C$ depends only on $r$ and $N$. \par

\noindent The proposition is proved in a much more complicated manner in [J]; there, Pad\'e-type approximants are constructed. This is the key ingredient in the original proof of Josefson's theorem that locally pluripolar sets are globally pluripolar.

\sect{Appendix on complex Monge-Amp\`ere operator} For simplicity, we work in $\CC^2$ with variables $(z,w)$. We use the notation
$d=\partial +\bar \partial$ and $d^c =i(\bar \partial - \partial)$
where, for a $C^1$ function $u$,
$$\partial u:={\partial u\over \partial z}dz+{\partial u\over \partial w}dw,  \quad \bar \partial u:={\partial u\over \partial \bar z}d\bar z+ {\partial u\over \partial \bar w}d\bar w$$ (recall section 7) so that $dd^c =2i \partial \bar \partial$. For a $C^2$ function $u$,
$$(dd^cu)^2=16 \bigl[ {\partial^2 u \over \partial z  \partial \bar z}{\partial^2 u \over \partial w  \partial \bar w}-{\partial^2 u \over \partial z \partial \bar w}{\partial^2 u \over \partial w  \partial \bar z}\bigr]{i\over 2}dz\wedge d\bar z\wedge {i\over 2}dw\wedge d\bar w$$
is, up to a positive constant, the determinant of the complex
Hessian of $u$ times the volume form on $\CC^2$. Thus if $u$ is also
psh, $(dd^cu)^2$ is a positive measure which is absolutely
continuous with respect to Lebesgue measure. If $u$ is psh in an
open set $D$ and locally bounded there, then
$(dd^cu)^2$ is a positive measure in $D$ (cf.\ [BT1]).

To see this, we first recall that a psh function $u$ in $D$ is an usc function $u$ in $D$ which is
subharmonic on components of $D\cap L$ for complex affine lines $L$.
In particular, $u$ is a locally integrable function in $D$ such that
$$dd^cu=2i\bigl[{\partial^2 u \over \partial z
\partial \bar z}dz\wedge d\bar z+{\partial^2 u \over \partial w
\partial \bar w}dw\wedge d\bar w+{\partial^2 u \over \partial z
\partial \bar w}dz\wedge d\bar w+ {\partial^2 u \over \partial \bar
z \partial w}d\bar z\wedge dw\bigr] $$ is a positive $(1,1)$
current (dual to $(1,1)$ forms); i.e., a $(1,1)$ form with
distribution coefficients. The derivatives are
to be interpreted in the distribution sense and are actually measures; i.e., they act on compactly supported continuous functions. Here, a $(1,1)$ current
$T$ on a domain $D$ in $\CC^2$ is positive if $T$ applied to $i\beta \wedge
\bar \beta$ is a positive distribution for all $(1,0)$ forms $\beta
= adz+bdw$ with $a,b\in C^{\infty}_0(D)$ (smooth functions having
compact support in $D$). Writing the action of a current $T$ on a
form $\psi$ as $\inpro{T,\psi}$, this means that
$$\inpro{T, \phi (i\beta \wedge \bar \beta)}\geq 0 \quad \hbox{for all} \ \phi \in C^{\infty}_0(D) \ \hbox{with} \ \phi \geq 0.$$
For a discussion of currents and the general definition of
positivity, we refer the reader to Klimek [K], section 3.3.

Following [BT1], we now define $(dd^cv)^2$ for a psh $v$ in $D$
if $v\in L^{\infty}_{loc}(D)$ using the fact that $dd^cv$ is a
positive $(1,1)$ current with measure coefficients. First note that
if $v$ were of class $C^2$, given $\phi\in C^{\infty}_0(D)$, we have
$$\eqalign{\int_D \phi (dd^cv)^2&\;=\;  -\int_D d\phi \wedge d^c v \wedge dd^cv\cr
&\;= \hbox{(exercise!)}\ \ -\int_D dv \wedge d^c \phi \wedge dd^cv= \int_D v dd^c \phi \wedge dd^c v}$$
since all boundary integrals vanish. The applications of Stokes'
theorem are justified if $v$ is smooth; for arbitrary psh $v$ in $D$
with $v\in  L^{\infty}_{loc}(D)$, these formal calculations serve as
motivation to {\it define} $(dd^cv)^2$ as a positive measure
(precisely, a positive current of bidegree $(2,2)$ and hence a
positive measure) via
$$\inpro{(dd^cv)^2,\phi}:=\int_D v dd^c \phi \wedge dd^c v.$$
This defines $(dd^cv)^2$ as a $(2,2)$ current (acting on $(0,0)$
forms; i.e., test functions) since $vdd^cv$ has measure
coefficients. We refer the reader to [BT1] or [K] (p.~113)
for the verification of positivity of $(dd^cv)^2$.

The next result shows that for a nonpluripolar compact set $K\subset \CC^2$, the Monge-Amp\`ere measure $(dd^cV_K^*)^2$ associated to the $L$-extremal function $V_K^*$ of $K$ plays the role of the equilibrium measure $\Delta g_K$ associated to the Green function $g_K$ of a nonpolar compact set $K \subset \CC$. A plurisubharmonic function $u$ in a domain $D$ satisfying the property that {\sl for any $D'$ relatively compact in $D$, and any $v$ psh in $\bar D'$, if $u\geq v$ on $\partial D'$, then $u\geq v$ on $D'$}, is called \dword{maximal} in $D$. Bedford and Taylor showed that for locally bounded psh $u$, $u$ is maximal in $D$ if and only if $(dd^cu)^2=0$ in $D$. Thus, maximal psh functions are the ``correct'' analogue of harmonic functions in $\CC$. The following result (cf. [BT2], Corollary 9.4) shows that $V_K^*$ is maximal outside of $\hat K$.

\proclaim Proposition. \ Let $K$ be a nonpluripolar compact set in $\CC^N$ . Then we have $(dd^cV_K^*)^N=0$ outside of $K$.

Similarly, if $D$ is a bounded open neighborhood of $K$, the relative extremal function satisfies $(dd^c\omega^*(\cdot,K,D))^N=0$ in $D\setminus K$.

\sect{A few open problems} Here are a few open problems.
\item {1.} (Section 6) {\sl Theorem ([A], [Bl2]) has a sharp version for balls; e.g., if $K$ is the closed unit ball of $\RR^N$, then $f\in C(K)$ extends to be harmonic in the ball of radius $R>1$ if and only if $\limsup_{n\to \infty} d_n(f,K)^{1/n} \leq 1/R$ (cf.\ [BL1]). For $N>2$, are there compact sets other than balls for which a sharp version of this theorem holds?}

\item {2.} (Section 7) {\sl Let $D$ be a smoothly bounded pseudoconvex domain. Does $D$ possess the Mergelyan property if and only if there are pseudoconvex domains $D_j$ with $\bar D \subset D_j$ such that $\bar D = \cap_j D_j$?}

\item {3.} (Section 8) {\sl Does there exist a compact, polynomially convex subset $K$ of the unit sphere in $\CC^2$ such that $P(K)\not= C(K)$?}

\item {4.} (Section 9) {\sl For Fekete arrays on a nonpluripolar compact set $K\subset \CC^N, \ N\geq 2$, do the normalized discrete measures $\mu_n$ converge weak-* to the Monge-Amp\'ere measure $\mu_K :=(dd^cV_K^*)^N$ of the $L$-extremal function $V_K^*$? Verify this for {\bf any} single nonpluripolar compact set $K$! More generally, is the conclusion true for arrays satisfying} $$\limsup_{n\to \infty}\Lambda_n^{1/n}\leq 1?$$

\item {5.} (Section 10) {\sl Let $K\subset \RR^N, \ N\geq 3$, be a compact, convex set with nonempty interior. Which such sets $K$ admit \dword{harmonic} extremal arrays, i.e., arrays $\{{\cal A}_d\}_{d=1,2,\ldots}$ in $K$ such that  ${\cal K}_d(f)$ converges to $f$ uniformly on $K$ for each $f$ harmonic in a neighborhood of $K$?}

\item {6.} (Section 12) {\sl If $K\subset \CC^N$ satisfies a Markov inequality does $K$ have property (HCP)? Is $K$ necessarily regular? If $N>1$, is $K$ necessarily nonpluripolar?}




\References
\def\AJM{Amer.\ J. Math.}
\def\CA{Constr.\ Approx.}
\def\IUMJ{Indiana Univ.\ Math.\ J.}
\def\JAT{J. Approx.\ Theory}
\def\MA{Math.\ Ann.}
\def\MSNS{Mat.\ Sb.\ (N.S.)}
\def\PAMS{Proc.\ Amer.\ Math.\ Soc.}
\def\SJMA{SIAM J. Math.\ Anal.}

\parindent1.5em
\def\nextitem#1{\item{#1}}
\def\bibitem{}
\nextitem {[AT]}
\refJ  Alexander, H., Taylor, B. A.;
Comparison of two capacities in $\CC^n$;
{\sl Math. Z.}; 186(3); 1984; 407--417; MR0744831 (85k:32034)

\nextitem {[AW]}
\refB Alexander, H., Wermer, J.;
Several complex variables and Banach algebras, Third edition;
Graduate Texts in Mathematics, 35. Springer-Verlag (New York); 1998; xii+253 pp.; MR1482798 (98g:32002)

\nextitem {[AIW1]}
\refJ Anderson, J., Izzo, A., Wermer, J.;
Polynomial approximation on three-dimensional real-analytic submanifolds of $\CC^n$;
\PAMS; 129(8); 2001; 2395--2402; MR1823924 (2002d:32021)

\nextitem {[AIW2]}
\refJ Anderson, J., Izzo, A., Wermer, J.;
Polynomial approximation on real-analytic varieties in $\CC^n$;
\PAMS; 132(5); 2004; 1495--1500; \hfill\break MR2053357 (2005d:32017)

\nextitem {[AP]}
\refJ Andersson, M., Passare, M.;
Complex Kergin interpolation;
\JAT; 64(3); 1991; 214--225; MR1091471 (92b:41003).

\nextitem {[APS]}
\refB Andersson, M., Passare, M., Sigurdsson, R.;
Complex convexity and analytic functionals;
Progress in Mathematics, 225. BirkhŠuser Verlag, (Basel); 2004;
xii+160 pp.; ISBN: 3-7643-2420-1, MR2060426 (2005a:32011)

\nextitem {[A]}
\refJ Andrievskii, V;
Uniform harmonic approximation on compact sets in $\RR^k,\;k\geq 3$;
\SJMA; 24(3); 1993; 216--222; MR1199535 (93m:41010)

\nextitem {[BBL1]}
\refJ Bagby, T., Bos, L., Levenberg, N.;
Quantitative approximation theorems for elliptic operators;
\JAT; 85(3); 1996; 69--87; MR1382051 (97h:41060)

\nextitem {[BBL2]}
\refJ Bagby, T., Bos, L., Levenberg, N.;
Multivariate simultaneous approximation;
\CA; 18(3); 2002; 569--577; MR1920286 (2003f:41026)

\nextitem {[BL1]}
\refJ Bagby, T., Levenberg, N.;
Bernstein theorems;
New Zealand J. Math.; 22(3); 1993; 1--20; MR1244005 (95j:41014)

\nextitem {[BL2]}
\refQ Bagby, T., Levenberg, N.;
Bernstein theorems for harmonic functions;
(Methods of approximation theory in complex analysis and mathematical physics (Leningrad, 1991)),
xxx (ed.), {\sl Lecture Notes in Math.}, {\bf 1550}, Springer (Berlin); 1993; 7--18; MR1322287 (95m:41009)

\nextitem {[BL3]}
\refJ Bagby, T., Levenberg, N.;
Bernstein theorems for elliptic equations;
\JAT; 78(3); 1994; 190--212; MR1285258 (96b:41036)

\nextitem {[Ba1]}
\refJ Baran, M.;
Bernstein type theorems for compact sets in $\RR^n$;
\JAT; 69(3); 1992; 156--166; MR1160251 (93e:41021)

\nextitem {[Ba2]}
\refJ Baran, M.;
Bernstein type theorems for compact sets in $\RR^n$ revisited;
\JAT; 79(3); 1994; 190--198;  MR1302342 (95h:41023)

\nextitem {[Ba3]}
\refJ Baran, M.;
Complex equilibrium measure and Bernstein type theorems for compact sets in $\RR^n$;
\PAMS; 123(3); 1995; 485--494; MR1219719 (95c:31006)

\nextitem {[Ba4]}
\refJ Baran, M.;
Conjugate norms in $\CC^n$ and related geometrical problems;
{\sl Dissertationes Math. (Rozprawy Mat.)}; 377; 1998; 67 pp; MR1657095 (2000c:32086)

\nextitem {[BP]}
\refJ Baran, M., Ple\'sniak, W.;
Markov's exponent of compact sets in $\CC^n$;
\PAMS; 123(3); 1995; 2785--2791; MR1301486 (95k:41022)

\nextitem {[BF]}
\refQ Bedford, E., Fornaess, J.-E.;
Approximation on pseudoconvex domains;
(Complex approximation {(Proc.\ Conf., Quebec, 1978)}), xxx (ed.),
{\sl Progr.\ Math.}, 4, Birkh\"auser, Boston (Mass.); 1980; 18--31; MR0578636 (82b:32019)

\nextitem {[BT1]}
\refJ Bedford, E., Taylor, B. A.;
The Dirichlet problem for a complex Monge-Amp\'ere equation;
{\sl Invent.\ Math.}; 37(3); 1976; 1--44; MR0445006 (56 \#3351)

\nextitem {[BT2]}
\refJ Bedford, E., Taylor, B. A.;
A new capacity for plurisubharmonic functions;
{\sl Acta Math.}; 149; 1982; 1--40; MR0674165 (84d:32024)

\nextitem {[Ber]} \refR Berndtsson, B.; A remark on approximation
on totally real sets; arXiv:math.CV/0608058; 2006;

\nextitem {[Bi]}
\refJ Bishop, E.;
Mappings of partially analytic spaces;
\AJM; 83; 1961; 209--242; MR0123732 (23 \#A1054)

\nextitem {[Bl1]}
\refJ Bloom, T.;
Polynomial interpolation;
{\sl Bol.\ Soc.\ Brasil.\ Mat.}; 10(3); 1979; 75--86;
\hfill\break MR0607007 (82h:41001)

\nextitem {[Bl2]}
\refJ Bloom, T.; Kergin interpolation of entire functions on $\CC^n$;
{\sl Duke Math.\ J.}; 48(3); 1981; 69--83; MR0610176 (83k:32005)

\nextitem {[Bl3]}
\refJ Bloom, T.;
On the convergence of multivariable Lagrange interpolants;
\CA;  5(3); 1989; 415--435; MR1014307 (90m:32032)

\nextitem{[Bl4]}
\refJ Bloom, T.;
Interpolation at discrete subsets of $\CC^n$;
\IUMJ; 39(3); 1990; 1223--1243; MR1087190 (91k:32015)

\nextitem {[Bl5]}
\refJ Bloom, T.;
Orthogonal polynomials in $\CC^ n$;
\IUMJ;  46(3); 1997; 427--452;

MR1481598 (98j:32006)

\nextitem {[Bl6]}
\refJ Bloom, T.;
Some applications of the Robin function to multivariable approximation theory;
\JAT; 92(3); 1998; 1--21; MR1492855 (98k:32021)

\nextitem {[Bl7]}
\refJ  Bloom, T.;
On the convergence in capacity of rational approximants;
\CA; 17(3); 2001; 91--102; MR1794803 (2001j:41012)

\nextitem {[Bl8]} \refJ Bloom, T.;
On families of polynomials which approximate
the pluricomplex Green function;
\IUMJ; 50(4); 2001; 1545--1566; MR1889070  (2003a:32055)

\nextitem {[BBCL]}
\refJ Bloom, T., Bos, L., Christensen, C., Levenberg, N.;
Polynomial interpolation of holomorphic functions in $\CC$ and $\CC^ n$;
{\sl Rocky Mountain J. Math.}; 22(3); 1992; 441--470; MR1180711 (93i:32016)

\nextitem {[BC1]}
\refJ Bloom, T., Calvi, J.-P.;
A continuity property of multivariate Lagrange interpolation;
{\sl Math.\ Comp.}; 66(220); 1997; 1561--1577; MR1422785 (98a:41001)

\nextitem {[BC2]}
\refJ Bloom, T., Calvi, J.-P.;
Kergin interpolants of holomorphic functions;
\CA; 13(3); 1997; 569--583; MR1466066 (98h:32021)

\nextitem {[BC3]}
\refJ Bloom, T., Calvi, J.-P.;
The distribution of extremal points for Kergin interpolation: real case;
{\sl Ann.\ Inst.\ Fourier} (Grenoble); 48(3); 1998; 205--222; MR1614898 (99c:32015)

\nextitem {[BlL1]}
\refJ Bloom, T., Levenberg, N.;
Lagrange interpolation of entire functions in $\CC^ 2$;
{\sl New Zealand J. Math.};  22(3); 1993; 65--73; MR1244023 (94i:32016)

\nextitem {[BlL2]}
\refJ Bloom, T., Levenberg, N.;
Distribution of nodes on algebraic curves in $\CC^N$;
{\sl Ann.\ Inst.\ Fourier} (Grenoble); 53(3); 2003; 1365--1385;
MR2032937 (2004j:32035)

\nextitem {[BlL3]}
\refJ Bloom, T., Levenberg, N.;
Weighted pluripotential theory in $\CC^N$;
\AJM; 125(3); 2003; 57--103; MR1953518 (2003k:32045)

\nextitem {[Bo2]}
\refJ Bos, L.; On Kergin interpolation in the disk;
\JAT; 37(3); 1983; 251--261;
\hfill\break MR0693012 (85b:41001)

\nextitem {[BCDVX]} \refJ  Bos, L.,  Caliari, M., De Marchi, S.,
Vianello, M., Xu, Y.; Bivariate Lagrange interpolation at the
Padua points: the generating curve approach; \JAT; xxx; 200x;
xxx--xxx;

\nextitem {[BC]}
\refJ Bos, L., Calvi, J.-P.;
Kergin interpolants at the roots of unity approximate $C^2$ functions;
{\sl J. Anal.\ Math.}; 72; 1997; 203; MR1482995 (98j:41001)

\nextitem {[BLMT]}
\refJ Bos, L., Levenberg, N., Milman, P., Taylor, B. A.;
Tangential Markov inequalities characterize algebraic submanifolds of $\RR^N$;
\IUMJ; 44(3); 1995; 115--138; MR1336434 (96i:41009)

\nextitem {[BLW]} \refR Bos, L., Levenberg, N., Waldron, S.;
Metrics associated to multivariate polynomial inequalities II;
preprint; 2006;

\nextitem {[BW]}
\refQ Bos, L., Waldron, S.;
On the structure of Kergin interpolation for points in general position;
(Recent progress in multivariate approximation {(Witten-Bommerholz, 2000)}),
xxx (ed.),  {\sl Internat.\ Ser.\ Numer.\ Math.}, {\bf 137}, Birkh\"auser (Basel); 2001; 75--87; MR1877498 (2003h:41041)

\nextitem {[BB]}
\refJ Bruna, J., Burgu\'es, J..;
Holomorphic approximation in $C^m$-norms on totally real compact sets in $\CC^n$;
\MA; 269(3); 1984; 103--117; MR0756779 (86c:32014)

\nextitem {[Ca]}
\refB  Carleson, L.;
Selected problems on exceptional sets; Van Nostrand Mathematical Studies, No.~13 D. Van Nostrand Co., Inc. (Princeton, N.J.-Toronto, Ont.-London); 1967; v+151 pp., MR0225986 (37 \#1576)

\nextitem {[C1]}
\refJ Chirka, E. M.;
Meromorphic continuation, and the rate of rational approximations in $\CC^N$. (Russian);
\MSNS; 99(141); 1976; 615--625; MR0412472 (54 \#598)

\nextitem {[C2]} \refJ Chirka, E. M.; Rational approximations of
holomorphic functions with singularities of finite order.
(Russian); \MSNS; 100(142); 1976; 137--155; MR0417423 (54 \#5473)

\nextitem {[CLP]}
\refJ Coman, D., Levenberg, N., Poletsky, E.;
Smooth submanifolds intersecting any analytic curve in a discrete set;
\MA; 332(3); 2005; 55--65; MR2139250 (2005m:32066)

\nextitem {[DF]}
\refJ Diederich, K., Fornaess, J.-E.;
A smooth curve in $\CC^{2}$ which is not a pluripolar set;
{\sl Duke Math. J.}; 49(4); 1982; 931--936; MR0683008 (85b:32025)

\nextitem {[DL]}
\refQ Duval, J., Levenberg, N.;
Large polynomial hulls with no analytic structure;
(Complex analysis and geometry {(Trento, 1995)}), xxx (ed.),
 {\sl Pitman Res.\ Notes Math.} Ser., 366, Longman (Harlow); 1997; 119--122;
MR1477444 (99g:32023)

\nextitem {[G1]}
\refJ Gonchar, A. A.;
A local condition for the single-valuedness of analytic functions of several variables. (Russian);
\MSNS; 93(135); 1974; 296--313, 327; MR0589892 (58 \# 28632)

\nextitem {[G2]}
\refJ Gonchar, A. A.;
On a theorem of Saff. (Russian);
\MSNS; 94(136); 1975; 152--157; MR0396965 (53 \#825)

\nextitem {[GMS]}
\refJ G\"otz, M., Maymeskul, V. V., Saff, E. B.;
Asymptotic distribution of nodes for near-optimal polynomial interpolation on certain curves in $\RR^2$;
\CA; 18(3); 2002; 255--283; MR1890499 (2003e:30007)

\nextitem {[HW]}
\refJ Harvey, F. R., Wells, R. O.;
Holomorphic approximation and hyperfunction theory on a $C^{1}$\ totally real submanifold of a complex manifold;
{\sl Math. Ann.}; 197; 1972; 287--318; MR0310278 (46 \#9379)

\nextitem {[He]} G. M. Henkin [1985] The method of integral
representations in complex analysis. Current problems in
mathematics. Fundamental directions, Vol. 7, in a  translation of
{\sl Sovremennye problemy matematiki. Fundamentalnye napravleniya,
Tom 7}, Akad.\ Nauk SSSR, Vsesoyuz.\ Inst.\ Nauchn.\ i Tekhn.\
Inform., Moscow [MR0850489 (87f:32003)]. Translation by P. M.
Gauthier. Translation edited by A. G. Vitushkin. Encyclopaedia of
Mathematical Sciences, 7. Springer-Verlag (Berlin); 1990; vi+248
pp. ISBN: 3-540-17004-9, MR1043689 (90j:32003)

\nextitem {[H\"o]}
\refB H\"ormander, L.;
An introduction to complex analysis in several variables;
Third edition. North-Holland Mathematical Library, 7. North-Holland Publishing Co. (Amsterdam);
1990; xii+254 pp. ISBN: 0-444-88446-7, MR1045639 (91a:32001)

\nextitem {[H\"oW]}
\refJ H\"ormander, L.,  Wermer, J.;
Uniform approximation on compact sets in $\CC^n$;
{\sl Math.\ Scand.}; 23; 1968; 5--21 (1969); MR0254275 (40 \#7484)

\nextitem {[Iv]}
\refQ Ivashkovich, S. M.;
Rational curves and extensions of holomorphic mappings;
(Several complex variables and complex geometry, Part 1 {(Santa Cruz, CA, 1989)}), xxx (ed.),
{\sl Proc.\ Sympos.\ Pure Math.}, {\bf 52}, Part 1, Amer.\ Math.\ Soc.  (Providence, RI); 1991;  93--104; MR1128517 (92k:32022)

\nextitem {[I]}
\refJ Izzo, A.;
Failure of polynomial approximation on polynomially convex subsets of the sphere;
{\sl Bull. London Math.\ Soc.}; 28(3); 1996; 393--397; MR1384828 (98d:32017)

\nextitem {[JP]}
\refB Jarnicki, M., Pflug, P.;
Extension of holomorphic functions;
de Gruyter Expositions in Mathematics, 34. Walter de Gruyter \& Co.  (Berlin);
2000; x+487 pp. ISBN: 3-11-015363-7, MR1797263 (2001k:32017)

\nextitem {[J]}
\refJ Josefson, B.;
On the equivalence between locally polar and globally polar sets for plurisubharmonic functions on $\CC^n$;
{\sl Ark. Mat.}; 16(1); 1978; 109--115; MR0590078 (58 \#28669)

\nextitem {[Ka]}
\refQ Kallin, E.;
Polynomial convexity: The three spheres problem;
(1965 Proc.\ Conf.\ Complex Analysis {(Minneapolis, 1964)}), xxx (ed.),
Springer (Berlin); 1965; 301--304; MR0179383 (31 \#3631)

\nextitem {[Kh]} \refJ Khudauiberganov, G.; On the polynomial and
rational convexity of the union of compact sets in $\CC^n$.
(Russian); Izv.\ Vyssh.\ Uchebn.\ Zaved.\ Mat.; 297(2); 1987;
70--74; MR0889199 (88g:32031)

\nextitem {[K]}
\refB Klimek, M.;  Pluripotential theory;
London Math.\ Society Monographs. New Series, 6. Oxford Science Publications.
The Clarendon Press, Oxford University Press (New York); 1991; xiv+266 pp. ISBN 0-19-853568-6 MR1150978 (93h:32021)

\nextitem {[Le]}
\refQ Lelong, P.;
Fonctions plurisousharmoniques{;} mesures de Radon associ\'ees. Applications aux fonctions analytiques. (French);
(Colloque sur les fonctions de plusieurs variables, tenu \'a Bruxelles, 1953),
xxx (ed.), Georges Thone, {(Li\`ege);} Masson \& Cie. (Paris); 1953; 21--40; MR0061682 (15,865a)

\nextitem {[LT]}
\refQ  Levenberg, N., Taylor, B. A.;
Comparison of capacities in $\CC^n$;
(Complex analysis {(Toulouse, 1983)}), xxx (ed.), Lecture Notes in Math., 1094, Springer (Berlin); 1984; 162--172; MR0773108 (86g:32023)

\nextitem {[MM]} \refJ  Micchelli, C. A., Milman, P.;
A formula for Kergin interpolation in $\RR^k$;
\JAT; 29(3); 1980; 294--296; MR0598723 (82h:41008)

\nextitem {[Pl1]}
\refJ Ple\'sniak, W.;
Markov's inequality and the existence of an extension operator for $C\sp \infty$ functions;
\JAT; 61(3); 1990; 106--117; MR1047152 (91h:46065)

\nextitem {[Pl2]}
\refQ Ple\'sniak, W.;
Recent progress in multivariate Markov inequality;
(Approximation theory), xxx (ed.),  Monogr.\ Textbooks Pure Appl.\ Math., 212, Dekker (New York); 1998;  449--464; MR1625243 (99g:41016)

\nextitem {[Ran]}
\refB Range, R. M.;
Holomorphic functions and integral representations in several complex variables;
Graduate Texts in Mathematics, 108. Springer-Verlag (New York);  1986; xx+386 pp. ISBN: 0-387-96259-X, MR0847923 (87i:32001)

\nextitem {[RS]}
\refJ Range, R. M., Siu, Y. T.;
$C^{k}$ approximation by holomorphic functions and $\bar \partial $-closed forms on $C^{k}$ submanifolds of a complex manifold;
\MA; 210; 1974; 105; MR0350068 (50 \#2561)

\nextitem {[Ra]}
\refB Ransford, T.;
Potential theory in the complex plane;
London Mathematical Society Student Texts, 28. Cambridge University Press (Cambridge); 1995; x+232 pp. MR1334766 (96e:31001)

\nextitem {[Sa1]}
\refJ Sadullaev, A.;
Rational approximations and pluripolar sets, (Russian);
\MSNS; 119(161); 1982; 96--118; MR0672412 (84d:32026)

\nextitem {[Sa2]}
\refJ Sadullaev, A.;
A criterion for fast rational approximation in $\CC^n$, (Russian);
\MSNS; 125(167); 1984; 269--279; MR0764481 (86b:32006)

\nextitem {[Sa3]} \refJ Sadullaev, A.; An estimate for polynomials
on analytic sets; {\sl Math.\ USSR-Izv.}; 20; 1983; 493--502;

\nextitem {[Sa4]}
\refJ Sadullaev, A.;
Plurisubharmonic measures and capacities on complex manifolds; {\sl Russian Math Surveys};  36(4); 1981; 61--119; MR0629683 (83c:32026)

\nextitem {[S]}
\refJ Saff, E.;
Regions of meromorphy determined by the degree of best rational approximation;
{\sl Proc. Amer. Math. Soc.}; 29; 1971; 30--38; MR0281930 (43 \#7644)

\nextitem {[SaT]}
\refB Saff, E., Totik, V.;
Logarithmic potentials with external fields;
Springer-Verlag (Berlin); 1997; 505 pp.
ISBN: 3-540-57078-0.

\nextitem {[Sha]}
\refB Shabat, B. V.;
Introduction to Complex Analysis, Part II: Functions of Several Variables;
Amer.\ Math.\ Society Translations of Mathematical Monographs, 110. American Mathematical Society (Providence, RI); 1992; x+371 pp. ISBN: 0-8218-4611-6.

\nextitem {[Sh]}
\refJ Shirokov, N. A.;
The Jackson-Bernstein theorem and strictly convex domains in $\CC^n$,
(Russian);
{\sl Dokl.\ Akad.\ Nauk SSSR}; 276(3); 1984; 1079--1081; MR0753192 (85m:32012)

\nextitem {[Si]}
\refJ Siciak, J.;
Extremal plurisubharmonic functions in $\CC^n$;
{\sl Ann. Polon.\ Math.}; 39; 1981; 175--211; MR0617459 (83e:32018)

\nextitem {[St]}
\refJ Stahl, H.;
The convergence of Pad\'e approximants to functions with branch points;
\JAT; 91; 1997; 139--204; MR1484040 (99a:41017)

\nextitem {[ST]}
\refB Stahl, H., Totik, V.;
General orthogonal polynomials;
Cambridge Univ.\ Press (Cambridge); 1992; MR1163828 (93d:42029)

\nextitem {[St1]}
\refB Stout, E. L.;
The theory of uniform algebras;
Bogden \& Quigley, Inc. (Tarrytown-on-Hudson, N. Y.); 1971; x+509 pp. MR0423083 (54 \#11066)

\nextitem {[St2]}
\refJ Stout, E. L.;
Holomorphic approximation on compact, holomorphically convex, real analytic varieties;
{\sl Proc. Amer. Math. Soc.}; 134 (8); 2006; 2302-2308; MR2213703

\nextitem {[Wei]}
\refQ Weinstock, B.;
Uniform approximation and the Cauchy-Fantappi\'e integral;
(Several complex variables {(Proc.\ Sympos.\ Pure Math., Vol.~XXX.\ Part 2,
Williams Coll., Williamstown, Mass., 1975)}), xxx (ed.), Amer.\ Math.\ Soc. (Providence, R.I.); 1977; 187--191; MR0460720 (57 \#713)

\nextitem {[We]}
\refJ Wermer, J.;
Approximation on a disk;
\MA; 155; 1964; 331; MR0165386 (29 \#2670)

\nextitem {[W]}
\refJ W\'ojcik, A.;
On zeros of polynomials of best approximation to holomorphic and $C^{\infty}$ functions;
{\sl Monatsh.\ Math.};  105(3); 1988; 75--81;  MR0928760 (89c:30103)

\nextitem {[Z1]}
\refJ Zaharjuta, V. P.;
Transfinite diameter, Chebyshev constants and capacity for a compactum in $\CC^n$ (Russian);
\MSNS; 96(138);  1975; 374--389; MR0486623 (58 \#6342)

\nextitem {[Z2]}
\refJ Zaharjuta, V. P.;
Extremal plurisubharmonic functions, orthogonal polynomials, and the Bern\v ste\u\i n-Walsh theorem for functions of several complex variables (Russian);
{\sl Ann. Polon. Math.}; 33 (1-2); 1976/77; 137--148; MR0444988 (56 \#3333)



{

\bigskip\obeylines
Norman Levenberg
Department of Mathematics, Indiana University
Bloomington, IN 47405 USA
{\tt nlevenbe@indiana.edu}

}
\bye